\documentclass[11pt]{amsart}
\usepackage{amsmath}
\usepackage{amsxtra}
\usepackage{amscd}
\usepackage{amsthm}
\usepackage{amsfonts}
\usepackage{amssymb}
\usepackage{eucal}

\newtheorem{thm}{Theorem}
\newtheorem{conj}{Conjecture}
\newtheorem{cor}{Corollary}
\newtheorem{lem}{Lemma}

\theoremstyle{remark}
\newtheorem{remark}{Remark}

\theoremstyle{definition}

\numberwithin{equation}{section}

\newcommand{\bean}{\begin{eqnarray}}
\newcommand{\eean}{\end{eqnarray}}
\newcommand{\be}{\begin{displaymath}}
\newcommand{\ee}{\end{displaymath}}
\newcommand{\bea}{\begin{eqnarray*}}
\newcommand{\eea}{\end{eqnarray*}}

\newcommand{\thmref}[1]{Theorem~\ref{#1}}
\newcommand{\secref}[1]{Section~\ref{#1}}
\newcommand{\lemref}[1]{Lemma~\ref{#1}}

\newcommand{\corref}[1]{Corollary~\ref{#1}}
\newcommand{\remref}[1]{Remark~\ref{#1}}
\newcommand{\conjref}[1]{Conjecture~\ref{#1}}
\newcommand{\nc}{\newcommand}
\nc{\on}{\operatorname}
\nc{\ch}{\mbox{ch}}
\nc{\Z}{{\mathbb Z}}
\nc{\C}{{\mathbb C}}
\nc{\pone}{{\mathbb P}^1}
\nc{\pa}{\partial}
\nc{\F}{{\mathcal F}}
\nc{\arr}{\rightarrow}
\nc{\larr}{\longrightarrow}
\nc{\al}{\alpha}
\nc{\ri}{\rangle}
\nc{\lef}{\langle}
\nc{\W}{{\mathcal W}}
\nc{\la}{\lambda}
\nc{\ep}{\epsilon}
\nc{\su}{\widehat{{\mathfrak s}{\mathfrak l}}_2}
\nc{\sw}{{\mathfrak s}{\mathfrak l}}
\nc{\g}{{\mathfrak g}}
\nc{\h}{{\mathfrak h}}
\nc{\n}{{\mathfrak n}}
\nc{\N}{\widehat{\n}}
\nc{\G}{\widehat{\g}}
\nc{\De}{\Delta}
\nc{\gt}{\widetilde{\g}}
\nc{\Ga}{\Gamma}
\nc{\one}{{\mathbf 1}}
\nc{\z}{{\mathfrak Z}}
\nc{\La}{\Lambda}
\nc{\wt}{\widetilde}
\nc{\wh}{\widehat}
\nc{\cri}{_{\kappa_c}}
\nc{\kk}{h^\vee}
\nc{\sun}{\widehat{\sw}_N}
\nc{\si}{\sigma}
\nc{\el}{\ell}
\nc{\bi}{\bibitem}
\nc{\om}{\omega}
\nc{\ol}{\overline}
\nc{\ds}{\displaystyle}
\nc{\dzz}{\frac{dz}{z}}
\nc{\Res}{\on{Res}}
\nc{\mc}{\mathcal}
\nc{\Cal}{\mathcal}
\nc{\bb}{{\mathfrak b}}
\nc{\ot}{\otimes}
\nc{\R}{{\mc R}}
\nc{\yy}{{\mc Y}}
\nc{\ga}{\gamma}

\nc{\us}{\underset}
\nc{\opl}{\oplus}
\nc{\beq}{\begin{equation}}
\nc{\Fq}{{\mathcal F}}
\nc{\Mq}{{\mathcal M}}
\nc{\Rep}{\on{Rep}}
\nc{\sssec}{\subsubsection}
\nc{\ssec}{\subsection}
\nc{\lan}{\langle}
\nc{\ran}{\rangle}

\nc{\D}{\mathcal D}
\nc{\Vect}{\on{Vect}}
\nc{\ghat}{\G}
\nc{\T}{\mc T}
\nc{\Tloc}{\T^\g_{\on{loc}}}
\nc{\vac}{|0\ran}
\nc{\Wick}{{\mb :}}
\nc{\mb}{\mathbf}
\nc{\delz}{\partial_z}
\nc{\K}{{\cali K}}
\nc{\cali}{\mathcal}
\nc{\li}{\mathfrak l}
\nc{\lt}{\widetilde{\li}}
\nc{\astar}{a^*}
\nc{\cA}{{\mc A}}
\nc{\ka}{\kappa}

\nc{\OO}{{\mc O}}
\nc{\AutO}{\on{Aut}\OO}
\nc{\DerO}{\on{Der}\OO}
\nc{\DerpO}{\on{Der}_+\OO}
\nc{\Au}{{\mc A}ut}
\nc{\mf}{\mathfrak}
\nc{\V}{{\mathbb V}}
\nc{\hh}{\wh{\h}}

\nc{\pp}{{\mathfrak p}}
\nc{\mm}{{\mathfrak m}}
\nc{\rr}{{\mathfrak r}}
\nc{\ket}{\rangle}
\nc{\zz}{{\mathfrak z}}
\nc{\gr}{\on{gr}}
\nc{\Spe}{\on{Spec}}
\nc{\rv}{\crho}
\nc{\can}{\on{can}}
\nc{\CC}{\on{Op}_G(D))}
\nc{\Op}{\on{Op}_G(D)}
\nc{\MOp}{\on{MOp}_G(D)}
\nc{\Db}{{\mathbb D}}
\nc{\ww}{w}

\nc{\af}{{\mathbb A}^1}
\nc{\bs}{\backslash}
\nc{\laa}{(\la_i)}
\nc{\zn}{(z_i)}

\nc{\cla}{\check{\la}}
\nc{\cmu}{\check{\mu}}
\nc{\crho}{\check{\rho}}
\nc{\chal}{\check{\al}}
\nc{\cc}{{\mathfrak c}}

\nc{\M}{{\mathbb M}}

\nc{\ZZ}{{\mc Z}}

\nc{\UU}{{\mathbb U}}

\nc{\Conn}{\on{Conn}(\Omega^{\crho})}
\nc{\Con}{\on{Conn}(\Omega^{-\rho})}
\nc{\Co}{\on{Conn}(\Omega^{\rho})}

\nc{\ppart}{(\!(t)\!)}
\nc{\pparl}{(\!(\la)\!)}
\nc{\zpart}{(\!(z)\!)}
\nc{\ppzi}{(\!(t-z_i)\!)}
\nc{\ppinf}{(\!(t^{-1})\!)}
\nc{\Ind}{\on{Ind}}
\nc{\I}{{\mathbb I}}

\nc{\Bun}{\on{Bun}}

\newcommand {\IC}{\mathbb{C}}

\newcommand {\A}{{\mathcal A}}

\newcommand {\hatg}{\wh{\g}}

\newcommand {\ind}{\operatorname{Ind}}

\nc{\gtil}{\wt{\g}}
\nc{\ntil}{\wt{\n}}
\nc{\htil}{\wt{\h}}
\nc{\gbar}{\ol{\g}}
\nc{\nbar}{\ol{\n}}
\nc{\bbar}{\ol{\bb}}
\nc{\lhat}{\wh{\mf l}}

\begin{document}

\title{Quantization of soliton systems and Langlands duality}

\author{Boris Feigin}

\address{Landau Institute for Theoretical Physics, Kosygina St 2,
Moscow 117940, Russia}

\author{Edward Frenkel}\thanks{B.F. was supported by the grants RFBR
05-01-01007, RFBR 05-01-02934 and NSh-6358.2006.2. E.F. was
supported by the DARPA Program ``Focus Areas in Theoretical
Mathematics'' and the NSF grant DMS-0303529.}

\address{Department of Mathematics, University of California,
  Berkeley, CA 94720, USA}

\date{May 2007. Revised October 2009.}

\begin{abstract}

We consider the problem of quantization of classical soliton
integrable systems, such as the KdV hierarchy, in the framework of a
general formalism of Gaudin models associated to affine Kac--Moody
algebras. Our experience with the Gaudin models associated to
finite-dimensional simple Lie algebras suggests that the common
eigenvalues of the mutually commuting quantum Hamiltonians in a model
associated to an affine algebra $\ghat$ should be encoded by affine
opers associated to the Langlands dual affine algebra $^L \ghat$. This
leads us to some concrete predictions for the spectra of the quantum
Hamiltonians of the soliton systems. In particular, for the KdV system
the corresponding affine opers may be expressed as Schr\"odinger
operators with spectral parameter, and our predictions in
this case match those recently made by Bazhanov, Lukyanov and
Zamolodchikov. This suggests that this and other recently found
examples of the correspondence between quantum integrals of motion and
differential operators may be viewed as special cases of the
Langlands duality.

\end{abstract}

\maketitle

\section{Introduction}

Soliton equations, such as the celebrated KdV hierarchy, are
infinite-dimensional classical integrable systems possessing infinite
collections of Poisson commuting Hamiltonians. It is natural to try to
{\em quantize} these systems. This means, in particular, constructing
mutually commuting quantum Hamiltonians in the appropriate
associative algebra (such as the completed enveloping algebra of the
Virasoro algebra in the case of the KdV hierarchy), whose symbols are
the classical Hamiltonians, and finding the spectra of these quantum
Hamiltonians on representations of this algebra. This problem has
acquired additional importance after the work of A. Zamolodchikov
\cite{Zam}, in which he showed that these quantum Hamiltonians may be
identified with integrals of motion of certain deformations of
conformal field theories.

The existence of local commuting quantum Hamiltonians has been
established in \cite{FF:toda,FF:laws} in the case of the generalized
KdV hierarchies and in \cite{FF:ih} in the case of the generalized
AKNS hierarchies. The corresponding non-local Hamiltonians have been
constructed in \cite{BLZ1,BHK} for the KdV hierarchies (we will
discuss below an extension of this construction to the AKNS
hierarchies). A study of the problem of diagonalization of the quantum
KdV Hamiltonians was initiated in the important series of works
\cite{BLZ1}--\cite{BLZ3} by V. Bazhanov, S. Lukyanov and
A. Zamolodchikov.

These works culminated in a remarkable conjecture made in \cite{BLZ}
(building on the earlier works \cite{DT,BLZ4}), which related the
spectra of the {\em quantum KdV Hamiltonians} to certain
one-dimensional {\em Schr\"odinger operators}. This connection appears
to be rather mysterious. In this paper we suggest that it may be
viewed as a special case of a broader picture which fits in the
general framework of integrable systems of Gaudin type and the
Langlands correspondence.

The first step in our approach is the realization that the quantum
soliton systems may be viewed as special cases of {\em Gaudin
models} associated to affine Kac--Moody algebras. At the classical
level, this is quite clear. Indeed, it has been known since
\cite{RS,DS} that the algebra of classical Poisson commuting
Hamiltonians of the soliton hierarchies may be obtained using the
``shift of argument'' method, applied to an affine Kac--Moody
algebra. Its finite-dimensional counterpart is the Poisson commuting
``shift of argument'' subalgebra constructed in \cite{MF}. As
explained in \cite{FFT}, the corresponding Hamiltonian system is the
simplest example of a classical Gaudin model with irregular
singularities (which is in turn a special case of the Hitchin systems
\cite{Hitch,Beauville}). Recently, a quantization of this algebra, as
well as its generalizations corresponding to other Gaudin models, has
been constructed in \cite{Ryb,FFT}. We expect that the soliton
integrable systems may be quantized along the same lines, that is, in
the framework of the quantization of Gaudin models.

The second step in our approach is the link between the spectra of the
quantum {\em Gaudin Hamiltonians} and differential operators known as
{\em opers}, discovered in our earlier works \cite{FFR,F:faro,FFT}. We
have shown there that the spectra of the Gaudin Hamiltonians
associated to a simple finite-dimensional Lie algebra $\g$ are encoded
by $^L \g$-opers (where $^L \g$ is the Langlands dual Lie algebra of
$\g$) on $\pone$ with finitely many singular points.

The key point of our analysis \cite{FFR,F:faro,FFT} of the Gaudin
models associated to a simple Lie algebra $\g$ is the identification
of the commutative algebra of quantum Gaudin Hamiltonians with a
quotient of the center $Z(\ghat)$ of the completed enveloping algebra
of the affine Kac--Moody algebra $\ghat$ at the critical level. Now we
consider the Gaudin systems associated to an affine algebra
$\ghat$. Therefore it is natural to assume that the corresponding
Hamiltonians come from the center $Z(\widehat{\ghat})$ of the
completed enveloping algebra of the {\em double affine algebra}
$\widehat{\ghat}$. At present, we do not have a precise definition of
this enveloping algebra or its center. However, we postulate that it
exhibits the same salient features as its affine counterpart
$Z(\ghat)$.

According to a theorem of \cite{FF:gd,F:wak}, $Z(\ghat)$ is isomorphic
to the algebra of functions on the space $\on{Op}_{^L \g}(D^\times)$
of $^L \g$-opers on the (formal) punctured disc $D^\times$. Hence we
expect that $Z(\widehat{\ghat})$ is related to the space of {\em
affine opers} $\on{Op}_{^L \ghat}(D^\times)$ associated to $^L \ghat$,
defined in \cite{F:flag} (see also \cite{BF}). Here $^L \ghat$ is the
Langlands dual Lie algebra to the affine algebra $\ghat$, that is, its
Cartan matrix is the transpose of that of $\ghat$ (so the dual of an
untwisted affine algebra might be a twisted affine algebra). Testing
this conjectural relation was in fact one of the main motivations for
the present work. Since the structure of the center
$Z(\widehat{\ghat})$ seems out of reach at the moment, we turn the
tables and try to use the available information about the affine
Gaudin models and the quantum soliton systems to gain insights into
the structure of $Z(\widehat{\ghat})$.

The upshot of all this is that the spectra of the quantum $\ghat$-KdV
Hamiltonians should be represented by the affine $^L \ghat$-opers on
$\pone$ with prescribed singularities. In the case when $\g=\sw_2$
these opers may be represented by Schr\"odinger operators with a
spectral parameter. Thus, we suggest that the connection between the
spectra of quantum KdV Hamiltonians and Schr\"odinger operators
conjectured in \cite{BLZ} comes about as follows:

\begin{itemize}

\item The quantum KdV system is interpreted as a generalized Gaudin
  model associated to the affine Kac--Moody algebra $\su$ (or $\ghat$
  for the more general quantum $\ghat$-KdV systems).

\item The spectra of the Gaudin Hamiltonians associated to an
  affine Kac--Moody algebra $\ghat$ are encoded by the affine $^L
  \ghat$-opers.

\item Since the quantum KdV Hamiltonians may be interpreted as the
  Gaudin Hamiltonians associated to $\su$, their spectra should be
  encoded by $\su$-opers, which are nothing but Schr\"odinger
  operators with spectral parameter: \linebreak $\pa_t^2 - v(t) -
  \la$. More generally, we expect that the spectra of the quantum
  $\ghat$-KdV Hamiltonians are encoded by affine $^L
  \ghat$-opers. These opers should satisfy an important {\em no
  monodromy} property.
\end{itemize}

\medskip

Thus, schematically, the correspondence we propose looks as follows:

$$
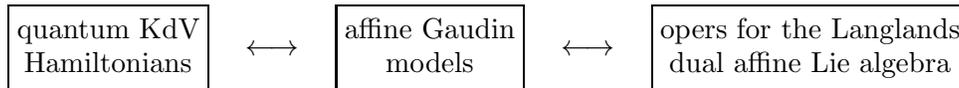
\boxed{\begin{matrix} \text{quantum KdV} \\
\text{Hamiltonians} \end{matrix}} \quad \longleftrightarrow \quad
\boxed{\begin{matrix} \text{affine Gaudin} \\
    \text{models} \end{matrix}} \quad \longleftrightarrow \quad
\boxed{\begin{matrix} \text{opers for the Langlands}  \\
    \text{dual affine Lie algebra}  \end{matrix}}
$$

\medskip

The appearance of the Langlands dual Lie algebra suggests that this
connection (and perhaps the more general ODE/IM correspondence of
Dorey, e.a., \cite{DDT,Dorey1,Dorey2}) should be viewed as a special
case of Langlands duality.

\medskip

In order to test this proposal, we compare our $\su$-opers with the
Schr\"odinger operators obtained in \cite{DT,BLZ4,BLZ}. We find that
they match perfectly, after a simple change of variable. Our formulas
easily generalize to an arbitrary affine Lie algebra $\ghat$, and
again we find that they are consistent with the known formulas for the
differential operators associated to the spectra of the quantum
Hamiltonians in \cite{DDT,BHK,Dorey1,Dorey2}.

In \cite{DT,BLZ4,DT1,DDT,BHK,BLZ} the connection between the spectra
of the quantum KdV Hamiltonians and differential operators is
explained by the fact that certain generating functions of the former
(the so-called $Q$-operators) satisfy the same functional relations as
connection coefficients (or Stokes multipliers) of the latter. These
functional relations are the so-called $T-Q$ relations generalizing
Baxter's famous relation. In \cite{Dorey1,Dorey2} (and earlier works
referenced there) the functional relations of this type were taken as
the basis of an ``ODE/IM correspondence'' for more general solvable
models.

In contrast, in this paper we do not attempt to explain the connection
between the spectra of quantum Hamiltonians and affine opers. Rather,
our goal is to emphasize that the connection between them fits in the
general scheme of quantum Gaudin models and is in fact a manifestation
of the Langlands duality. This is based on our analysis
\cite{FFR,F:faro,FFT} in the finite-dimensional case and on the fact
that in the affine case our predictions match those obtained in
\cite{DT,BLZ4,DDT,BHK,BLZ}. This already has a number of consequences:

\begin{enumerate}

\item the appearance of the Langlands dual affine Lie algebra
  (something that would be difficult to see based on the most studied
  case of $\wh\sw_n$, which is self-dual);

\item a uniform way of expressing the differential operators encoding
  the spectra of the quantum Hamiltonians as opers associated to
  arbitrary affine Lie algebras (again, this is not obvious
  from the previously known examples);

\item a Bethe Ansatz method for constructing eigenvectors of the
  quantum Hamiltonians.

\end{enumerate}

In particular, we conjecture explicit Bethe Ansatz formulas
(generalizing the Bethe Ansatz in the finite-dimensional Gaudin
models) for eigenvectors of the quantum KdV Hamiltonians in the
irreducible representations of the Virasoro algebra from the unitary
minimal models.

We expect that the appearance of the $T-Q$ relations of
\cite{DT,BLZ4,DT1,DDT,BHK,BLZ} may also be explained in the framework
of the formalism developed in this paper. More specifically, we
believe that it comes from the connection between quantum integrals of
motion and the center of a quantum affine algebra at the critical
level. We plan to address these issues in a future work.

\medskip

The paper is organized as follows. In \secref{shift fd} we recall the
results of \cite{Ryb,FFT} on the quantization of the shift of argument
Gaudin algebras in the finite-dimensional case and the identification
of their spectra with $^L \g$-opers. We then consider in \secref{shift
aff} the affine analogue of this construction which corresponds to the
quantization of the Hamiltonians of the soliton hierarchies such as
the $\g$-AKNS. We discuss in detail the classical Hamiltonians
obtained by expansion of the monodromy matrix with respect to the
spectral parameter (following \cite{RS}) and the quantum analogues of
both local and non-local Hamiltonians. We then construct eigenvectors
of these Hamiltonians using an affine analogue of the Bethe Ansatz
method. In order to describe the spectra of the quantum Hamiltonians,
we introduce in \secref{affine opers} the affine opers, following
\cite{F:flag}. We conjecture that the affine opers that encode the
spectra of the quantum Hamiltonians on representations of $\ghat$ are
$^L \ghat$-opers on $\pone$ satisfying certain analytic properties (in
particular, to integrable representations of $\ghat$ correspond $^L
\ghat$-opers without monodromy).

Next, in \secref{kdv}, we apply the same analysis to the generalized
KdV hierarchies. We state our conjectures describing the spectra of
the quantum KdV Hamiltonians in terms of the affine opers
corresponding to the Langlands dual affine Lie algebra $^L \ghat$. We
write down explicit formulas for these opers in the case when
$\ghat=\su$ and show that after a simple change of variables they give
rise to the Schr\"odinger operators constructed in \cite{BLZ}. We then
discuss the generalization of these formulas to the case of
$\wh\sw_n$ and to general affine Lie algebras.

Motivated by this study of the ``shift of argument'' soliton systems,
corresponding to the simplest affine Gaudin models with irregular
singularities, we introduce in \secref{reg sing} the affine analogues
of the Gaudin models with regular singularities. We show that in this
case we also have interesting classical integrable systems, which may
be quantized by analogy with the finite-dimensional case. Again, we
conjecture that the spectra of the corresponding quantum Hamiltonians
are described by the affine $^L \ghat$-opers on $\pone$ with regular
singularities at finitely many points. We also construct the
eigenvectors of these Hamiltonians using Bethe Ansatz. We show that a
particular case of this model is closely related to the GKO coset
construction \cite{GKO}. We conjecture that in this case our Bethe
vectors are eigenvectors of all quantum KdV Hamiltonians on
irreducible modules over the Virasoro algebra from the unitary minimal
models.

We summarize our results and conjectures in \secref{conclusion}.

\bigskip

\noindent{\bf Acknowledgments.} This project was conceived and to a
large extent carried out during our visits to Kyoto University in
2004--2007. We thank our hosts M. Kashiwara, T. Miwa and H. Nakajima,
as well as K. Hasegawa, for their hospitality.

We also thank V. Bazhanov for useful discussions.

\section{The shift of argument Gaudin models in the finite-dimensional
  case}    \label{shift fd}

We start by summarizing what is known about the simplest Gaudin models
with irregular singularities, following \cite{Ryb,FFT}. The
corresponding classical Hamiltonian system may be viewed as an
example of a Hitchin system. The Poisson commutative algebra of
classical Hamiltonians is the so-called ``shift of argument''
subalgebra $\ol{\mc A}_\chi$ of the symmetric algebra $S(\g)$ of a
simple Lie algebra $\g$, introduced by A.S. Mishchenko and
A.T. Fomenko \cite{MF} (see also \cite{Man}). The corresponding
algebra of quantum Hamiltonians is a commutative subalgebra of the
universal enveloping algebra $U(\g)$. It has been constructed in
\cite{Ryb,FFT} using the center of the completed universal enveloping
algebra of the affine Kac--Moody algebra $\ghat$ at the critical
level. This allowed us to identify in \cite{FFT} the spectrum of this
algebra with the space of $^L \g$-opers on $\pone$ with regular
singularity at $0 \in \pone$ and irregular singularity at $\infty \in
\pone$ (here $^L \g$ is the Langlands dual Lie algebra to $\g$).

We also consider here analogous commutative subalgebras in the
Hamiltonian reductions of $U(\g)$ known as finite ${\mc W}$-algebras.

In the following sections we will use this construction as a prototype
for the description of the quantum ``shift of argument'' subalgebra
when the finite-dimensional Lie algebra $\g$ is replaced by its
affinization $\ghat$. It is well-known (see \cite{RS,DS}) that the
corresponding classical systems give rise to hierarchies of soliton
equations such as $\g$-AKNS and $\ghat$-KdV. In other words, the
classical shift of argument subalgebra is generated by the classical
Hamiltonians of these hierarchies. Therefore the corresponding quantum
Hamiltonians generate the quantum shift of argument subalgebra, and
this suggests that their spectra are encoded by affine $^L
\ghat$-opers. We will find that in the case of the quantum KdV
Hamiltonians this matches the predictions of \cite{BLZ}.

\subsection{Classical system}    \label{class fd}

The Gaudin models with irregular singularities for finite-dimensional
simple Lie algebras have been introduced and studied in
\cite{FFT}. Here we will focus on the particular case related to the
so-called shift of argument subalgebra $\ol{\mc A}_\chi \subset
S(\g)$ and its quantization ${\mc A}_\chi \subset U(\g)$, as it is
most closely related to the structures emerging in soliton equations
(such as the KdV and $n$-wave hierarchies) and their quantization,
once we switch from $\g$ to the affine algebra $\ghat$. But other
Gaudin models with irregular singularities may also have interesting
analogues in the affine case.

Our starting point is a classical Hitchin system (see
\cite{Hitch,Beauville}) in which the phase space is the space of the
following $\g^*$-valued one-forms on $\pone$:
\begin{equation}    \label{eta1}
\eta = \left( \frac{A}{z} + \chi\right) dz, \qquad A \in \g^*,
\end{equation}
where $\chi \in \g^*$ is fixed. This space is therefore isomorphic to
$\g^*$. A one-form like this has a pole of order $2$ at the point
$\infty$ with the leading coefficient $-\chi$, and it also has a pole
of oder $1$ at the point $0$.

We define the classical shift of argument algebra $\ol{\A}_\chi$ as
the subalgebra of $S(\g) = \on{Fun} \g^*$ generated by the
coefficients of the polynomials in $z^{-1}$ obtained by evaluating the
invariant polynomials $P \in (\on{Fun} \g^*)^G,
i=1,\ldots,\ell$, on $\eta$. In other words, $\ol{\A}_\chi$ is the
subalgebra of $\on{Fun} \g^*$ generated by the iterated directional
derivatives $D^{i}_\chi P$ of invariant polynomials $P\in(\on{Fun}
\g^*)^G$ in the direction $\chi$, where
\begin{equation}    \label{dir der}
D_\chi P(y)=\left.\frac{d}{du}\right|_{u=0}P(y+u\chi).
\end{equation}
Equivalently, $\ol{\A}_\chi$ is the subalgebra of $S\g$ generated by
the shifted polynomials
\begin{equation}\label{eq:shift}
P_{u\chi}(y)=P(y+u\chi)
\end{equation}
where $P$ varies in $(\on{Fun}\g^*)^G$ and $u\in\C$, hence its name.

The algebra $\ol{\A}_\chi$ is Poisson commutative with respect to the
Kirillov--Kostant Poisson structure on $\g^*$. This follows from the
general results on the Hitchin systems with singularities.

The algebra $\ol{\A}_\chi$ was first introduced in \cite{MF}. It has
been shown in \cite{MF,Ko2,FFT} that for regular $\chi$ the
algebra $\ol{\mc A}_\chi$ is a maximal Poisson commutative subalgebra
of $S(\g)$ which is a free polynomial algebra in $\dim\bb$ generators
\begin{equation}    \label{iterated}
D^{n_i}_\chi \ol{P}_i, \qquad i=1,\ldots,\ell; \; n_i=0,\ldots,d_i,
\end{equation}
where $\ol{P}_i$ is a generator of $S(\g)^{\g}$ of degree $d_i+1$ and
$D_\chi$ is the derivative in the direction of $\chi$ given by formula
\eqref{dir der}.\footnote{We recall that $a \in \g$ is called
{\em regular} if its centralizer in $\g$ has the smallest possible
dimension; namely, $\ell$, the rank of $\g$. We identify $\g^*$ with
$\g$ by using a non-degenerate invariant inner product.}

Note that for any $\chi$ the algebra $\ol{\mc A}_\chi$ belongs to the
centralizer $S(\g)^{\g_\chi}$ of $\g_\chi$ in $S(\g)$, where $\g_\chi
\subset \g$ is the stabilizer of $\chi$. Let $\psi \in \g_\chi^*$ be a
character of $\g_\chi$, or, in other words, a one-point coadjoint
orbit in $\g^*_\chi$. Set
$$
S^\psi_\chi(\g) = (S(\g)/(\g_\chi - \psi(\g_\chi)))^{\g_\chi},
$$
where $(\g_\chi - \psi(\g_\chi))$ is the ideal in $S(\g)$ generated by
the elements $A-\psi(A), A \in \g_\chi$. The algebra $S^\psi_\chi(\g)$
is a Poisson algebra obtained by Hamiltonian reduction of $\g^*$ with
respect to the adjoint action of $G_\chi$ and the one-point orbit
$\psi \in \g^*_\chi$. Let $\ol{\mc A}^\psi_\chi$ be the projection of
$\ol{\mc A}_\chi$ onto $S^\psi_\chi(\g)$. This is a Poisson
commutative subalgebra of $S^\psi_\chi(\g)$.

For regular $\chi$ the Lie algebra $\g_\chi$ is commutative, and its
symmetric algebra $S(\g_\chi)$ is contained in ${\mc A}_\chi$. In
fact, a basis of $\g_\chi$ may be obtained by taking the maximal
non-constant directional derivatives $D^{d_i}_\chi \ol{P}_i,
i=1,\ldots,\ell$, among the generators \eqref{iterated}. Therefore for
any $\psi \in \g^*_\chi$ the algebra $\ol{\mc A}^\psi_\chi$ is a free
polynomial algebra generated by the (non-zero) images of
$$
D^{n_i}_\chi \ol{P}_i, \qquad i=1,\ldots,\ell; \; n_i=0,\ldots,d_i-1,
$$
in $S^\psi_\chi(\g)$.

For nilpotent elements $\chi$ the algebras $S^\psi_\chi(\g)$ are known
as the (classical) {\em finite ${\mc W}$-algebras}, and the above
construction produces Poisson commutative subalgebras in these ${\mc
W}$-algebras.

For example, suppose that $\chi$ corresponds to the maximal root
generator $f_\theta \in \n_- \subset \g$ under an identification $\g
\simeq \g^*$ given by a non-degenerate invariant inner product on
$\g$. Then $\g_\chi = \n_-$. Let $\psi$ be a principal character on
$\n_-$ taking non-zero values on the simple root generators $f_i$. We
have $S^\psi_\chi(\g) = S(\g)^{\g}$ in this case, and therefore it is
clear that $\ol{\mc A}^\psi_\chi = S(\g)^{\g}$ as well. However, the
affine analogue of this construction is non-trivial: the analogue of
$S^\psi_\chi(\g)$ is the classical ${\mc W}$-algebra obtained by the
Drinfeld--Sokolov reduction of $\ghat^*_1$, and the analogue of
$\ol{\mc A}^\psi_\chi$ is its Poisson commutative subalgebra generated
by the (local and non-local) Hamiltonians of the $\ghat$-KdV
hierarchy (see \secref{kdv}).

\subsection{Quantization}    \label{quant fd}

Let us now discuss the quantization of the algebra $\ol{\mc A}_\chi$
in the finite-dimensional case. Here by quantization we understand a
commutative subalgebra ${\mc A}_\chi$ of $U(\g)$ such that the
associated graded $\on{gr} {\mc A}_\chi$ of ${\mc A}_\chi$ (with
respect to the standard filtration on $U(\g)$) is a Poisson
commutative subalgebra of $\on{gr} U(\g) = S(\g)$ which contains
$\ol{\mc A}_\chi$. Such a quantization has been constructed in
\cite{Ryb,FFT} in a uniform way for an arbitrary simple Lie algebra
$\g$ (following the procedure that had been used for constructing the
algebra of quantum Hamiltonians in the Gaudin models with regular
singularities in \cite{FFR,F:faro}, see \secref{RS fd}
below).\footnote{For classical Lie algebras the quantization of
$\ol{\mc A}_\chi$ had been constructed previously in a different way
in \cite{NO,Tar,CT}.}

Let $\ghat$ be the affine Kac--Moody algebra associated to $\g$ (see
\secref{aff km}), and
\begin{equation}    \label{vacuum}
\V_{\on{crit}} = \ind_{\g[[t]] \oplus \C {\mb 1}}^{\hatg}\IC_{-h^\vee}
\end{equation}
the vacuum $\ghat$-module of critical level $k=-h^\vee$, where
$h^\vee$ is the dual Coxeter number of $\g$ (in this formula $\g[[t]]$
acts by $0$ on the one-dimensional module $\IC_{-h^\vee}$, and ${\mb
1}$ acts by multiplication by $-h^\vee$). The quantum shift of
argument algebra ${\mc A}_\chi$ was obtained in \cite{Ryb,FFT} as a
quotient of the algebra $\on{End}_{\ghat} \V_{\on{crit}}$. According
to a theorem of \cite{FF:gd,F:wak}, the latter algebra is isomorphic
to the algebra of functions on the space $\on{Op}_{^L \g}(D)$ of $^L
\g$-opers on the (formal) disc. Here $^L \g$ is the Langlands dual Lie
algebra to $\g$, whose Cartan matrix is the transpose of the Cartan
matrix of $\g$. The notion of oper was introduced by Beilinson and
Drinfeld \cite{BD,BD:opers} following Drinfeld and Sokolov
\cite{DS}. We refer the reader to the definition given in
\cite{F:flag,F:faro}, which is well-suited for our present goals. (We
will recall this definition below in \secref{affine opers}.)

Thus, we have an isomorphism
\begin{equation}    \label{isom}
\on{End}_{\ghat} \V_{\on{crit}} \simeq \on{Fun} \on{Op}_{^L \g}(D).
\end{equation}
We note that $\on{End}_{\ghat} \V_{\on{crit}}$ is the quotient of the
center $Z(\ghat)$ of the completed enveloping algebra of $\ghat$ at
the critical level, which is isomorphic to the algebra of functions on
the space $\on{Op}_{^L \g}(D^\times)$ of opers on the (formal)
punctured disc (see, e.g., \cite{F:book}, Sect. 4.3).

Using the isomorphism \eqref{isom} and general results on coinvariants
from \cite{FB}, we have obtained the following description of ${\mc
A}_\chi$ in \cite{FFT}:

\begin{thm}
For regular $\chi \in \g^*$ the algebra ${\mc A}_\chi$ is canonically
isomorphic to the algebra of functions on the space $\on{Op}_{^L
\g}(\pone)_{\pi(\chi)}$ of $^L \g$-opers on $\pone$ with regular
singularity at the point $0$ and with irregular singularity of order
$2$ at $\infty$, with the $2$-residue $\pi(\chi)$:
\begin{equation}    \label{ident Achi}
{\mc A}_\chi \simeq \on{Op}_{^L \g}(\pone)_{\pi(\chi)}.
\end{equation}
\end{thm}

For example, in the case when $\g=\sw_2$ and $\chi \in \h^*$, the space
$\on{Op}_{^L \g}(\pone)_{\pi(\chi)}$ consists of second order
differential operators of the form
\begin{equation}    \label{fuchsian1}
\pa_z^2 - \chi^2 + \frac{v_{-1}}{z} + \frac{v_{-2}}{z^2},
\qquad v_{-1}, v_{-2} \in \C.
\end{equation}
Note that they have a pole of order $4$ at $z=\infty$ with the leading
term $-\chi^2$. For $\g=\sw_n$ the role of operators \eqref{fuchsian1}
is played by the following $n$th order differential operators:
$$
(-\pa_z)^n - v_1(z) (-\pa_z)^{n-2} - \ldots + v_{n-2}(z) \pa_z -
v_{n-1}(z),
$$
where $v_j(z)$ has the form
$$
v_j(z) = P_j(\chi) + \sum_{k=1}^{j+1} v_{j,k} z^{-k},
$$
and $P_j$ is a particular symmetric polynomial of degree $j+1$.

For an explicit description of $\on{Op}_{^L \g}(\pone)_{\pi(\chi)}$
for general simple Lie algebras, see \cite{FFT}.

\medskip

It has been proved in \cite{Ryb,FFT} that $\on{gr} {\mc
A}_\chi$ contains $\ol{\mc A}_\chi$, with the equality for regular
$\chi$. For general $\chi \in \g^*$ this equality is a
conjecture (see Conjecture 1 of \cite{FFT}).

\medskip

Let us summarize: the quantization ${\mc A}_\chi$ of the classical
Poisson commutative algebra $\ol{\mc A}_\chi$ is a commutative
subalgebra in a quantum (non-commutative) algebra $U(\g)$, but it
turns out to be isomorphic to another classical commutative algebra,
namely the algebra of functions on the space $\on{Op}_{^L
\g}(\pone)_{\pi(\chi)}$. A conceptual understanding of this
isomorphism is obtained by realizing ${\mc A}_\chi$ as a quotient of
the center of an enveloping algebra of the affine Kac--Moody algebra
$\ghat$ and using the identification of the latter with the algebra of
functions on opers. The non-triviality of this isomorphism is
underscored by the fact that the Lie algebra $\g$ gets replaced by its
Langlands dual Lie algebra $^L \g$ under this isomorphism. In fact,
this isomorphism may be viewed as an example of the geometric
Langlands correspondence, as explained in \cite{F:icmp,F:rev,FFT}.

\medskip

The Poisson commutative subalgebras $\ol{\mc A}^\psi_\chi
\subset S^\psi_\chi(\g)$ may also be quantized. Namely, by
construction, ${\mc A}_\chi$ is contained in $U(\g)^{\g_\chi}$. For a
character $\psi: \g_\chi \to \C$, let ${\mc A}^\psi_\chi$ be the
projection of ${\mc A}_\chi$ onto
$$
U^\psi_\chi(\g) = (U(\g)/(\g_\chi-\psi(\g_\chi)))^{\g_\chi},
$$
where $(\g_\chi-\psi(\g_\chi))$ is the left ideal in $U(\g)$ generated
by $A-\psi(A), A \in \g_\chi$. Thus, $U^\psi_\chi(\g)$ is the result
of quantum Hamiltonian (or BRST) reduction of $U(\g)$ with respect to
$\g_\chi$ and its character $\psi$. In particular, if $\chi$ is
nilpotent, $U^\psi_\chi(\g)$ is a quantum finite ${\mc W}$-algebra,
and so ${\mc A}^\psi_\chi$ gives rise to its commutative
subalgebra. For example, if $\chi$ corresponds to the maximal root
generator $f_\theta \in \n_- \subset \g$ and $\psi$ is a principal
character on $\n_-$ (as in \secref{class fd}), then ${\mc A}^\psi_\chi
= U^\psi_\chi(\g) = Z(\g)$, the center of $U(\g)$.

\subsection{The spectra of the quantum Hamiltonians}    \label{spectra}

Suppose that $\chi \in \g^* \simeq \g$ is regular semi-simple. Without
loss of generality we may assume that $\chi$ is in a fixed Cartan
subalgebra $\h \subset \g$. The algebra ${\mc A}_\chi$ then contains
$\h$. It is natural to try to describe the spectrum of ${\mc A}_\chi$
on a given $\g$-module $M$. Note that since $\h \subset {\mc A}_\chi$,
we obtain that the action of ${\mc A}_\chi$ will preserve the weight
decomposition of $M$. The isomorphism \eqref{ident Achi} implies that
the common eigenvalues of ${\mc A}_\chi$ on $M$ are encoded by $^L
\ghat$-opers on $\pone$ with regular singularity at the point $0$ and
with singularity of order $2$ at $\infty$, with the $2$-residue
$\pi(\chi)$ (which is, by definition, an element of $\h^*/W$, see,
e.g., \cite{FFT}). In \cite{FFT}, Theorem 5.10, we have obtained the
following more precise description of the common spectra of ${\mc
A}_\chi$ on the Verma modules and irreducible finite-dimensional
representations of $\g$.

\begin{thm}    \label{spectra fd} \

{\em (1)} If $M=M_\nu$, the Verma module with highest weight $\nu \in
\h^*$, then the common eigenvalues of ${\mc A}_\chi$ on $M$ are
encoded by $^L \g$-opers on $\pone$ with singularity of order $2$
at $\infty$, with the $2$-residue $\pi(\chi)$ and with regular
singularity at the point $0$ and the residue (which is, by definition,
an element of $\h^*/W$, see, e.g., \cite{FFT}) equal to
$\varpi(-\nu-\rho)$, where $\varpi$ is the projection $\h^* \to
\h^*/W$.

{\em (2)} If $M=V_\nu$ is an irreducible finite-dimensional
representation of $\g$, then the common eigenvalues of ${\mc A}_\chi$
on $M$ are encoded by the $^L \g$-opers in part (1) satisfying the
additional condition that they have trivial monodromy.
\end{thm}

According to Conjecture 2 of \cite{FFT}, {\em all} opers described in
part (2) of this theorem are realized in the spectrum of ${\mc
A}_\chi$ on $V_\nu$.\footnote{This conjecture has been proved by
L. Rybnikov and the authors in {\em Opers with irregular singularity
and spectra of the shift of argument subalgebra}, Preprint
arXiv:0712.1183.}

This is analogous to the description of the spectra of the
Hamiltonians in the Gaudin models with regular singularities from
\cite{FFR,F:faro} which we will recall in \secref{RS fd}.

In \cite{FFT} we have constructed eigenvectors of ${\mc A}_\chi$ by
Bethe Ansatz. Suppose that $M=M_\nu$. Then the Bethe eigenvectors have
the form
\begin{multline}    \label{bethe for one mod}
\phi(w_1^{i_1},\ldots,w_m^{i_m}) = \\ \sum_{\sigma \in S_m}
\frac{f_{i_{\sigma(1)}} f_{i_{\sigma(2)}} \ldots
f_{i_{\sigma(m)}}}{(w_{{\sigma(1)}} -
w_{{\sigma(2)}})(w_{{\sigma(2)}} - w_{{\sigma(3)}}) \ldots
(w_{{\sigma(m-1)}}-w_{{\sigma(m)}}) w_{{\sigma(m)}}} v_{\nu},
\end{multline}
where $v_\nu$ is the highest weight vector of $M_\nu$, and the sum is
over all permutations on $m$ letters. This vector has the weight
$$
\nu - \sum_{j=1}^m \al_{i_j}.
$$

The following result is proved in \cite{FFT}, Sect. 6.6 (for the
definition of Miura transformation, see \cite{FFT} and \secref{affine
opers} below).

\begin{thm}
The Bethe vector $\phi(w_1^{i_1},\ldots,w_m^{i_m})$ is an eigenvector
of ${\mc A}_\chi$ if the complex numbers $w_j$ satisfy the Bethe
Ansatz equations\footnote{Note that our $\chi$ here corresponds to
$-\chi$ in \cite{FFT}.}
\begin{equation}    \label{one mod}
\frac{\langle \chal_{i_j},\nu \rangle}{w_j} -
\sum_{s\neq j} \frac{\langle \chal_{i_j},\al_{i_s} \rangle}{w_j-w_s} +
\langle \chal_{i_j},\chi \rangle = 0, \quad j=1,\ldots,m.
\end{equation}

The corresponding eigenvalues of elements of ${\mc A}_\chi$ are
encoded by the $^L\g$-opers in $\on{Op}_{^L \g}(\pone)_{\pi(\chi)}$
obtained by applying the Miura transformation to the connection
\begin{equation}    \label{olnabla}
\ol\nabla = \pa_z - \chi - \frac{\nu}{z} + \sum_{j=1}^m
\frac{\al_{i_j}}{z-w_j}
\end{equation}
on the $^L H$-bundle $\Omega^{\rho}$ on $\pone$
\end{thm}

Conjecturally, these Bethe eigenvectors form an eigenbasis of $M_\nu$
for generic $\chi$ and $\nu$, and a basis of $V_\nu$ for a dominant
integral $\nu$ and generic $\chi$ (in this case the above weights
should be in the set of weights of $V_\nu$).

Finally, we note that some of the Hamiltonians in ${\mc A}_\chi$ may
be constructed explicitly.  Assume again that $\chi \in \h \subset \g$
is regular semi-simple. In this case the span of all elements of ${\mc
A}_\chi$ of order less than or equal to $1$ (with respect to the
standard filtration on $U(\g)$) is the direct sum of the Cartan
subalgebra $\h = \g_\chi \subset \g \subset U(\g)$ and the scalars $\C
\subset U(\g)$. Next, we have the following quadratic elements:
\begin{equation}\label{eq:T Ham}
T_\gamma(\chi)= \sum_{\alpha\in\Delta_+}
\frac{(\alpha,\gamma)}{(\alpha,\chi)} \left( f_\alpha e_\alpha +
e_\alpha f_\alpha \right), \qquad
\ga \in \h^*,
\end{equation}
where the elements $e_\al,f_\al$ generate the $\sw_2$ subalgebra
$\g_\al$ of $\g$ corresponding to $\al \in \De_+$, and they are
normalized in such a way that $\ka_0(e_\al,f_\al)=1$, where $\ka_0$ is
a non-degenerate invariant inner product on $\g$. These are the {\em
DMT Hamiltonians} introduced by C. De Concini (unpublished) and by
J. Millson and V. Toledano Laredo \cite{MTL,TL} (see also
\cite{FMTV}).\footnote{We note that quasi-classical versions of these
Hamiltonians have appeared in \cite{Boalch} in the study of
isomonodromic equations with irregular singularities.} Explicit
formulas for higher order Hamiltonians are unknown in
general. However, for $\g=\sw_n$ explicit formulas have been found in
\cite{Tar,CT}.

\section{The shift of argument affine Gaudin models}    \label{shift aff}

We now try to generalize the above construction to the situation where
a finite-dimensional simple Lie algebra $\g$ is replaced by the
corresponding (extended) affine Kac--Moody algebra $\gtil$.

In what follows, we will focus on the case when $\gtil$ is an
untwisted affine algebra. However, all of our definitions, results and
conjectures presented below generalize in a straightforward way to the
case of twisted affine algebras.

\subsection{Affine Kac--Moody algebras}    \label{aff km}

Let
\begin{equation}\label{central extension}
0 \to \C {\mb 1} \to \ghat \to \g \ppart \to 0
\end{equation}
be the Kac--Moody central extension of $\g\ppart$ with the commutation
relations
\begin{align*}
\label{comm rel}
[A \otimes f(t),B \otimes g(t)] &= [A,B] \otimes fg - \ka_0(A,B)
\on{Res}_{t=0} f dg \cdot {\mb 1}, \\ [{\mb 1},A \otimes
f(t)] &= 0,
\end{align*}
and $\ka_0$ is the invariant inner product on $\g$ normalized in such
a way that the maximal root has squared length $2$. (We will sometimes
write $(A,B)$ instead of $\ka_0(A,B)$.) Let $\gtil$ be the
corresponding {\em extended affine Kac--Moody algebra}
$$
\gtil = \C d \ltimes \ghat,
$$
where $d$ acts trivially on ${\mb 1}$ and as the outer derivation
$t\pa_t$ on $\g\ppart$:
$$
[d,A \otimes f(t)] = A \otimes t\pa_t f(t).
$$

Let us fix a Cartan decomposition
$$
\g = \n_- \oplus \h \oplus \n_+.
$$
Then the Lie algebra $\gtil$ has the following Cartan decomposition:
$$
\gtil = \ntil_- \oplus \htil \oplus \ntil_+,
$$
where
\begin{align*}
\ntil_+ &= (\n_+ \otimes 1) \oplus (\g \otimes t\C[[t]]),
\\ \ntil_- &= (\n_- \otimes 1) \oplus (\g \otimes t^{-1}
\C[t^{-1}]), \\ \htil &= (\h \otimes 1) \oplus \C{\mb 1} \oplus \C
d
\end{align*}

The Lie algebra $\gtil$ is a symmetrizable Kac--Moody algebra (see
\cite{Kac}). We denote by $\wt{I}$ the set of vertices of the Dynkin
diagram of $\gtil$. Thus, $\wt{I} = I \cup \{ 0 \}$, where $I = \{
1,\ldots,\ell \}$ is the set of vertices of the Dynkin diagram of
$\g$. We have the generators $e_i$ of $\ntil_+$, $f_i$ of $\ntil_-$,
$h_i=\chal_i$ and ${\mb 1}$ of $\htil$, where in all cases $i \in
\wt{I}$. We also denote by $\al_i \in \htil^*, i \in \wt{I}$, the
simple roots of $\gtil$.

A generalization of the shift of argument method to affine Kac--Moody
algebras is well-known to give rise, classically, to Poisson commuting
Hamiltonians of a large class of soliton equations. This point of view
on soliton equations was explained by A. Reyman and
M. Semenov-Tian-Shansky \cite{RS} (see also related works
\cite{Man1,RSF}) and has been widely used ever since. In particular,
this construction was subsequently generalized by V. Drinfeld and
V. Sokolov \cite{DS} to include Hamiltonian reductions (this yields
the generalized KdV hierarchies, as we explain in \secref{kdv}
below). In this section and the next we will discuss the affine Gaudin
models, both classical and quantum, corresponding to the regular
semi-simple shift. They realize the $\g$-AKNS soliton models (or
$n$-wave systems). Then in \secref{kdv} we will consider the
$\ghat$-KdV models corresponding to the Drinfeld--Sokolov reduction of
the shift by the nilpotent maximal root element.

\subsection{The classical model}    \label{aff
  rs}

Consider the case when the shift is a regular semi-simple element
$$
\chi \in \h \subset \g \subset \gtil.
$$
The phase space is then an affine analogue of the space of one-forms
\eqref{eta1}, so it consists of the one-forms
\begin{equation}    \label{eta2}
\eta = \left( \frac{A}{z} + \chi\right) dz, \qquad A \in \gtil^*.
\end{equation}
Instead of the full dual space $\gtil^*$ will consider its codimension
$2$ affine subspace which consists of the functionals taking the value
$k \neq 0$ on the central element ${\mb 1}$ and $0$ on the element
$d$. In other words, this is the affine hyperplane $\ghat^*_k \subset
\ghat^*$ which consists of the linear functionals on $\ghat$ taking
value $k$ on ${\mb 1}$. It is endowed with the Kirillov--Kostant
Poisson structure, and it is well-known that it may be identified with
the space of $k$-connections
\begin{equation}    \label{conn}
k\pa_t + A(t), \qquad A(t) \in \g^*\ppart \simeq \g\ppart
\end{equation}
(where the second isomorphism is obtained using the non-degenerate
invariant inner product $\ka_0$ on $\g$) on the punctured disc, so
that the adjoint action of $G\ppart$ on $\ghat^*_k$ becomes the action
of $G\ppart$ on $k$-connections by gauge transformations (see
\cite{RS,FB}).

By rescaling $A(t) \mapsto \frac{1}{k} A(t)$, we may identify
$\ghat^*_k$ and $\ghat^*_1$, so without loss of generality we may, and
will, set $k=1$.

Thus, we replace $\eta$ in \eqref{eta2} by
$$
\eta = \left( \frac{\pa_t + A(t)}{z} + \chi\right) dz, \qquad A(t) \in
\g\ppart
$$
and further by the equivalent expression
\begin{equation}    \label{eta L}
L = \pa_t + A(t) + \chi z.
\end{equation}
We now recognize it as the $L$-operator of the so-called the $\g$-AKNS
hierarchy, also known as the generalized $n$-wave hierarchy (see,
e.g., \cite{RS,DS}).

The Poisson commuting Hamiltonians of the corresponding integrable
system are obtained by expanding in $z$ (or $z^{-1}$) the invariants
of the monodromy of the connection. Denote by $M(z) \in G$ (where $G$
is the connected simply-connected Lie group with the Lie algebra $\g$)
the formal monodromy matrix of the connection operator $L$ given by
formula \eqref{eta L}. More precisely, $M(z)$ is defined as follows:
its action on any finite-dimensional representation $V$ of $G$ is
obtained by computing the monodromy of the matrix differential
operator $\pa_t + A(t) + \chi z$ acting on $V$.\footnote{Strictly
speaking, this monodromy is defined for connections that converge on a
small punctured disc around $0$; however, the coefficients in the
expansions that we will consider are also well-defined when $A(t)$ is
a formal Laurent power series.}

Now, for any function $\varphi$ on $G$ which is invariant under the
adjoint action of $G$, the Hamiltonian
$$
H_{\varphi}(z) = \varphi(M(z))
$$
is a gauge invariant function of $L$ given by formula \eqref{eta L},
depending on the spectral parameter $z$. One shows in the same way as
for the Hitchin systems in the finite-dimensional case that
\begin{equation}    \label{comm ham}
\{ H_{\varphi}(z),H_{\psi}(y) \} = 0, \qquad \forall
\varphi,\psi \in \on{Fun}(G)^{G}, \; z,y \in \C^\times.
\end{equation}
Therefore we may obtain infinite collections of commuting Hamiltonians
of the affine Hitchin model by expanding the Hamiltonians
$H_{\varphi}(z)$ around the points $0$ and $\infty \in \pone$.

The asymptotic expansion of $H_\varphi(z)$ at $z=\infty$ gives rise to
the local Hamiltonians, in the sense that they may be expressed as
integrals of differential polynomials in the matrix coefficients of
$A(t)$ (in the adjoint representation). There are standard methods for
constructing these Hamiltonians explicitly, see
\cite{RS,DS,FF:ih}. These Hamiltonians generate the generalized
$\g$-AKNS (or $n$-wave) hierarchy. Since the algebra of invariant
functions on $G$ is a polynomial algebra with $\ell = \on{rank}(\g)$
generators, we find that there are $\ell$ local Hamiltonians of each
positive spin.

On the other hand, the expansion of $H_\varphi(z)$ around $z=0$ gives
rise to the non-local Hamiltonians. Let us consider the simplest
examples of such Hamiltonians obtained by taking the $z$-linear terms
of $H_\varphi(z)$.

The computation of these linear terms is based on the following
observation which is well-known in the study of the inverse scattering
problem (for a proof, see, e.g., \cite{Reyman}, Lemma 4.2.4).

\begin{lem}    \label{linear term}
Suppose that we have a connection of the form
\begin{equation}    \label{int sys}
\pa_t + A^{(0)}(t) + \ep A^{(1)}(t) + \ldots, \qquad A^{(i)}(t) \in
\g\ppart.
\end{equation}
Then the $\ep$-linear term in the expansion of $\varphi(M(\ep))$,
where $M(\ep) \in G$ is the monodromy of this connection and $\varphi$
is an invariant function on $G$, is equal to
\begin{equation}    \label{ep lin}
\int (A^{(1)}(t),\Psi_\varphi(t)) dt
\end{equation}
(the integral is taken over a contour around $0$), where $\Psi_\varphi(t)$
is a solution of the equation
\begin{equation}    \label{fund sol}
\pa_t \Psi_\varphi(t) + [A^{(0)}(t)),\Psi_\varphi(t)] = 0,
\end{equation}
which is single-valued (i.e., has no monodromy around $0$), so that
the above contour integral is well-defined.
\end{lem}

More concretely, let us write
$$
A^{(0)}(t) = \sum_a J^{(0)}_a(t) J^a, \qquad A^{(1)}(t) = \sum_a
J^{(1)}_a(t) J^a,
$$
where $\{ J^a \}$ is a basis of $\g$. Let
$\Psi_b^a(t)$ be the matrix elements of the fundamental solution of
the equation \eqref{fund sol}, which may be written explicitly as
follows:
$$
\pa_t \Psi_b^a(t) + \sum_c J^{(0)}_c(t) \sum_d f^{cb}_d \Psi_b^d(t) =
0,
$$
where $\{ f^{cb}_d \}$ are the structure constants of $\g$.
Then the integral in \eqref{ep lin} is equal to
\begin{equation}    \label{ep lin1}
\sum_a \int J^{(1)}_a(t) \Psi_\varphi^a(t) dt,
\end{equation}
where
$$
\Psi_\varphi^a(t) = \sum_b c_b(\varphi) J^b \Psi_b^a(t),
$$
with the coefficients $c_b(\varphi)$ determined by $\varphi$.

Since the algebra of invariant functions on $G$ is isomorphic to the
polynomial algebra in $\ell = \on{rank}(\g)$ generators, the elements
$\Psi_\varphi^a(t)$, where $\varphi$ are invariant functions on $G$,
span an $\ell$-dimensional subspace in the space of all solutions of
\eqref{fund sol}. The corresponding integrals \eqref{ep lin} give rise
to commuting Hamiltonians of an integrable system on the space of
connections of the form \eqref{int sys} (in which $\ep$ is taken as
the spectral parameter).

We now apply this in our case, which fits in the statement of
\lemref{linear term} with $\ep=z$, $A^{(0)}(t) = A(t)$ and $A^{(1)}(t)
= \chi$. Then we obtain that the $z$-linear term in the expansion of
$\varphi(M(z))$, where $\varphi \in (\on{Fun} G)^G$, is given by the
formula
\begin{equation}    \label{H1}
H_\varphi(\chi) = \int \langle \chi,\Psi_\varphi(t) \rangle dt.
\end{equation}
The Poisson commutativity \eqref{comm ham} of the Hamiltonians
$H_{\varphi}$ then implies that
$$
\{ H_{\varphi}(\chi),H_{\psi}(\chi) \} = 0, \qquad \varphi,\psi \in
(\on{Fun} G)^G.
$$

Higher terms in the $z$-expansion of $\varphi(M(z))$ are
given by more complicated formulas, involving multiple integrals of
solutions of the equation
$$
\pa_t \Psi(t) + [A(t),\Psi(t)] = 0
$$
(in the adjoint representation of $\g$). The solutions $\Psi(t)$ of
this equation, and hence all of these Hamiltonians, are {\em
non-local}, in the sense that they cannot be realized as integrals of
differential polynomials in the matrix elements of $A(t)$.

Finally, observe that both local and non-local classical Hamiltonians
Poisson commute with the Heisenberg-Poisson algebra generated by the
Fourier coefficients of the Cartan components of $A(t)$. Indeed, we
obtain our Hamiltonians in effect by differentiating invariant
functions on $\ghat^*_1$ (namely, $\varphi$ of the monodromy) in the
direction of $\chi \in \h^* \subset \ghat^*_1$. Therefore these
Hamiltonians are invariant under the action of the centralizer of
$\chi$ in $\ghat$, which is precisely the (homogeneous) Heisenberg
subalgebra $\wh\h$ of $\ghat$. This is equivalent to the fact that the
Hamiltonians Poisson commute with the Heisenberg-Poisson algebra
generated by $\wh\h$ (in the case of local Hamiltonians this
commutativity is discussed from a different point of view in
\cite{FF:ih}).

\medskip

We now wish to quantize these Hamiltonians. The following remark might
be useful for understanding what follows.

\medskip

Recall that the Hamiltonians of the classical shift of argument Gaudin
model for a finite-dimensional Lie algebra $\g$ are obtained by taking
the derivatives of invariant polynomials on $\g^*$. Each invariant
polynomial $P$ on $\g^*$ has finite degree. Therefore, when we shift
the argument by $z\chi$ and expand in powers of $z$, we obtain a
finite polynomial in $z$, whose coefficients are the derivatives of
the $P$ in the direction of $\chi$. In the affine case the role of
$\g^*$ is played by $\ghat^*_k$, which is the space of connections
\eqref{conn}. The role of invariant polynomials on $\g^*$ is now
played by the spectral invariants $\varphi(M)$ of the monodromy of the
connection. As functions on $\ghat^*_k$, they have ``infinite
degree''. When we apply the shift of argument by $z\chi$, the $n$th
coefficient of the $z$-expansion of the corresponding functions
$\varphi(M(z))$ will be equal to the $n$th derivative of this function
in the direction of $\chi$. This way we obtain the non-local classical
Hamiltonians. On the other hand, it is instructive to think of the
coefficients in the asymptotic $z^{-1}$-expansion as expressions
obtained by differentiating this function $(\infty - n)$ times. These
are the local classical Hamiltonians.

Now, the quantization of the derivatives of the invariant functions on
$\g^*$ in the direction of $\chi$ gives rise to the quantum
Hamiltonians of the shift of argument Gaudin model, as discussed in
\secref{shift fd}. As explained in \cite{Ryb,FFT}, to construct these
Hamiltonians we really need to use the completed enveloping algebra of
$\g\zpart$ (more precisely, of its critical central extension). The
quantum Hamiltonians are constructed, roughly speaking, by taking
derivatives of central elements of this algebra. Likewise, the
quantization of the spectral invariants of the monodromy of the
connection gives rise to the quantum Hamiltonians of the affine Gaudin
models (quantum soliton systems). These Hamiltonians should also be
interpreted as coming from the center of a completed enveloping
algebra, but now of the {\em double affine algebra} $\ghat\zpart$ (or
some extension of it). We expect that this center contains quantum
analogues of the spectral invariants of the monodromy of the
connection, and the local and non-local quantum Hamiltonians may be
obtained by taking the ``derivatives'' of these central elements (as
in the finite-dimensional case), either finitely many times or
infinitely many times.

While at present we do not have a precise definition of the completed
enveloping algebra of $\ghat\zpart$ or its center, the quantum
Hamiltonians of the affine Gaudin models do exist and may be
constructed explicitly, as explained in the next section. We hope that
by analyzing these Hamiltonians we may gain insights into the
structure of the center of the enveloping algebra of $\ghat\zpart$.

\subsection{Quantum Hamiltonians}    \label{affine
shift quantum}

The quantization of the local Hamiltonians has been discussed in our
paper \cite{FF:ih}, in which we explained how to construct them using
the free field (Wakimoto) realization of affine algebras. We refer the
reader to \cite{FF:ih} for more details. According to Conjecture 1 of
\cite{FF:ih}, all of the classical local Hamiltonians of the $\g$-AKNS
hierarchy may be quantized. Therefore we expect to have $\ell$
mutually commuting quantum local Hamiltonians of each positive spin.

In addition, for non-zero values of level $k$ these quantum
Hamiltonians should commute with the (homogeneous) Heisenberg
subalgebra $\wh\h \subset \ghat$, which is the central extension of
$\h\ppart$.

Explicit formulas for the quantum $\g$-AKNS Hamiltonians are unknown
in general. But we can write down those of spin $1$. First of all, we
have the following formula, generalizing formula \eqref{eq:T Ham} for
the DMT Hamiltonians (and using the same notation):
\begin{align}    \label{aff T Ham}
\wh{T}_\gamma(\chi) &= \sum_{\alpha\in\Delta_+}
\frac{(\alpha,\gamma)}{(\alpha,\chi)} \int \Wick e_\alpha(w)
f_\alpha(w) \Wick \, w dw \\ \notag &= \sum_{\wh\al \in
\wh\De^{\on{re}}_+} \frac{(\wh\alpha,\gamma)}{(\wh\alpha,\chi)}
f_{\wh\al} e_{\wh\al}, \qquad \ga \in \h^*,
\end{align}
where $\wh\De^{\on{re}}_+$ is the set of positive real roots of
$\ghat$. Here we set
$$
e_\al(w) = \sum_{n \in \Z} (e_{\al} \otimes t^n) w^{-n-1},
$$
and similarly for $f_\al(w)$, and assume that $\chi$ is regular so
that $(\alpha,\chi) \neq 0$ for all $\al \in \De_+$. The normal
ordering is defined as in \cite{FB}.

These operators mutually commute, but the drawback is that they do not
commute with $\wh\h$. However, we can modify them slightly to make
them commute with $\wh\h$ provided that the level $k$ is
non-zero. Namely, let us set
\begin{multline}    \label{wt T}
\wt{T}_\gamma(\chi) = \\ \sum_{\alpha\in\Delta_+}
\frac{(\alpha,\gamma)}{(\alpha,\chi)} \int \Wick \left( e_\alpha(w)
f_\alpha(w) + f_\alpha(w) e_\alpha(w) - \frac{1}{k} h_\al(w)^2 \right)
\Wick \, w dw, \qquad \ga \in \h^*,
\end{multline}
where $h_\al \in \h$ is the element corresponding to the root $\al
\in \h^*$ via the normalized bilinear form.

The following lemma will be proved in the Appendix.\footnote{The
operators $\wh{T}_\gamma(\chi)$ may be obtained as certain limits of
the connection operators studied in \cite{FMTV}, and part (1) of
\lemref{comm DMT} may be derived from the commutativity of these
operators proved in \cite{FMTV} (see \secref{gen km} below for more
details).}

\begin{lem}    \label{comm DMT} \hfill

{\em (1)} For any regular $\chi \in \h^*$ the Hamiltonians
$\wh{T}_\gamma(\chi), \ga \in \h^*$, commute with each other.

{\em (2)}  For any regular $\chi \in \h^*$ the Hamiltonians
$\wt{T}_\gamma(\chi), \ga \in \h^*$, commute with each other and with
the Heisenberg subalgebra $\wh\h \subset \ghat$.

\end{lem}

\medskip

Next, we discuss the quantization of the non-local Hamiltonians. In
order to do this, we need to construct quantum analogues of the
solutions $\Psi(t)$ of the differential equation
\begin{equation}    \label{sol}
(k \pa_t + A(t)) \Psi(t) = 0,
\end{equation}
taking values in the adjoint representation of $\g$. It has been shown
by V. Knizhnik and A. Zamolodchikov in \cite{KZ} that these quantum
analogues may be obtained from the primary fields of the WZW model
corresponding to the adjoint representation. We will now recall how
this works for an arbitrary finite-dimensional representations $V$ of
$\g$. In what follows we will assume that $k$ is a complex number not
equal to the critical value $-h^\vee$.

By definition, for given highest weight modules $M_1$
and $M_2$ over $\ghat$ of the same level $k$, the primary field
corresponding to $V$ is a formal power series
$$
\Phi_V(w) = \sum_{n \in \Z} \Phi_{V,(n)} w^{-n}: M_1 \otimes V \to
M_2,
$$
which is $\ghat$-invariant, where $\ghat$ acts on $V$ by evaluation at
$w$,
$$
J^a_n \overset{\on{def}}= J^a \otimes t^n \mapsto w^n J^a.
$$
This is equivalent to the Fourier coefficients $\Phi_{V,(n)}$ having
the following commutation relations with the affine algebra $\ghat$:
$$
[J^a_m,\Phi_{V,(n)}] = J^a \cdot \Phi_{V,(n+m)}.
$$
Here on the left hand side we consider the action of $J^a_n$ on $M_1$
and $M_2$, and on the right hand side the action of $J^a$ on $V$.

If $M_1$ and $M_2$ have highest weights $\la_1$ and $\la_2$ and $V =
V_\mu$ has highest weight $\mu$, then the corresponding conformal
dimensions (the highest weights of the Segal--Sugawara Virasoro
algebra) are equal to $\Delta(\la_1), \Delta(\la_2)$, and
$\Delta(\mu)$, respectively, where $\Delta(\nu) = C(\mu)/(k+h^\vee)$,
and $C(\nu)$ is the scalar by which the normalized Casimir element of
$U(\g)$ acts on $V_\nu$. The true primary field is then
$$
\wt\Phi_V(w) = w^{\Delta(\la_2)-\Delta(\la_1)-\Delta(\mu)} \Phi_V(w).
$$
As explained in \cite{KZ}, it satisfies the following differential
equation:
\begin{equation}    \label{quantum sol}
(k+h^\vee) \pa_w \Phi_V(w) - \sum_a \Wick J^a(w) (J_a \cdot \Phi_V)(w)
\Wick  = 0,
\end{equation}
which may be viewed as a quantization of the differential equation
\eqref{sol} written in the representation $V$ (note the ``quantum
shift'' $k \mapsto k+h^\vee$). Therefore we will consider the primary
fields $\Psi_V(w)$ as quantizations of solutions of the equation
\eqref{sol}. Moreover, we see that the monodromy of this
$\wt\Phi_V(w)$ is equal to $\exp 2\pi
i(\Delta(\la_2)-\Delta(\la_1)-\Delta(\mu))$.

Suppose now that $M_1$ and $M_2$ are Verma modules $M_{\la_1,k}$ and
$M_{\la_2,k}$. It is known that for generic $k$ the space of primary
fields of type $V = V_\mu$ acting from $M_{\la_1,k}$ to $M_{\la_2,k}$
has the dimension equal to the dimension of the component of weight
$\la_2-\la_1$ in $V_\mu$. Moreover, it is known that this space is
isomorphic to the component of weight $\la_2-\la_1$ in the irreducible
representation of the quantum group $U_q(\g)$ with highest weight
$\mu$, where $q=\exp(\pi i/(k+h^\vee))$.

In particular, we find that if $M_1=M_2=M_\la$, then the space of
primary fields $\Phi(w) = \Phi_{\g}(w)$ of type $\g$ (the adjoint
representation) acting from $M_\la$ to itself is isomorphic to the
weight $0$ subspace of $\g$, which is the Cartan subalgebra of
$\g$. It therefore has dimension $\ell$, the rank of $\g$. It is clear
that this subspace is precisely the quantization of the space of
single-valued solutions of the equation \eqref{sol}, spanned by
$\Psi_\varphi, \varphi \in (\on{Fun} G)^G$. Let us choose a basis
$\Phi_p(w), p=1,\ldots,\ell$, of this space of primary fields.

Then the vector space of the quantum Hamiltonians that quantize the
Hamiltonians $H_p(\chi)$ given by formula \eqref{H1} is just the
span of the $0$th Fourier coefficients of the fields $\Phi_p(w)$,
\begin{equation}    \label{tilde Psi}
H_p(\chi) = \int \Phi_p(w) \; \frac{dw}{w} = \int \wt\Phi_p(w)
w^{h^\vee/(k+h^\vee)} \; \frac{dw}{w}
\end{equation}
of the primary fields $\Phi_p(w): M_\la \to M_\la$. Here we use
the fact that
$$
\Delta(\theta) = h^\vee/(k+h^\vee),
$$
where $\theta$ is the maximal root, which is the
highest weight of the adjoint representation of $\g$.

We claim that the Hamiltonians $H_p(\chi), p = 1,\ldots,\ell$,
commute with each other and with the local Hamiltonians \eqref{aff T
Ham}. This is checked by a direct calculation which will be presented
in \cite{FF:next}.

The quantization of the classical non-local Hamiltonians,
corresponding to the higher terms in the $z$-expansion of the
monodromy $\varphi(M(z))$ is expressed by multiple integrals of the
adjoint primary fields $\wt\Phi(w) w^{h^\vee/(k+h^\vee)}$ (including
the components acting between different Verma modules). We will now
briefly sketch the construction of these integrals. The details will
appear in \cite{FF:next}.

This construction is analogous to the construction of \cite{BLZ1,BHK},
where non-local quantum KdV Hamiltonians were constructed using
bosonic screening operators corresponding to the simple roots of
$\ghat$. These operators satisfy the $q$-Serre relations, and hence,
in a certain sense, generate the quantized enveloping algebra
$U_q(\ntil_-)$, where $\ntil_-$ is the lower nilpotent Lie subalgebra
of $\ghat$ (see \cite{BMP,FF:laws} for the precise definition of what
this means). The non-local quantum KdV Hamiltonians are then obtained
from elements of $U_q(\ntil_-)$ which correspond to singular vectors
in the Verma modules over $U_q(\ghat)$ of {\em critical level}.

More precisely, as explained in \cite{FF:laws}, whenever we have a
singular vector of weight $\mu$ in the Verma module $M^q_{\la}$ over
$U_q(\ghat)$ of highest weight $\la$ and level $0$, we have a
well-defined operator $\pi_\la \to \pi_\mu$, where $\pi_\la$ denotes
the Fock representation of a Heisenberg Lie algebra with highest
weight $\la$. This operator is constructed as follows. Let us write
the singular vector as $P v_\la$, where $v_\la \in M^q_\la$ is the
highest weight vector and $P \in U_q(\ntil_-)$. This operator is
obtained by substituting the bosonic screening operators
$$
S_i = \int S_i(w) \; \frac{dw}{w}, \qquad i=0,\ldots,\ell,
$$
where the $S_i(w) = e^{\beta \phi_i(w)}$ are the bosonic primary
fields corresponding to the simple roots in $\ghat$, into $P$ instead of
the generators $f_i, i=0,\ldots,\ell$. Here $q = \exp(\pi i
\beta^2)$. Now, if we replace the screening operator $S_0$ by the
operator
\begin{equation}    \label{S0}
\wt{S}_0 = \int S_0(w) \; w^{h^\vee \beta^2} \; \frac{dw}{w},
\end{equation}
and keep $S_i, i=1,\ldots,\ell$, unchanged, then the same procedure
will work, except that we will now assign operators acting on Fock
representations to singular vectors in the Verma modules of {\em
critical level} $k=-h^\vee$ instead of level $k=0$. There are many
more of those than at $k=0$, and one can check that they commute with
each other. These are the non-local quantum Hamiltonians of
\cite{BLZ1,BHK,FFS}.

The non-local Hamiltonians of the affine shift of argument Gaudin
model may be constructed in a similar fashion. The most direct way to
do it is to replace the Virasoro (or, more generally, ${\mc
W}$-algebra) screening operators $S_i$ with the screening operators
for the affine algebra $\ghat$ arising from the free field (Wakimoto)
realization of $\ghat$. These screening operators have been defined in
\cite{FF:ih}, where it was shown that they satisfy the same Serre
relations of $U_q(\ntil_-)$. Hence, substituting them in the formulas
discussed above, we construct non-local commuting Hamiltonians of our
model.

However, we want to stress that it is not necessary to use the free
field realization to define these quantum non-local Hamiltonians. They
may be defined intrinsically in terms of the adjoint primary fields of
the WZW model associated to $\ghat$. In order to do that, we need to
replace the algebra $U_q(\ntil_-)$ by a smaller quantum algebra, which
is a quantization of $U(\gtil_-)$, where $\gtil_- =
t^{-1}\g[t^{-1}]$. To explain this more precisely, let us recall that
the algebra $U_q(\ntil_-)$ may be constructed using certain rank one
local systems on the configuration spaces (see
\cite{SV1,Var}). Moreover, homology classes with coefficients in these
local systems correspond to singular vectors discussed above. These
local systems describe the monodromy properties of the above screening
currents, and the non-local Hamiltonians may be interpreted as
integrals of products of the screening currents over homology cycles
of these local systems (see \cite{FF:laws} for more details).

Likewise, the quantization of $U(\gtil_-)$ that we need here, which we
denote by $U_q(\gtil_-)$, may be constructed using local systems on
the configuration spaces that capture the monodromy properties of the
adjoint primary fields $\Phi_{\g}(w)$. These local systems are defined
on the trivial vector bundles on the configuration spaces $\C^N \bs
\on{diag}$ with the fibers $\g^{\otimes N}$ and the connection given
by the corresponding KZ operators. The monodromies of these local
systems around the diagonals are given by the $R$-matrices of
$U_q(\g)$ in the adjoint representation. Now, the commuting non-local
Hamiltonians quantizing the coefficients in the $z$-expansion of
$\varphi(M(z))$ are constructed by integrating products of the primary
fields over homology classes with coefficients in these local systems,
which in turn correspond to singular vectors in
$U_q(\gtil_-)$-modules. Since we want to use singular vectors at the
critical level $k=-h^\vee$ (and not $k=0$), we need to insert the
factor $w^{h^\vee/(k+h^\vee)}$ in the integral \eqref{tilde Psi}. This
is analogous to inserting the factor $w^{h^\vee \beta^2}$ in formula
\eqref{S0}.

The simplest of these non-local Hamiltonians are the operators
$H_p(\chi)$ given by formula \eqref{tilde Psi}. They correspond to the
quantization of the $z$-linear term in the expansion of
$\varphi(M(z))$. Higher degree terms correspond to integrals over more
complicated homology classes with coefficients in the above local
system. Applying the free field realization of $\ghat$, we may express
these Hamiltonians as contour integrals of the $\ghat$ screening
currents and obtain formulas similar to those of \cite{BLZ1,BHK}.

\subsection{Bethe Ansatz}

Next, we discuss the eigenvectors of the quantum Hamiltonians. The
Bethe Ansatz described in \secref{spectra} has a straightforward
generalization to this case. The formula \eqref{bethe for one mod} for
the Bethe vector and the Bethe Ansatz equations \eqref{one mod}
easily generalize to the affine case. We just need to enlarge the set
$I$ to $\wt{I}$ and include the generator $f_0$ (corresponding to the
simple root $\al_0$ of $\gtil$) and view $\nu$ as an element of
$\htil^*$. Then for any collection of distinct non-zero complex
numbers $w_j, j=1,\ldots,m$, and labels $i_1,\ldots,i_m \in \wt{I}$ of
the Dynkin diagram of $\gtil$, the vector
$\phi(w_1^{i_1},\ldots,w_m^{i_m})$ given by formula \eqref{bethe for
one mod} belongs to the Verma module $M_\nu$ over $\ghat$.

We expect that if the Bethe Ansatz equations \eqref{one mod} are
satisfied, then the vector $\phi(w_1^{i_1},\ldots,w_m^{i_m})$ is an
eigenvector of the quantum Hamiltonians discussed above (both local
and non-local). Since we expect that these Hamiltonians commute with
the Heisenberg Lie agebra $\wh\h$, it is natural to assume that these
vectors are also invariant under the Lie subalgebra $t\h[[t]] \subset
\wh\h$ (we have checked that this is indeed true in the simplest
examples). Thus, we arrive at the following conjecture.

\begin{conj}    \label{h inv}
If the Bethe Ansatz equations \eqref{one mod} are satisfied, then the
vector $\phi(w_1^{i_1},\ldots,w_m^{i_m})$ is a $t\h[[t]]$-invariant
eigenvector of all local and non-local quantum affine Gaudin
Hamiltonians, in particular, the Hamiltonians $\wt{T}_\ga, \ga \in
\h^*$, given by formula \eqref{wt T}, and $H_p(\chi),
p=1,\ldots,\ell$, given by formula \eqref{tilde Psi}.
\end{conj}

Note that $\phi(w_1^{i_1},\ldots,w_m^{i_m})$ is a homogeneous vector
of weight $\la - \sum_{j=1}^m \al_{i_j}$. It is clear from formulas
\eqref{wt T} and \eqref{aff T Ham} that any homogeneous
$t\h[[t]]$-invariant eigenvector of $\wt{T}_\ga, \ga \in \h^*$, is
automatically an eigenvector of $\wh{T}_\ga, \ga \in \h^*$. Therefore
\conjref{h inv} implies the following statement.

\begin{conj}    \label{phi eig}
The vector $\phi(w_1^{i_1},\ldots,w_m^{i_m})$ is an eigenvector of the
affine DMT Hamiltonians \eqref{aff T Ham} if the Bethe Ansatz
equations \eqref{one mod} are satisfied.
\end{conj}

To prove this conjecture, we cannot use the argument used in
\cite{FFT} because we do not know that the affine DMT Hamiltonians
arise by the same mechanism as the DMT Hamiltonians in the
finite-dimensional case (that is, from the center of an enveloping
algebra of the affinized Lie algebra). But this conjecture can
probably be derived from the results of \cite{FMTV} (see \secref{gen
km} below for more details).

\medskip

The next question is to describe the common eigenvalues of the quantum
Hamiltonians that occur on highest weight $\ghat$-modules, such as
Verma modules and irreducible integrable representations. In
particular, we wish to describe the eigenvalues on the Bethe
eigenvectors.

According to \thmref{spectra fd}, in the finite-dimensional case the
common eigenvalues of the quantum shift of argument algebra ${\mc
A}_\chi$ are encoded by $^L \g$-opers on $\pone$ satisfying certain
conditions. In particular, the common eigenvalues on the Bethe vectors
\eqref{bethe for one mod} are encoded by the $^L \g$-opers obtained by
applying the Miura transformation to the connection
\eqref{olnabla}. We wish to generalize these results to the affine
case. This will be done in the next section.

\section{Affine opers}    \label{affine opers}

The notion of $\g$-opers, where $\g$ is a finite-dimensional simple
Lie algebra, was introduced by Beilinson and Drinfeld
\cite{BD,BD:opers} following Drinfeld and Sokolov \cite{DS}. Denote by
$G$ the Lie group of inner automorphisms of the Lie algebra $\g$ and
$B$ its Borel subgroup. A $\g$-{\em oper} on a smooth curve $X$ is by
definition a triple $(\F,\nabla,\F_B)$, where $\F$ is a principal
$G$-bundle $\F$ on $X$, $\nabla$ is a connection on $\F$ and $\F_B$
is a $B$-reduction of $\F$ such that locally on $X$, if we choose a
local coordinate $z$ and trivialize $\F_B$, the connection acquires
the form
\begin{equation}    \label{form of nabla}
\nabla = \pa_z + \sum_{i=1}^\ell \psi_i(z) f_i + {\mb v}(z),
\end{equation}
where each $\psi_i(z)$ is a nowhere vanishing function, and ${\mb
v}(z)$ is a $\bb$-valued function.

In this section we generalize the notion of $\g$-opers and define
$\ghat$-opers, where $\ghat$ is an affine Kac--Moody algebra,
following \cite{F:flag}. We then use these affine opers to give a
conjectural description of the spectra of the quantum Hamiltonians of
the affine Gaudin model introduced in the previous section.

\subsection{Definition}

We now consider another copy of the affine Kac--Moody algebra $\gtil$,
which we will use to describe the spectra of the quantum affine Gaudin
Hamiltonians. It is important to realize that this affine Kac--Moody
algebra has nothing to do with the affine Kac--Moody algebra used in
the construction of the Gaudin models above (in fact, we expect that
the two are Langlands dual to each other). To emphasize this
difference, we will now use $\la$ as the loop variable. Thus, in this
section we set
$$
\gtil = \C d \ltimes \ghat,
$$
where
$$
0 \to \C {\mb 1} \to \ghat \to \g \pparl \to 0.
$$
It will be convenient to use the following slightly unusual
convention for the Cartan decomposition:
$$
\gtil = \ntil_- \oplus \htil \oplus \ntil_+,
$$
where
\begin{align*}
\ntil_+ &= (\n_+ \otimes 1) \oplus (\g \otimes
\la^{-1}\C[\la^{-1}]), \\ \ntil_- &= (\n_- \otimes 1) \oplus (\g
\otimes \la \C[[\la]]), \\ \htil &= (\h \otimes 1) \oplus \C{\mb
1} \oplus \C d.
\end{align*}

Now let $\wt{G}$ be the Lie group (more precisely, ind-group scheme)
associated to $\gtil$ (see, e.g., \cite{Kash} for the precise
definition). We will denote by $\ol{G}$ the quotient of $\wt{G}$ by
its center. Thus, the Lie algebra of $\ol{G}$ is
$$
\gbar = \C d \ltimes \g\pparl.
$$
The group $\ol{G}$ comes with the lower unipotent and Borel subgroups
$\ol{N}_-$ and $\ol{B}_-$ (which are proalgebraic groups)
corresponding to the Lie subalgebras
$$
\nbar_- = \ntil_- \qquad \on{and} \qquad \bbar_- =
(\h \oplus \C d) \oplus \ntil_-
$$
of $\gbar$, respectively, and the upper unipotent and Borel subgroups
$\ol{N}_+$ and $\ol{B}_+$ (which are group ind-schemes) corresponding
to
$$
\nbar_+ = \ntil_+ \qquad \on{and} \qquad \bbar_+ = (\h \oplus \C d)
\oplus \ntil_+,
$$
respectively.

Let $X$ be a smooth algebraic curve. We will use the notion of
principal $\ol{G}$-bundle on $X$ and a connection on it given in
\cite{F:flag}. Roughly speaking, such a bundle may be trivialized on
open subsets of $X$ (in either Zariski or analytic topology). Having
chosen such a trivialization and a local coordinate $z$ on a an open
subset $U$, a connection is by definition a first order differential
operator $\nabla = \pa_z + A(z)$, where $A(z) \in \gtil$. If we change
the trivialization by $g(z) \in \wt{G}$, $\nabla$ transforms by the
usual formula
$$
\nabla \mapsto \pa_z + g A(z) g^{-1} - (\pa_z g) g^{-1}.
$$
Under a change of coordinates $z = \varphi(s)$ the operator $\nabla$
also transforms in the usual way:
$$
\nabla \mapsto \pa_s + \varphi'(s) A(\varphi(s)).
$$

Now, by analogy with the finite-dimensional case
\cite{DS,BD,BD:opers}, a $\ghat$-{\em oper} on $X$ may be defined as a
triple $(\F,\nabla,\F_{\ol{B}_+})$, where $\F$ is a principal
$\ol{G}$-bundle $\F$ on $X$, $\nabla$ is a connection on $\F$ and
$\F_{\ol{B}_+}$ is a $\ol{B}_+$-reduction of $\F$, such that locally,
with any choice of local coordinate $z$ and trivialization of
$\F_{\ol{B}_+}$, the connection operator has the form
\begin{equation}    \label{oper new}
\nabla = \pa_z + \sum_{i=1}^\ell \psi_i(z) f_i + {\mb v}(z),
\end{equation}
where each $\psi_i(z)$ is a nowhere vanishing function, and ${\mb
v}(z)$ is a $\bbar_+$-valued function.

We note that another version of affine opers, introduced earlier in
\cite{BF}, has $\bbar_+$ and $\bbar_-$ switched. For comparison of the
two definitions and the explanation as to why the definition given
here is the one relevant to the Gaudin model, see \cite{F:flag}.

However, this naive definition appears to be inadequate. The reason is
that for ind-groups, such as $\ol{G}$, a connection does not
necessarily give rise to a trivialization of $\F$, even locally
analytically. This is because the equation
$$
(\pa_z + A(z)) \Psi(z) = 0,
$$
where $A(z)$ takes values in $\ol\g \simeq \C d \ltimes \g\pparl$,
does not {\em a priori} have local solutions with values in $\ol{G}
\simeq \C^\times \ltimes G\pparl$ (we may well obtain power series in
$\la$ unbounded from below, which are not well-defined as elements of
$\ol{G}$). This is related to the fact that we do not have an
exponential map from the Lie algebra $\ol\g$ to the group
$\ol{G}$. This means that the usual correspondence between flat
connections on $\ol{G}$-bundles and $\ol{G}$-local systems does not
exist. In particular, the notion of monodromy of an abstract oper as
defined above does not make sense.

Our experience with opers suggests that their monodromy properties are
important in the context of Gaudin models (see
\cite{F:flag,F:faro}). In order to have well-defined monodromy, we
will assume that, locally on $X$, the oper bundle $\F$ admits a
reduction to a {\em proalgebraic} subgroup of $\ol{G}$, such as
$G[[\la]]$, which is preserved by the connection $\nabla$. Then the
connection $\nabla$ does give rise to local analytic trivializations
of this $G[[\la]]$-subbundle, and hence of the $\ol{G}$-bundle $\F$
as well, so the monodromy of the oper connection is
well-defined. Thus, from now on by a $\ghat$-oper on $X$ we will
understand a triple $(\F,\nabla,\F_{\ol{B}_+})$ as above such that,
locally on $X$, $\F$ admits a horizontal reduction to $G[[\la]]$ or to
another ``maximal compact'' subgroup of $G\pparl$. We denote the set
of $\ghat$-opers on $X$ by $\on{Op}_{\ghat}(X)$.

Let
$$
p_{-1} = \sum_{i \in \wt{I}} f_i = \sum_{i \in I} f_i +
e_{\theta} \otimes \la.
$$
We will say that a $\ghat$-oper on $X \bs x$ has regular singularity
at $x$ if the connection operator may be brought to the form
\begin{equation}    \label{oper with RS2}
\nabla = \pa_z + \frac{1}{z} \left( p_{-1} + {\mb v}(z) \right),
\qquad {\mb v}(z) \in \bbar_+[[z]],
\end{equation}
on the punctured disc $D_x^\times$ near $x$, were $z$ is a local
coordinate at $x$. We define the residue of such an oper as in the
finite-dimensional case (see, e.g., \cite{F:flag}), so that it takes
values in $\ol\h/W_{\on{aff}}$. We will denote by $\ol\varpi$ the
projection $\ol\h \to \ol\h/W_{\on{aff}}$.

Let $\nu$ be a dominant integral coweight such that $\nu_i = \langle
\nu,\chal_i \rangle \in \Z_+$ for all $i \in \wt{I}$. We will say that
a $\ghat$-oper on $X$ is $\nu$-{\em regular} at $x \in X$ if it may be
brought to the form
\begin{equation}    \label{psi la new1}
\nabla = \pa_z + \sum_{i \in \wt{I}} z^{\nu_i} f_i + {\mb v}(z),
\end{equation}
where ${\mb v}(z) \in \bbar_+[[z]]$, on $D_x^\times$.

We also define opers with singularity of order $2$ and their
$2$-residue similarly to the finite-dimensional case (see
\cite{FFT}).

\subsection{Miura opers and Miura transformation}
\label{miura}

Let $\ol{H}$ be the subgroup of $\ol{G}$ corresponding to the Lie
algebra $\ol\h = \C d \oplus \h$. Let $\crho$ be the element of
$\ol\h$ such that $\langle \al_i,\crho \rangle = 1$ for all $i \in
\wt{I}$. It defines a homomorphism $\C^\times \to \ol{H}$. Let
$\Omega^{\crho}$ be the $\ol{H}$-bundle on $X$ obtained as the
push-forward of the $\C^\times$-bundle on $X$ corresponding to the
canonical line bundle on $X$. Connections on $\Omega^{\crho}$ are
described as follows. If we choose a local coordinate $z$ on an open
subset of $X$, then we trivialize $\Omega^{\crho}$ and represent the
connection as an operator
\begin{equation}    \label{cartan conn}
\ol\nabla = \pa_z + {\mb u}(z), \qquad {\mb u}(z) \in \ol\h.
\end{equation}
If $s$ is another coordinate such that $z=\varphi(s)$, then this
connection will be represented by the operator
\begin{equation}    \label{trans for conn1}
\pa_s + \varphi'(s) {\mb u}(\varphi(s)) - \crho \cdot
\frac{\varphi''(s)}{\varphi'(s)}.
\end{equation}
We will call them {\em Cartan connections}.

We define the {\em Miura transformation} from Cartan connections to
$\ghat$-opers taking an operator \eqref{cartan conn} to the gauge
equivalence class of
\begin{equation}    \label{miura oper}
\nabla = \pa_z + p_{-1} + {\mb u}(z).
\end{equation}

The Miura transformation may be understood more conceptually in terms
of the {\em Miura opers}. By definition \cite{F:wak}, a Miura
$\ghat$-oper on $X$ is a quadruple
$(\F,\nabla,\F_{\ol{B}_+},\F_{\ol{B}_-})$, where
$(\F,\nabla,\F_{\ol{B}_+})$ is a $\ghat$-oper on $X$ and
$\F_{\ol{B}_-}$ is a $\ol{B}_-$-reduction of $\F$ which is preserved
by $\nabla$. It is called generic if the reductions $\F_{\ol{B}_+}$
and $\F_{\ol{B}_-}$ are in generic relative position (see
\cite{F:wak}).  We have a map from Cartan connections to generic
Miura opers taking an operator \eqref{cartan conn} to the quadruple in
which $\F$ is the trivial bundle with the connection \eqref{miura
oper} and the tautological reductions $\F_{\ol{B}_+},\F_{\ol{B}_-}$
satisfying the above conditions. This map sets up a bijection
between generic Miura opers and Cartan connections (see
\cite{F:wak}). The inverse map sends the quadruple as above to the
$\ol{H}$-bundle $\F_{\ol{B}_-}/\ol{N}_-$ with the connection induced
by $\nabla$.

An arbitrary Miura oper on $X$ is always generic away from finitely
many points of $X$. At each of the remaining points the reductions
$\F_{\ol{B}_+}$ and $\F_{\ol{B}_-}$ are in relative position labeled
by an element $w \neq 1$ of the (affine) Weyl group $W_{\on{aff}}$ of
$\ol{G}$. The corresponding Cartan connection has a singularity at
such a point with the residue equal to $\crho-w(\crho)$ (see
\cite{F:flag} for more details).

The Miura transformation is thus defined as the composition of the map
from the Cartan connections to generic Miura opers and the forgetful
map from Miura opers to opers.

\subsection{Parabolic Miura opers}    \label{par miura}

Finally, we introduce the notion of {\em parabolic Miura oper}. By
definition, a parabolic Miura $\ghat$-oper on $X$ is a quadruple
$(\F,\nabla,\F_{\ol{B}_+},\F_{\ol{G}_-})$, where
$(\F,\nabla,\F_{\ol{B}_+})$ is a $\ghat$-oper on $X$ and
$\F_{\ol{G}_-}$ is a reduction of $\F$ to $\ol{G}_- = \C^\times
\ltimes G[[\la]]$ which is preserved by $\nabla$. We will call it
generic if $\F_{\ol{B}_+}$ and $\F_{\ol{G}_-}$ are in generic relative
position. We then have the following map from the space
$\on{Op}_{\g}(X)$ of $\g$-opers on $X$ to the space of generic
parabolic Miura $\ghat$-opers: it sends a $\g$-oper represented by the
connection
$$
\ol\nabla = \pa_z + \sum_{i \in I} f_i + {\mb v}(z), \qquad {\mb v}(z)
\in \bb_+,
$$
to the quadruple in which $\F$ is the trivial bundle,
$$
\nabla = \ol\nabla + f_0,
$$
and $\F_{\ol{B}_+},\F_{\ol{G}_-}$ are the tautological reductions of
$\F$. However, this map is not bijective, because we are missing the
one-dimensional subgroup $\C^\times$ of $\ol{G}_-$ corresponding to
the element $d \in \ol\h$.

In order to obtain a bijection, consider the Cartesian product of
$\on{Op}_{\g}(X)$ and the space $\on{Conn}_X(\Omega^{h})$, where $h =
\langle \delta,\crho \rangle$ is the Coxeter number of $\g$.  Define a
map from $\on{Op}_{\g}(X) \times \on{Conn}_X(\Omega^{h})$ to the
space of generic parabolic Miura $\ghat$-opers sending a pair
consisting of an oper $\ol\nabla \in \on{Op}_{\g}(X)$ and a connection
$\ol\nabla' \in \on{Conn}_X(\Omega^{h})$ of the form
$$
\ol\nabla' = \pa_z + \ga(z),
$$
to the connection
$$
\nabla = \ol\nabla + f_0 + \ga(z) d.
$$
This map is already bijective. To construct the inverse map, observe
that
$$
\ol{G}_-/G^{(1)} \simeq G \times \C^\times
$$
(where $G^{(1)}$ is the first congruence subgroup of $G[[\la]]$). To a
quadruple $(\F,\nabla,\F_{\ol{B}_+},\F_{\ol{G}_-})$ as above, we
associate a $G \times \C^\times$-bundle $\F/G^{(1)}$ together with a
connection and a $B$-reduction $\F_{\ol{B}_+}/G^{(1)}$. This bundle is
a product of the $\g$-oper bundle and the $\C^\times$-bundle
corresponding to $\Omega^h$. Therefore these data give rise to a point
in $\on{Op}_{\g}(X) \times \on{Conn}_X(\Omega^{h})$.

Taking the composition with the forgetful map from parabolic Miura
opers to opers, we now obtain a map from $\on{Op}_{\g}(X) \times
\on{Conn}_X(\Omega^{h})$ to $\ghat$-opers, which is a parabolic
version of the Miura transformation. This gives us a convenient way to
represent $\ghat$-opers in terms of $\g$-opers, which we will exploit
in what follows.

Note that any parabolic Miura $\ghat$-oper on a curve $X$ is generic
away from finitely many points of $X$. At each of those points the
relative position of $\F_{\ol{B}_+}$ and $\F_{\ol{G}_-}$ is measured
by a coset $\ol{w} \in W \bs W_{\on{aff}}$, and the corresponding
$\g$-oper has regular singularity with the residue
$\varpi(-\ol{w(\crho)})$, where the bar stands for the projection
$\ol\h \to \h$ along $\C d$ and $\varpi$ is the projection $\h \to
\h/W$, as before (note that $\varpi(-\ol{w_1(\crho)}) =
\varpi(-\ol{w_2(\crho)})$ if $w_1 \sim w_2$ in $W \bs W_{\on{aff}}$,
so this formula is well-defined).

\begin{remark}
There is an analogous notion of parabolic Miura oper in the
finite-dimen\-sional case. Let $P$ be a parabolic subgroup of
$G$. Then we define a Miura $(\g,P)$-oper to be a quadruple
$(\F,\nabla,\F_B,\F'_P)$, where $(\F,\nabla,\F_B)$ is a $\g$-oper (see
the beginning of this section) and $\F'_P$ is a $P$-reduction of $\F$
which is horizontal with respect to $\nabla$. (In the affine case the
role of $P$ is played by $\ol{G}_-$.) We have the obvious notion of
generic Miura $(\g,P)$-oper.

Let $L = P/[P,P]$. Then we define, in the same way as above, a
bijective map from the space of $L$-opers to generic Miura
$(\g,P)$-opers. Composing it with the forgetful map from Miura
$(\g,P)$-opers to $\g$-opers, we obtain a parabolic analogue of the
Miura transformation map from $L$-opers to $\g$-opers. The ordinary
Miura opers correspond to the case $P=B$.\qed
\end{remark}

\subsection{Spectra of quantum Gaudin Hamiltonians and affine opers}
\label{spectra aff}

Now let $^L \ghat$ be the Langlands dual Lie algebra of $\ghat$, whose
Cartan matrix is the transpose of the Cartan matrix of $\ghat$ (thus,
$^L\ghat$ is a twisted affine Lie algebra if $\ghat$ is non-simply
laced). We will use the upper index $L$ to denote all of the above
groups corresponding to $^L \ghat$, such as $^L \h$, $^L \ol{H}$,
etc. In particular, we have a canonical identification $^L \h = \h^*$,
and so we will consider a weight $\nu \in \h^*$ as an element of $^L
\h$. By abusing notation, we will denote by the same symbol its
projection onto $^L \ol\h = {}^L \wt\h/\C {\mb 1} = {}^L \h \oplus \C
d$.

As in the finite-dimensional case (see \cite{FFT} and
\secref{spectra}), we attach to each solution of the Bethe Ansatz
equation \eqref{one mod}, where $i_j \in \wt{I}$, the Cartan
$^L \ol{H}$-connection
\begin{equation}    \label{olnabla1}
\ol\nabla = \pa_z - \chi - \frac{\nu}{z} + \sum_{j=1}^m
\frac{\al_{i_j}}{z-w_j}.
\end{equation}
In the same way as in \cite{F:flag}, Lemma 2.10, we prove the
following:

\begin{lem}    \label{bethe and mon}
The Bethe Ansatz equations \eqref{one mod} are
equivalent to the condition that the $^L \ghat$-oper $\nabla =
\ol\nabla + p_{-1}$ obtained by applying the Miura transformation to
the connection \eqref{olnabla1} is regular at the points $w_j,
j=1,\ldots,m$, or, equivalently, has no monodromies around these
points.
\end{lem}

Note that this $^L \ghat$-oper has regular singularity at the point
$z$ with the residue and $\infty$ with the residue determined by
$\nu$, and order $2$ singularity at $\infty$ with the $2$-residue
$\pi(\chi)$.

This suggests that the true parameters for the common eigenvalues of
the quantum Gaudin Hamiltonians are not the Cartan connections
\eqref{new miura}, but the corresponding $^L \ghat$-opers obtained by
applying the Miura transformation. Thus, we propose the following
description of the spectra of the quantum affine Gaudin Hamiltonians
(including the affine DMT Hamiltonians $\wh{T}_\ga(\chi)$ and the
non-local Hamiltonians $H_p(\chi)$):

\bigskip

\hspace*{4mm} \parbox{.9\linewidth}{{\em The common eigenvalues of the
quantum affine Gaudin Hamiltonians on the subspace of
$t\h[[t]]$-invariant vectors in an irreducible $\ghat$-module with
highest weight $\nu$ and level $k \neq 0$ over $\ghat$ are encoded by
$^L \ghat$-opers on $\pone$ with regular singularity at $0$ with
residue $\ol\varpi(-\nu-\rho)$ and singularity of order $2$ at
$\infty$ with $2$-residue $\pi(\chi)$.}}

\bigskip

Here it is important to note that the projective line $\pone$, on
which these $^L \ghat$-opers are defined, has a global coordinate $z$,
which is nothing but the {\em spectral parameter} of the affine Gaudin
model (in particular, $z$ is the spectral parameter in the
$L$-operator \eqref{eta L}).

\subsection{The case of $\su$}

What do these $^L \ghat$-opers look like? Consider the case when
$\ghat = {}^L \ghat = \su$. The Cartan subalgebra $^L \wt\h$ has a
basis consisting of
$$
h = \begin{pmatrix} 1 & 0 \\ 0 & -1 \end{pmatrix},
$$
${\mb 1}$, and $d$, and $^L \ol\h$ is spanned by $h$ and $d$ (recall
that $d$ acts as $\la \pa_\la$). Without loss of generality we may,
and will, consider only highest weights $\nu$ such that $\langle \nu,d
\rangle = 0$. They are determined by the numbers $k = \langle \nu,{\mb
1} \rangle$ and $\ell = \frac{1}{2} \langle \nu,h \rangle$. Then the
corresponding element of $^L \ol\h = \C h \oplus \C d$ is represented
by $\ell h + k d$.

We have two simple roots, $\al_1$ and $\al_0$, which in $^L \ol\h$ are
represented as $\al_1 = h$ and $\al_0 = -h$. We therefore divide the
set of points $\{ w_1,\ldots,w_m \}$ into two subsets:
$\{w^1_1,\ldots,w^1_{m_1}\}$ and $\{ w^0_1,\ldots,w^0_{m_0} \}$. The
connection \eqref{olnabla1} gives rise to the Miura $\su$-oper
\begin{equation}    \label{su oper}
\pa_z + \begin{pmatrix} u(z) & \la \\ 1 & -u(z) \end{pmatrix} +
\frac{k}{z} d,
\end{equation}
where
\begin{equation}    \label{ut1}
u(z) = - \chi - \frac{\ell}{z} + \sum_{j=1}^{m_1} \frac{1}{z-w^1_j} -
\sum_{j=1}^{m_0} \frac{1}{z-w^0_j}.
\end{equation}

Concretely, the statement of \lemref{bethe and mon} is realized as
follows. The connection in \eqref{su oper} has apparent singularities
at the points $z=w^1_j$ and $z=w^0_j$ (in addition to $z=0$ and
$z=\infty$), but the condition of \lemref{bethe and mon} is that these
singularities are ``fake''; that is, they can be removed by a gauge
transformation (by an element of $\ol{B}_+$). This imposes certain
conditions on the points $w^1_j$ and $w^0_j$, which are in fact the
Bethe Ansatz equations \eqref{one mod}.

Let us see more precisely how this works. This is in fact a local
question, which involves only the restrictions of the Miura oper
\eqref{su oper} to the formal discs around these points.

Suppose that we have a Miura $\su$-oper of the form
\begin{equation}    \label{su oper1}
\pa_z + \begin{pmatrix} u(z) & \la \\ 1 & -u(z) \end{pmatrix} +
x(z) d
\end{equation}
on the formal disc around a point $w$ such that
$$
u(z) = \sum_{n\geq -1} u_n (z-w)^n, \qquad x(z) = \sum_{n\geq
  0} x_n (z-w)^n.
$$
Suppose first that $u_{-1}=1$. This is the case when $u(z)$ is given
by formula \eqref{ut1} and $w=w^1_j$. Let us apply the gauge
transformation by
$$
\begin{pmatrix} 1 & -\frac{1}{z-w} \\ 0 & 1 \end{pmatrix}.
$$
Then we obtain the operator
$$
\pa_z + \frac{2}{z-w} \cdot \left. \left( u(z) - \frac{1}{z-w}
\right) \right|_{z=w} \begin{pmatrix} 0 & 1 \\ 0 & 0 \end{pmatrix}
+ \ldots,
$$
where the dots stand for the terms regular at $z=w$. Therefore this
operator is regular at $z=w$ if and only if the equation
\begin{equation}    \label{local bae1}
u_0 = 0
\end{equation}
is satisfied.

Next, consider the case when $u_{-1}=-1$, which corresponds to
$w=w^0_j$ in formula \eqref{ut1}. Applying the gauge
transformation by
$$
\begin{pmatrix} 1 & 0 \\ -\frac{1}{z-w} \la^{-1} & 1
\end{pmatrix},
$$
we obtain the operator whose expansion near $z=w$ has the
form
$$
\pa_z + \frac{2}{z-w} \cdot \left. \left( -u(z) -
\frac{1}{z-w} - \frac{x(z)}{2} \right) \right|_{z=w}
\begin{pmatrix} 0 & 0 \\ \la^{-1} & 0 \end{pmatrix} + \ldots,
$$
where the dots stand for terms regular at $z=w$. Therefore this
operator is regular at $z=w$ if and only if the equation
\begin{equation}    \label{local bae0}
u_0 + \frac{x_0}{2} = 0
\end{equation}
is satisfied.

Thus, we obtain that the connection \eqref{su oper} with $u_{-1}=1$
(resp., $u_{-1}=-1$) is gauge equivalent to a connection regular at $z=w$
if the equation \eqref{local bae1} (resp., \eqref{local bae0}) is
satisfied. Conversely, one shows, following the argument of
\cite{F:flag}, Lemma 2.10, that if this equation is not satisfied,
then the connection \eqref{su oper} has non-trivial monodromy around
$w$. Therefore it is impossible to make it regular by a gauge
transformation.

\bigskip

Let us apply this in the case of our Miura oper \eqref{su oper}. Then
equations \eqref{local bae1} for $w=w^1_j$ become
\begin{equation}    \label{bethe1}
\frac{\ell}{w_j^1} -
\sum_{s \neq j} \frac{1}{w_j^1-w^1_s} + \sum_{s=1}^{m_0}
\frac{1}{w^1_j-w^0_s} + \chi = 0, \qquad j=1,\ldots,m_1,
\end{equation}
and equations \eqref{local bae0} for $w=w^0_j$ become
\begin{equation}    \label{bethe2}
\frac{\frac{k}{2}-\ell}{w_j^0} + \sum_{s=1}^{m_1}
\frac{1}{w_j^0-w^1_s} - \sum_{s \neq j} \frac{1}{w^0_j-w^0_s} - \chi =
0, \qquad j=1,\ldots,m_0.
\end{equation}
These are precisely the Bethe Ansatz equations \eqref{one mod}.

\subsection{Second order operator}    \label{sec ord}

Now we wish to rewrite \eqref{su oper} as a second order differential
operator. To this end, we apply the gauge transformation by
\begin{equation}    \label{gauge eq}
\begin{pmatrix} 1 & -u(z) \\ 0 & 1 \end{pmatrix}.
\end{equation}
Then we obtain the following operator:
\begin{equation}    \label{oper bethe1}
\pa_z + \begin{pmatrix} 0 & v(z) + \la \\ 1 & 0 \end{pmatrix} +
\frac{k}{z} d,
\end{equation}
where
\begin{align} \notag
v(z) &= u(z)^2 + \pa_z u(z) \\ \label{form of conn1} &= \chi^2 +
\left( 2 \ell \chi + \sum_{j=1}^{m_0} \frac{k}{w^0_j} \right)
\frac{1}{z} + \frac{\ell(\ell+1)}{z^2} + \sum_{j=1}^{m_0}
\frac{2}{(z-w^0_j)^2} + \sum_{j=1}^{m_0} \frac{k}{w_j^0}
\frac{1}{z-w^0_j}.
\end{align}
In deriving this formula, we have used the Bethe Ansatz equations
\eqref{bethe1}--\eqref{bethe2}. They imply in particular that $v(z)$
has no singularities at $w^1_j, j=1,\ldots,m_1$.

In addition, the Bethe Ansatz equations \eqref{bethe2} imply that the
coefficients $v_{j,n}$ in the expansion of $v(z)$ at $z=w^0_j$,
$$
v(z) \sim \frac{2}{(z-w^0_j)^2} + \sum_{n \geq -1} v_{j,n}
(z-w^0_j)^n,
$$
satisfy the algebraic equation
\begin{equation}    \label{third order}
\frac{1}{4} \left( \frac{k}{w_j^0} \right)^3 - \frac{k}{w_j^0} v_{j,0}
+ v_{j,1} = 0, \qquad j=1,\ldots,m_0.
\end{equation}
Note also that $v_{j,-1} = k/w_j^0$.

For instance, if $k=0$, then we have $v_{j,-1}=0$ and the equation
\eqref{third order} reduces to $v_{j,1}=0$ for all
$j=1,\ldots,m_0$. Hence we obtain that in this case the condition is
that the coefficients in front of $(z-w^0_j)^{-1}$ and $(z-w^0_j)$ in
the expansion of $v(z)$ in Laurent series in $(z-w^0_j)$ should vanish
for all $j=1,\ldots,m_0$.

These equations guarantee that on the formal neighborhood of each of
the points $w_j^0, j=1,\ldots,m_0$, the $\su$-oper \eqref{oper
bethe1}, viewed as a connection on the trivial $\ol{G}$-bundle on the
punctured disc around this point, is gauge equivalent to a connection
without singularity. In other words, this $\su$-oper has no monodromy
around the points $w_j^0, j=1,\ldots,m_0$.

This is in fact a local condition on this oper which involves only its
restrictions to the formal discs around these points. Let us analyze
this local condition in more detail.

\bigskip

Suppose that we have an $\su$-oper of the form
\begin{equation}    \label{local oper bethe}
\pa_z + \begin{pmatrix} 0 & v(z) + \la \\ 1 & 0 \end{pmatrix} +
x(z) d,
\end{equation}
on the formal disc around a point $w$ such that
\begin{equation}    \label{vx}
v(z) = \frac{2}{(z-w)^2} + \sum_{n \geq -1} v_n (z-w)^n, \qquad x(z) =
\sum_{n\geq 0} x_n (z-w)^n.
\end{equation}
What is the condition that the connection \eqref{oper bethe1} on the
trivial $\ol{G}$-bundle on the punctured disc around the point $w$
(where $\ol{G} = \C^\times \ltimes G\pparl$) is gauge equivalent to a
connection without singularity?

This connection preserves the trivial $\ol{G}_-$-bundle, where
$\ol{G}_- = \C^\times \ltimes G[[\la]]$. Taking the quotient modulo
the subgroup $\C^\times \ltimes G^{(1)}$, where $G^{(1)}$ is the first
congruence subgroup of $G[[\la]]$, we obtain the connection
\begin{equation}    \label{local oper bethe1}
\pa_z + \begin{pmatrix} 0 & v(z) \\ 1 & 0 \end{pmatrix}
\end{equation}
on the trivial $G$-bundle, which we view as an $\sw_2$-oper. If the
connection \eqref{local oper bethe} has no monodromy around $z=w$,
then the same is true for the connection \eqref{local oper
bethe1}. This implies that the following equation on the coefficients
in the expansion of $v(z)$ at $z=w$ must hold:
\begin{equation}    \label{eq of no mon1}
\frac{1}{4} v_{-1}^3 - v_{0} v_{-1} + v_{1} = 0
\end{equation}
(see, e.g., \cite{FFT}, Sect. 4.4). It is easy to check that this
equation is satisfied if and only if $v(z)$ may be represented as the
Miura transformation
\begin{equation}    \label{miura again}
v(z) = u(z)^2 - \pa_z u(z),
\end{equation}
where
$$
u(z) = -\frac{1}{z-w} + \sum_{n \geq 0} u_n (z-w)^n.
$$

Thus, our connection \eqref{local oper bethe} is gauge equivalent,
under the gauge transformation \eqref{gauge eq}, to a connection of
the form \eqref{su oper1} with $u_{-1}=-1$. As explained in the
previous section, the latter is gauge equivalent to a connection
without singularity at $z=w$ if and only if equation \eqref{local
bae0} is satisfied. Substituting this into formula \eqref{miura
again}, we find that
\begin{equation}    \label{local bae2}
v_{-1} = x_0.
\end{equation}
As the result, we obtain the equation
\begin{equation}    \label{eq of no mon}
\frac{1}{4} x_0^3 - v_{0} x_0 + v_{1} = 0.
\end{equation}

Thus, we obtain that the condition for the operator \eqref{local oper
bethe} to be monodromy free is given by two equations: \eqref{eq of no
mon} and \eqref{local bae2}. Applying this to the $\su$-oper
\eqref{oper bethe1} on $\pone$ with $w=w^0_j$, we obtain the equations
\eqref{third order} and $v_{j,{-1}} = k/w_j^0$.

Recall that we obtained equation \eqref{third order} originally under
the assumption that $v(z)$ may be represented globally, on the entire
$\pone$, as the Miura transformation of \eqref{su oper}. But now we
see that we do not need to make this assumption. In fact, the local
condition of triviality of the monodromy around a point $z=w$ of any
$\su$-oper \eqref{local oper bethe} with $v(z)$ and $x(z)$ given by
formula \eqref{vx} is expressed by the equations \eqref{eq of no mon1}
and \eqref{local bae2}, which give rise to the equation \eqref{third
order} in the case of the oper \eqref{oper bethe1}.

\bigskip

Applying the (formal) gauge transformation by $z^{k d}$ to \eqref{oper
bethe1}, we obtain the operator
$$
\pa_z + \begin{pmatrix} 0 & v(z) + \la z^k \\
1 & 0 \end{pmatrix},
$$
which may be rewritten as the following second order operator with
spectral parameter
\begin{equation}    \label{final form rs}
\pa_z^2 - \chi^2 - \left( 2 \ell \chi + \sum_{j=1}^{m_0}
\frac{k}{w^0_j} \right) \frac{1}{z} - \frac{\ell(\ell+1)}{z^2} -
\sum_{j=1}^{m_0} \frac{2}{(z-w^0_j)^2} - \sum_{j=1}^{m_0}
\frac{k}{w_j^0} \frac{1}{z-w^0_j} - \la z^k.
\end{equation}
Now we may treat $\la$ as a complex parameter. The equations
\eqref{third order} then imply that the operator \eqref{final form rs}
has no monodromy around the points $w^0_j, j=1,\ldots,m_0$, for all
values of $\lambda$.

Let us summarize: a connection $\pa_z + u(z)$ given by formula
\eqref{ut} with the numbers $w^0_j, w^1_j$ satisfying the Bethe Ansatz
equations \eqref{bethe1}--\eqref{bethe2} gives rise to an $\su$-oper
\eqref{oper bethe1} on $\pone$ (equivalently, the second order
operator \eqref{final form rs}) such that the coefficients in its
expansion at $w^0_j, j=1,\ldots,m_0$, satisfy the equations
\eqref{third order}. As we have seen above, these equations are the
``no monodromy'' conditions which may be viewed as analogues of the
Bethe Ansatz equations \eqref{bethe2}. Our proposal for $\ghat = \su$
is that for generic $\chi$ and $\nu$ the common eigenvalues of the
quantum affine Gaudin Hamiltonians on the Verma module $M_\nu$ over
$\su$ are encoded by $\su$-opers of this form.

A solution of the Bethe Ansatz equations
\eqref{bethe1}--\eqref{bethe2} gives us a solution of the
equations \eqref{third order}, but we do not expect the converse to be
true in general unless we impose additional conditions on $v(z)$. For
instance, suppose that $v(z)$ is of the form \eqref{form of conn1} and
that in addition to having no monodromy around the points $w^0_j$,
which is ensured by the equations \eqref{third order}, the operator
\eqref{oper bethe1} has monodromy $\pm 1$ around $z=0$ (this means, in
particular, that $\ell$ is a non-negative half-integer). Then we
believe that such $v(z)$ does come from a solution of the Bethe Ansatz
equations \eqref{bethe1}--\eqref{bethe2}. This is a conjecture, which
is closely related to \conjref{bethe int} below.

Note that the operator \eqref{oper bethe1} also gives rise to a
parabolic Miura $\su$-oper corresponding to the $\sw_2$-oper
$$
\pa_z + \begin{pmatrix} 0 & v(z) \\ 1 & 0 \end{pmatrix}
$$
and the connection $\pa_z + \frac{k}{z}$ in
$\on{Conn}_{\pone}(\Omega^2)$, in the sense of the definition of
\secref{par miura}. The singular points $w^0_j, j=1,\ldots,m_0$, are
the points where the two reductions, $\F_{\ol{B}_+}$ and
$\F_{\ol{G}_-}$, of this parabolic Miura oper are in the relative
position $\ol{s_0} \in W \bs W_{\on{aff}}$, where $s_0 \in
W_{\on{aff}}$ is the simple reflection corresponding to the simple
root $\al_0$. This is the generic situation. For special values of
$\nu$ and $\chi$ the positions of some of the poles $w^0_j$ may
coalesce and these opers may degenerate further. In other words, there
may appear points on $\pone$ at which the relative positions of the
two reductions are given by elements of $W \bs W_{\on{aff}}$ other
than $\ol{s}_0$. If the relative position at a point $w$ is given by
the element $\ol{s_0(s_1 s_0)^n} \in W \bs W_{\on{aff}}, n \in \Z$,
then the expansion of the function $v(z)$ appearing in \eqref{form of
conn1} around $w$ will have the form
$$
v(w) \sim \frac{(2n+1)(2n+2)}{(z-w)^2} + \ldots .
$$
The coefficients in the expansion of this function at $w$ will
also have to satisfy an analogue of equation \eqref{third order},
which means that our $\su$-oper has no monodromy at $w$.

In addition, some of the points $w^0_j$ may tend to $0$. This
happens when the Verma module $M_\nu$ becomes reducible. The common
eigenvalues encoded be these degenerations of the opers \eqref{oper
bethe1} should correspond to the eigenvectors which appear in a
submodule of $M_\nu$.

Finally, note that we may rewrite \eqref{oper bethe1} as the following
partial differential operator:
$$
\left(\pa_z + \frac{k}{z} \la \pa_\la \right)^2 - v(z) -
\la.
$$
This formula underscores the fact that we are dealing here with two
independent variables: $z$ and $\lambda$. In our analysis we have
treated $z$ as a variable, but $\lambda$ as the parameter of the loop
algebra $\g\pparl$. Thus, we have considered an ordinary differential
operator \eqref{oper bethe1}, in $z$ only, but with values in an
infinite-dimensional Lie algebra; namely, the Lie algebra $\C d \oplus
\g\pparl$. However, considering it as a partial differential
operator, in both $z$ and $\lambda$, might also be useful in some
applications.

\subsection{Integrable representations}

Now we consider the common eigenvalues of the affine Gaudin
Hamiltonians on an integrable representation of $\ghat$. Recall that
for each integral dominant weight $\nu$, i.e., such that $\langle
\nu,\chal_i \rangle \in \Z_+, i \in \wt{I}$, so in particular $k =
\langle \nu,{\mb 1} \rangle \in \Z_+$, we have an irreducible
integrable representation $V_\nu$, which is the irreducible quotient
of the Verma module $M_\nu$. By analogy with the finite-dimensional
case (see \cite{FFT}, \thmref{spectra fd},(2) and the conjecture after
it), we propose the following description of the common eigenvalues of
the affine Gaudin Hamiltonians on the space of $t\h[[t]]$-invariant
vectors in integrable representations of $\ghat$. Note that this space
is a module over the corresponding generalized parafermion algebra,
which appears in the space of states of the $(\ghat,\wh\h)$ coset
model.

\begin{conj}    \label{bethe int}
The set of common eigenvalues of the affine Gaudin Hamiltonians on the
space of $t\h[[t]]$-invariant vectors in an integrable module
$V_\nu$ over $\ghat$ (counted without multiplicities) is in bijection
with the set of $^L \ghat$-opers on $\pone$ which are $\nu$-regular at
$0$, and have singularity of order $2$ at $\infty$ with the
$2$-residue $\pi(\chi)$.
\end{conj}

Note that all of these opers necessarily have {\em trivial monodromy}
on $\pone$.

\medskip

If the projection of the Bethe vector \eqref{genbv} onto $V_\nu$ is
non-zero, then the corresponding common eigenvalues of the affine
Gaudin Hamiltonians are recorded by the $^L \ghat$-oper obtained by
applying the Miura transformation to the Cartan connection \eqref{new
miura}. It then satisfies the conditions listed in \conjref{bethe
int}. This provides the first positive test for this conjecture.

It is natural to ask whether the Bethe Ansatz is complete in this
model for generic values of $\chi$. As in the finite-dimensional case
(see \cite{FFT}), we expect that it is complete, which gives us hope
that the Bethe eigenvectors may furnish an explicit family of bases
for all integrable representations of affine Kac--Moody algebras
(parameterized by $\chi$).

If the Bethe Ansatz is incomplete, then there are eigenvectors other
than the Bethe vectors introduced above. In this case we conjecture
that the common eigenvalues of the affine Gaudin Hamiltonians would
still be encoded by the $^L \ghat$-opers as in \conjref{bethe
int}. However, the eigenvalues corresponding to the non-Bethe
eigenvectors would be encoded by opers that cannot be expressed as the
Miura transformation of a connection of the form \eqref{olnabla1}. The
structure of these opers in the finite-dimensional case is explained
in \cite{F:faro}, Sect. 5.5. We expect that the same pattern holds in
the affine case as well. In particular, we expect that the opers
corresponding to the non-Bethe eigenvectors may still be expressed as
the Miura transformation, but of Cartan connections of the more
general form
\begin{equation}    \label{olnabla2}
\pa_z - \chi + \frac{\rho-y(\nu+\rho)}{z} + \sum_{j=1}^m
\frac{\rho-y_j(\rho)}{z-w_j},
\end{equation}
where $\rho \in {}^L \ol\h$ is defined as in \secref{par miura}, and
$y$ and $y_j, j=1,\ldots,m$, could be arbitrary elements of the affine
Weyl group corresponding to $^L \g$. They measure the
relative position between the two reductions of the oper bundle, the
oper reduction ${\mc F}_{B_+}$ and the horizontal reduction ${\mc
F}_{B_-}$ at the points $z=0$ and $z=w_j$, respectively.

Note that if $y=s_i$ is the simple reflection corresponding to the
$i$th simple root, then $\rho-y(\rho)=\al_i$. Thus, if $y=1$ and all
the $y_j$'s are simple reflections, then \eqref{olnabla2} becomes
\eqref{olnabla1}. We expect that this is the generic
situation. However, for special values of parameters some of the
common eigenvalues may be encoded by the opers corresponding the
Miura transformation of the more general connections
\eqref{olnabla2}. For instance, for $\g=\sw_2$ the weights
$\rho-y_j(\rho)$ have the form $\nu = 2n_j$, where $n_j$ is an
integer. The corresponding operator \eqref{final form rs} will then
have the leading term $-n_j(n_j+1)/(t-w_j)^2$ (formula \eqref{final
form rs} corresponds to the case when all $n_j=1$).

\subsection{Generalization to other Kac--Moody algebras}    \label{gen
  km}

It is natural to ask which of the structures discussed above may be
generalized to an arbitrary symmetrizable Kac--Moody algebra $\g$. Let
us fix $\chi \in \h^*$ such that
$$
(\al,\chi) \neq 0, \qquad \forall \al \in \De_+,
$$
where $\De_+$ is the set of positive roots of $\g$. The following
analogues of the DMT Hamiltonians \eqref{eq:T Ham} are well-defined on
any highest weight module over $\g$:
\begin{equation}    \label{DMT km}
T_\gamma(\chi) = \sum_{\alpha\in\Delta_+}
\frac{(\alpha,\gamma)}{(\alpha,\chi)} \; f^{(i)}_{\al} e^{(i)}_{\al},
\qquad \ga \in \h^*,
\end{equation}
where $\{ e^{(i)}_{\al} \}$ and $\{ f^{(i)}_{\al} \}$ form bases of
the root subspaces $\g_\al$ and $\g_{-\al}$, respectively, which are
dual with respect to a non-degenerate invariant inner product on
$\g$. These operators coincide with the connection operators of a
connection studied in \cite{FMTV}. The flatness of this connection,
established in \cite{FMTV} for an arbitrary symmetrizable Kac--Moody
algebra, implies that these operators mutually commute for different
$\ga$ (and fixed $\chi$).\footnote{This commutativity may probably
also be proved by generalizing the argument used in \cite{TL}. We
thank V. Toledano Laredo for a useful discussion of this point.}

However, note that for the affine Kac--Moody algebras these
Hamiltonians are different from the Hamiltonians \eqref{aff T
Ham}. While in formula \eqref{DMT km} the summation is over the set of
all positive roots of $\ghat$, in formula \eqref{aff T Ham} it is over
the set of positive real roots (and $\chi,\ga$ are elements of the
dual space to the Cartan subalgebra of the finite-dimensional Lie
algebra $\g$, not the extended Cartan subalgebra of $\gtil$). The
operators \eqref{aff T Ham} may be obtained from \eqref{DMT km} if we
choose $\ga$ to be orthogonal to the imaginary root $\delta$ and then
take the limit when $\chi$ becomes orthogonal to $\delta$ as well.

Consider the Bethe vectors $\phi(w_1^{i_1},\ldots,w_m^{i_m})$ given by
formula \eqref{bethe for one mod}, where $i_j \in I$, the set of
vertices of the Dynkin diagram of $\g$. We then have an analogue of
the statement of \conjref{phi eig} for the operators \eqref{DMT
km}. It may be derived from the results of \cite{FMTV} by considering
the asymptotics of the solutions of the dynamical KZ equations
constructed in \cite{FMTV} near the critical level, along the lines of
\cite{RV}.\footnote{We thank A. Varchenko for pointing this out.} By
taking the limit described in the previous paragraph, one can then
prove \conjref{phi eig}.

On the other hand, as explained in \cite{F:flag}, the notions of
opers, Cartan connections and Miura transformation generalize to
arbitrary symmetrizable Kac--Moody algebras in a straightforward
fashion. Moreover, to each solution of the Bethe Ansatz equations
\eqref{one mod} we may attach a Cartan connection by formula
\eqref{olnabla1}. Then its Miura transformation is an $^L \g$-oper
which satisfies the same properties as in the finite-dimensional and
affine cases. Therefore it is tempting to assume that analogues of our
conjectures relating $^L \g$-opers to the spectra of commuting
Hamiltonians (in particular, an analogue of \conjref{bethe int}) also
hold in this general context.

However, for Kac--Moody algebras other than finite-dimensional or
affine, it is not clear how to construct other commuting quantum
Hamiltonians. So at present it is not clear what is the commutative
algebra whose spectra are described by these $^L \g$-opers.

\section{Quantum KdV system}    \label{kdv}

In this section we consider another example of the shift of argument
affine Gaudin model, which corresponds to the $\ghat$-KdV
hierarchies. In this case the shift $\chi$ is the nilpotent element of
$\g \subset \ghat$ corresponding to the maximal root. The
corresponding shift of argument algebra is not so interesting by
itself, so we look instead at its image in the Hamiltonian reduction
with respect to the centralizer of $\chi$, which is
$\n\ppart$. Classically, this Hamiltonian reduction has been
introduced by Drinfeld and Sokolov \cite{DS}, and the corresponding
commuting local classical Hamiltonians are the Hamiltonians of the
$\ghat$-KdV hierarchy, as explained in \cite{DS}. These Hamiltonians
generate a commutative subalgebra in the Poisson algebra obtained by
the Drinfeld--Sokolov Hamiltonian reduction from the algebra of
functions on $\ghat^*_1$. This Poisson algebra is known as the
classical ${\mc W}$-{\em algebra}. Its quantization is the quantum
${\mc W}$-algebra ${\mc W}(\g)$ introduced in \cite{FL} for $\g=\sw_n$
and in \cite{FF:gd} for general $\g$ (see also
\cite{FF:ds,FF:laws,FB}). Thus, the quantum Gaudin Hamiltonians in
this case are commuting elements of ${\mc W}(\g)$. They are known as
the local quantum KdV Hamiltonians \cite{FF:laws}.

As in the case of regular semi-simple shift of argument, discussed in
the previous section, there are actually two types of Hamiltonians:
local and non-local. They are obtained by expanding invariant
functions of the monodromy at $z=\infty$ and $z=0$, respectively. The
existence of the local quantum Hamiltonians was established in
\cite{FF:toda,FF:laws}, though explicit formulas for them are unknown
in general. The non-local quantum Hamiltonians have been constructed
in \cite{BLZ1}. Their spectra on irreducible highest weight modules
over the Virasoro algebra have been studied in the series of papers
\cite{BLZ1}--\cite{BLZ}. In particular, in \cite{DT,BLZ4,BLZ} a
remarkable conjecture was made relating these spectra to certain
Schr\"odinger operators. In this section we show that this conjecture
fits in perfectly as a special case of our general proposal that the
spectra of quantum Hamiltonians in the affine Gaudin model associated
to $\ghat$ are given by opers for the Langlands dual affine Lie
algebra $^L\ghat$.

\subsection{The classical model}    \label{class kdv}

The phase space of the classical model is defined in the same way as
in \secref{aff rs}. It consists of the operators\footnote{In what
follows, in order to simplify the exposition, we will assume that
$\ghat$ is untwisted. But our discussion applies to twisted affine
algebras as well.}
\begin{equation}    \label{eta L1}
L = \pa_t + A(t) + \chi z, \qquad A(t) \in \g\ppart
\end{equation}
where now
$$
\chi = e_{\theta} \subset \n_+ \subset \g \subset \g\ppart
$$
is a nilpotent element corresponding to the maximal root of $\g$.
The centralizer of $\chi$ is therefore the Lie algebra $\n_+\ppart$,
As explained in \secref{class fd}, instead of constructing a
Poisson commutative subalgebra in the Poisson algebra of functions on
\eqref{eta L1}, we may construct a commutative subalgebra in its
Hamiltonian reduction with respect to $\n_+\ppart$ and its
character. We choose the character $\psi$ which corresponds to the
element
\begin{equation}    \label{olp}
\ol{p}_{-1} = \sum_{i=1}^\ell f_i \; \in \; \n_- \subset \g\ppart.
\end{equation}
The resulting Hamiltonian reduction is the {\em
Drinfeld--Sokolov reduction} \cite{DS}.

The reduced phase space of the Drinfeld--Sokolov reduction is
therefore the quotient of the space $\wt{\mc M}(\ghat)$ of operators
of the form
\begin{equation}    \label{reduced}
L' = \pa_t + p_{-1} + {\mb v}(t), \qquad {\mb v}(t) \in \bb_+\ppart,
\end{equation}
by the gauge action of the group $N_+\ppart$. Here
$$
p_{-1} = \sum_{i \in \wt{I}} f_i = \sum_{i=1}^\ell f_i +
e_{\theta} z
$$
is a ``principal nilpotent element'' of $\g\zpart$. According
to \cite{DS}, the action of $N_+\ppart$ on $\wt{\mc M}(\ghat)$ is
free. The resulting quotient space\footnote{Note that the action of
$N_+\ppart$, and hence the quotient, do not depend on whether we
include $e_{\theta} z$ in $p_{-1}$ or not.}
$$
{\mc M}(\ghat) = \wt{\mc M}(\ghat)/N_+\ppart
$$
is isomorphic to the space of $\g$-opers on the punctured disc
$D^\times = \on{Spec} \C\ppart$. The Poisson algebra of local
functionals on ${\mc M}(\ghat)$ is known as the classical ${\mc
W}$-algebra ${\mc W}(\g)$.

For example, ${\mc M}(\su)$ is the quotient of operators of the form
$$
\pa_t + \begin{pmatrix} a(t) & b(t) + z \\
  1 & -a(t) \end{pmatrix}, \qquad a(t),b(t) \in \C\ppart,
$$
by the upper triangular gauge transformations depending on $t$ (but
not on $z$). It is easy to see that each orbit contains a unique
operator of the form
$$
\pa_t + \begin{pmatrix} 0 & v(t) + z \\
  1 & 0 \end{pmatrix}, \qquad v(t) \in \C\ppart,
$$
and hence we may identify ${\mc M}(\su)$ with the space of such
operators, or, equivalently, with the space of second order
differential operators with spectral parameter
$$
\pa_t^2 - v(t) - z, \qquad v(t) \in \C\ppart.
$$

Similarly, the space ${\mc M}(\wh\sw_n)$ may be identified with the
space of $n$th order differential operators with spectral parameter
$$
(-\pa_t)^n - v_1(t) (-\pa_t)^{n-2} - \ldots + v_{n-2}(t) \pa_t -
v_{n-1}(t) - z.
$$
For other Lie algebras there also exist canonical representatives for
the orbits of the gauge action of $N_+\ppart$ on the space of
operators of the form \eqref{reduced} (see, e.g.,
\cite{F:faro}). Recall the element $\ol{p}_{-1} \in \n_-$ given by
formula \eqref{olp}. There exists a unique element of $\n_+$ of the
form
$$
\ol{p}_1 = \sum_{i=1}^\ell c_i e_i, \qquad c_i \in \C,
$$
such that $\ol{p}_1, \ol{p}_{-1}$, and $\ol{p}_0 =
[\ol{p}_1,\ol{p}_{-1}]$ form an $\sw_2$ triple. The element
$\frac{1}{2} \ol{p}_0 \in \h$ then defines the principal grading on
$\g$ such that $\deg \ol{p}_1 = 1, \deg \ol{p}_{-1} = -1$. Let
$$
V_{\can} = \bigoplus_{i \in E} V_{\can,i}
$$
be the space of $\on{ad}
\ol{p}_1$-invariants in $\n$, decomposed according to the principal
grading. Here
$$
E = \{ d_1,\ldots,d_\ell \}
$$
is the set of exponents of $\g$. Then $\ol{p}_1$ spans
$V_{\on{can},1}$. Choose a linear generator $\ol{p}_j$ of
$V_{\can,d_j}$ (if the multiplicity of $d_j$ is greater than one,
which happens only in the case $\ghat=D^{(1)}_{2n}, d_j=2n$, then we
choose linearly independent vectors in $V_{\on{can},d_j}$). The
following result is due to Drinfeld and Sokolov \cite{DS} (the proof
is reproduced in \cite{F:flag}, Lemma 2.1).

\begin{lem} \label{free}
The gauge action of $N_+\ppart$ on the space $\wt{\mc M}(\ghat)$ is
free, and each gauge equivalence class contains a unique operator of
the form $\pa_t + p_{-1} + {\mathbf v}(t)$, where ${\mathbf
v}(t) \in V_{\can}\ppart$, so that we can write
\begin{equation} \label{coeff fun}
{\mathbf v}(t) = \sum_{j=1}^\ell v_j(t) \cdot p_j, \qquad v_j(t) \in
\C\ppart.
\end{equation}
\end{lem}

Thus, each point of the reduced phase space ${\mc M}(\ghat)$ of the
Drinfeld--Sokolov reduction is canonically represented by an operator
$\pa_t + p_{-1} + {\mathbf v}(t)$, where ${\mathbf v}(t)$ is of the
form \eqref{coeff fun}.

\subsection{Classical Hamiltonians}    \label{class ham kdv}

The Poisson commuting Hamiltonians on the reduced phase space are
constructed in the same way as in the case of regular semi-simple
$\chi$: by expanding the monodromy of the operators \eqref{reduced}
with respect to the spectral parameter $z$. Namely, let $M(z) \in G$
be the monodromy matrix of the operator \eqref{eta L1}. For any
invariant function $\varphi$ on the group $G$ the corresponding
function $H_\varphi(z) = \varphi(M(z))$ on $\wt{\mc M}(\ghat)$ is
invariant under the gauge action of $N_+\ppart$ and hence gives rise
to a well-defined function on the quotient ${\mc M}(\ghat)$. We obtain
Poisson commuting Hamiltonians on ${\mc M}(\ghat)$ by expanding
$H_\varphi(z)$ in $z$ and $z^{-1}$.

The asymptotic expansion at $z=\infty$ yields the local Hamiltonians,
which generate the $\ghat$-KdV hierarchy of commuting Hamiltonian
flows on the Poisson manifold ${\mc M}(\ghat)$. These Hamiltonians
have the form
$$
H_n = \int P_n(v_j(t),v'_j(t),\ldots) dt, \qquad n = d_i + Nh, N \in
\Z_+,
$$
where $d_i \in E$ is an exponent of $\g$ and $h$ is the Coxeter
number. The integrand $P_n$ is a differential polynomial of degree
$n+1$, where we set $\deg v^{(m)}_j = d_j+m+1$.

On the other hand, the $z$-expansion at $z=0$ yields Poisson commuting
non-local Hamiltonians. It is instructive to consider the $z$-linear
term in the expansion. For that we apply \lemref{linear term} with
$\ep=z$,
$$
A^{(0)}(t) = \sum_{i=1}^\ell f_i + v(t), \qquad v(t) \in
\bb_+\ppart,
$$
and $A^{(1)}(t) = e_{\theta}$. We
obtain that the $z$-linear term in the expansion of $\varphi(M(z))$,
where $\varphi \in (\on{Fun} G)^G$, is given by the formula
\begin{equation}    \label{H 1}
H^{(1)}_\varphi = \int \left( e_{\theta},\Psi_\varphi(t)
\right) dt.
\end{equation}
This is the lowest weight component of the single-valued solution
$\Psi_\varphi(t) \in \g$ of the equation
\begin{equation}    \label{kdv eq}
\left( \pa_t + \sum_{i=1}^\ell f_i + v(t) \right) \Psi(t) = 0,
\end{equation}
with values in the adjoint representation of $\g$, corresponding to
the invariant function $\varphi$.

The commutativity \eqref{comm ham} of the Hamiltonians $H_{\varphi}$
then implies that
$$
\{ H^{(1)}_{\varphi},H^{(1)}_{\psi} \} = 0
$$
for all $\varphi,\psi \in (\on{Fun} G)^G$.

Higher order terms in the $z$-expansion of $\varphi(M(z))$ may be
expressed as multiple integrals of the lowest weight components of
solutions of the equation \eqref{kdv eq}.

For example, if $\g=\sw_2$, the equation \eqref{kdv eq} in the adjoint
representation is equivalent to the following scalar third order
differential equation:
\begin{equation}    \label{third}
(\pa_t^3 - 4 v \pa_t - 2 v') \Psi = 0
\end{equation}
(its solutions are products of solutions of the second order equation
$(\pa_t^2 - v) \psi = 0$). It has a unique (up to a scalar)
single-valued solution $\Psi_\varphi(t)$ (here we may choose as
$\varphi$ the trace function on $SL_2$, as it generates the algebra of
invariant functions on $SL_2$), and $H^{(1)}_\varphi = \int
\Psi_\varphi(t) dt$.

\begin{remark}    \label{miura real}
A convenient way to compute the higher order terms in the
$z$-expansion of $\varphi(M(z))$ is to realize the variables of the
KdV hierarchy in terms of the variables of the modified KdV (mKdV)
hierarchy. This provides a kind of ``free field realization'', also
known as the Miura transformation, for the commuting Hamiltonians.

Consider the space $\ol{\mc M}(\ghat)$ of operators of the form
\begin{equation}    \label{u t}
\pa_t + p_{-1} + {\mb u}(t), \qquad {\mb u}(t) \in \h\ppart.
\end{equation}
The natural map $\ol{\mc M}(\ghat) \to {\mc M}(\ghat)$ given by the
composition of the inclusion $\ol{\mc M}(\ghat) \to \wt{\mc M}(\ghat)$
and the projection $\wt{\mc M}(\ghat) \to {\mc M}(\ghat)$ is the Miura
transformation.\footnote{We have already encountered Miura
transformation above, but in a totally different context: namely, it
appeared in our study of the eigenvalues of the {\em quantum}
Hamiltonians. Here it appears in the study of the phase space on which
the {\em classical} Hamiltonians are defined.}
It is a Poisson map with respect to the Heisenberg--Poisson structure
on $\ol{\mc M}(\ghat)$.

Because the operator \eqref{u t} has such a simple structure, it is
easier to compute the monodromy matrix $M(z)$ for it (and hence for
$\varphi(M(z))$) than for the operators of the form \eqref{reduced} or
\eqref{coeff fun}. One can then write down more explicit integral
formulas for the coefficients in the $z$-expansion in $\varphi(M(z))$.
For example, in the case when $\g=\sw_2$ the operator \eqref{u t} has
the form
$$
\pa_t + \begin{pmatrix} u(t) & z \\ 1 & -u(t) \end{pmatrix}.
$$
The coefficients in the $z$-expansion of the trace of the monodromy of
this operator are written down explicitly in \cite{BLZ1}. They are
given by multiple integrals of $\exp(\pm 2\phi(t))$, where $\phi(t)$
is the anti-derivative of $u(t)$, that is, $u(t) = \phi'(t)$. These
are the classical analogues of the two screening currents discussed at
the end of \secref{affine shift quantum}.

Similar formulas may be obtained for other simple Lie algebras. \qed
\end{remark}

\subsection{Quantum Hamiltonians}    \label{kdv quantum}

Now we consider the quantization of the Hamiltonians of the
$\ghat$-KdV system. The problem of quantization of the local
Hamiltonians has been considered, in particular, in
\cite{Zam,EY,KM,FF:toda,FF:laws}. It is closely related to the problem
of quantization of the mKdV Hamiltonians and integrals of motion of
affine Toda field theories (which are in turn related to deformations
of conformal field theories). It was proved in \cite{FF:toda,FF:laws}
that all classical local $\ghat$-KdV Hamiltonians may be
quantized. The quantum local Hamiltonians generate a commutative
subalgebra in the quantum ${\mc W}$-algebra ${\mc W}_\beta(\g)$, where
$\beta$ is the parameter related to the central charge. This is an
associative algebra quantizing the classical (Poisson) ${\mc
W}$-algebra of local functionals on ${\mc M}(\ghat)$. It is obtained
by the quantum Drinfeld--Sokolov reduction of a completion of the
enveloping algebra of $\ghat$ and depends on one complex parameter,
the central charge $c$, which is a function of the level of $\ghat$
(see \cite{FF:gd,FB}). For example, ${\mc W}_\beta(\sw_2)$ is nothing
but the completed enveloping algebra of the Virasoro algebra. We will
discuss the quantum local KdV Hamiltonians in more detail in
\secref{two points local} below.

The quantum non-local Hamiltonians were introduced in \cite{BLZ1} in
the case of $\g=\sw_2$ and in \cite{BHK} for $\g=\sw_3$ (this
construction may be generalized to other $\g$ \cite{FFS}).

To give an idea as to what these quantum Hamiltonians look like, we
discuss the quantization of the Hamiltonians $H^{(1)}_{\varphi}$
introduced above. Consider first the case when $\g=\sw_2$. In this
case $H^{(1)}_{\varphi}$ is given by the integral of a single-valued
solutions $\Psi_\varphi(t)$ of the equation \eqref{third}, which is
unique up to a scalar. The quantum analogues of solutions of this
equations are the primary fields $\Phi_{(1,3)}(w)$ of the Virasoro
algebra. There are actually three primary fields,
$\Phi^\pm_{(1,3)}(w)$ and $\Phi^0_{(1,3)}(w)$, acting between Verma
modules over the Virasoro algebra, but only one of them,
$\Phi^0_{(1,3)}(w)$, acts from a Verma module to itself. This field is
the quantization of the unique (up to a scalar) single-valued solution
$\Phi_\varphi(t)$ of the equation \eqref{third}. Therefore the
quantization of the Hamiltonian $H^{(1)}_\varphi$ is given by the
$0$th Fourier component of the primary field
$\Phi^0_{(1,3)}(w)$. (Here we multiply the primary field by
$w^{\Delta_{(1,3)}}$, where $\Delta_{(1,3)}$ is the conformal
dimension, so that $\Phi^0_{(1,3)}(w)$ is a formal power series in
integer powers of $w$.)

To construct the quantizations of other coefficients in the
$z$-expansion of $\varphi(M(z))$, we may use the free field
realization of the Virasoro algebra in terms of the Heisenberg
algebra, quantizing the Miura transformation map between the Poisson
algebras of functions on ${\mc M}(\su)$ and $\ol{\mc M}(\su)$
described in \remref{miura real}. Then we need to quantize the
coefficients in the $z$-expansion of the monodromy matrix of the
operator \eqref{u t}. The resulting quantum non-local Hamiltonians
will be given by multiple integrals of the screening currents
quantizing the classical screening currents $\exp(\pm 2\phi(t))$
introduced in \remref{miura real}. This is the approach taken in
\cite{BLZ1}.

More precisely, the integrals of the two screening currents satisfy
the $q$-Serre relations of $U_q(\su)$, where $q = \exp(\pi i \beta^2)$
(see \cite{FF:laws}). The non-local Hamiltonians correspond to
singular vectors in the Verma modules over $U_q(\su)$ of critical
level (see the discussion in \secref{affine shift quantum}). In
particular, the simplest of these quantum Hamiltonians,
$H^{(1)}_\varphi$, is given by the $q$-commutator of the two screening
charges.

Alternatively, we may use the primary fields $\Phi^\pm_{(1,3)}(w)$ and
$\Phi^0_{(1,3)}(w)$, and the fact that their integrals generate a
quantum deformation of the algebra $U(\wt{\sw}_{2,-})$, as discussed in
\secref{affine shift quantum}.

This construction generalizes to other Lie algebras. For each $\g$
there are primary fields $\Phi_{(1,\on{adj})}(w)$ for the
corresponding quantum ${\mc W}$-algebra ${\mc W}_\beta(\g)$, which are
the analogues of the primary fields $\Phi_{(1,3)}(w)$. The vector
space of this fields is isomorphic to $\g$, and it contains an
$\ell$-dimensional subspace (isomorphic to $\h \subset \g$) of those
primary fields $\Phi_{(1,\on{adj}),p}(w)$ which act from each Verma
module over ${\mc W}_\beta(\g)$ to itself. These are the quantum
analogues of the single-valued solutions $\Psi_\varphi(t)$ of the
equation \eqref{kdv eq}. The $0$th Fourier coefficients of
$\Phi_{(1,\on{adj}),p}(w), p=1,\ldots,\ell$ (multiplied by
$w^{\Delta_{(1,\on{adj})}}$) are the quantizations of the classical
Hamiltonians $H^{(1)}_\varphi$ given by formula \eqref{H 1}. Other
quantum non-local Hamiltonians may be expressed as multiple integrals
of the primary fields $\Phi_{(1,\on{adj})}(w)$. Alternatively, we may
use the free field realization of ${\mc W}_\beta(\g)$ and express then
as multiple integrals of the corresponding screening currents
\cite{BHK,FFS}. Note that the quantum non-local Hamiltonians commute
with each other and with the local ones. We will discuss this in more
detail in \cite{FF:next}.

The upshot of the above discussion is that we have a large commutative
algebra acting on the Verma modules (and other highest weight modules)
over the Virasoro algebra. This algebra is generated by the local and
non-local quantum KdV Hamiltonians. We now look for a parameterization
of the common eigenvalues of these Hamiltonians. Our experience with
the shift of argument models in which the shift $\chi$ is regular
semi-simple suggests that those should be parameterized by affine $^L
\ghat$-opers of special kind (discussed in \secref{spectra aff}).

In the case at hand we have a nilpotent shift of argument. In
order to understand better what kind of opers correspond to this
shift, we first consider the analogous question in the
finite-dimensional case.

\subsection{Nilpotent shift in the finite-dimensional case}

Recall from \cite{FFT} that the shift of argument Gaudin algebra
${\mc A}_\chi$ may be obtained as a quotient of the algebra
$\on{End}_{\ghat} \I_{1,\chi,\ka_c}$, where
$$
\I_{1,\chi} = \Ind_{t\g[[t]] \oplus \C{\mb 1}}^{\ghat}
\C_\chi.
$$
Here $\C_\chi$ is the one-dimensional module over $t\g[[t]] \oplus
\C {\mb 1}$, on which $t\g[[t]]$ acts via the character
$$
t\g[[t]] \to t\g[[t]]/t^2\g[[t]] \simeq \g
\overset\chi\longrightarrow \C
$$
and the central element ${\mb 1}$ act as $-h^\vee$. Let us choose
$\chi$ to be the linear functional on $\g$ taking a non-zero value on
$e_{\theta}$ and zero on all other root generators. The
centralizer of this $\chi$ is the Lie subalgebra $\n_- \subset \g$.
We wish to understand the corresponding reduced algebra ${\mc
A}_\chi^\psi$ (see \secref{quant fd}), where $\psi$ is a principal
character of $\n_-$.

Hence we perform the quantum Drinfeld--Sokolov reduction of the module
$\I_{1,\chi}$ with respect to the Lie algebra $\n_-\ppart$ and its
character
\begin{equation}    \label{charac}
f_{\al} \otimes t^n \mapsto \begin{cases} 1, & \al=\al_i, n=0, \\
0, & \on{otherwise} \end{cases}
\end{equation}
(see \cite{FF:lmp}).

As the result, we obtain a module over the center $Z(\ghat)$ at the
critical level, which we will denote by
$H_{\on{DS}}(\I_{1,\chi})$. Recall that $Z(\ghat)$ is a completion of
the polynomial algebra with generators $S_{i,(n)}, i=1,\ldots,\ell$;
$n \in \Z$, where we use the convention, in which the generating
currents are
$$
S_i(w) = \sum_{n \in \Z} S_{i,(n)} w^{-n-1}
$$
(the first of then, $S_1(z)$ is the unnormalized Segal--Sugawara
current), see, e.g., \cite{FFT}. It is easy to see that the elements
$S_{i,(n_i)}, i=1,\ldots,\ell-1; n_i > d_i$, as well as
$S_{\ell,(n_\ell)}, n_\ell > d_\ell+1$, act by $0$ on the generating
vector of $\I_{1,\chi}$ and hence on $H_{\on{DS}}(\I_{1,\chi})$. But
$S_{\ell,(d_\ell+1)}$ contains
a monomial which is the product of $f_i \otimes 1, i=1,\ldots,\ell$,
and $e_{\theta} \otimes t$, with a non-zero coefficient. Each factor
of this monomial acts by a non-zero scalar on the generating vector;
for $f_i \otimes 1, i=1,\ldots,\ell$, this follows from formula
\eqref{charac}, and for $e_{\theta} \otimes t$ this follows from the
definition of $\I_{1,\chi}$. Therefore the monomial itself, and hence
$S_{\ell,(d_\ell+1)}$, act by a non-zero scalar on $\I_{1,\chi}$ and
on $H_{\on{DS}}(\I_{1,\chi})$.

Recall (see, e.g., \cite{F:faro}) that we have canonical
representatives for $^L \g$-opers on the punctured disc $D^\times$ of
the form\footnote{Note that in \secref{class kdv} have encountered
$\g$-opers as points of the phase space ${\mc M}(\ghat)$ of the
classical $\ghat$-KdV hierarchy. On the other hand, now we are
discussing the Hamiltonians in the quantum shift of argument Gaudin
model associated to $\g$. The spectra of those Hamiltonians are
expressed by $^L \g$-opers. It is important to realize that these are
two different spaces of opers! For more on this, see \secref{two
ways}.}
$$
\pa_z + p_{-1} + {\mb v}(z),
$$
where
\begin{equation}    \label{coeff fun1}
{\mathbf v}(z) = \sum_{j=1}^\ell v_j(z) \cdot \ol{p}_j, \qquad v_j(z)
\in \C\zpart.
\end{equation}

For $c \in \C$, let $\on{Op}^{\on{irr},c}_{^L \g}(D)$ be the space of
$^L \g$-opers on the punctured disc of the form \eqref{coeff fun1},
where $v_j(z) \in z^{-d_j-1} \C[[z]]$ for $j=1,\ldots,\ell-1$, and
$$
v_\ell(z) = c z^{-d_\ell-2} + \ldots.
$$
Thus, for $c=0$ this is the space of $^L \g$-opers $D^\times$ with
regular singularity, and for $c \neq0$ this is the space of opers with
the mildest possible irregular singularity at the origin.

\begin{lem}
We have an isomorphism
$$
H_{\on{DS}}(\I_{1,\chi}) \simeq \on{Fun} \on{Op}^{\on{irr},c}_{^L
\g}(D)
$$
for some non-zero value of $c$.
\end{lem}

The value $c$ depends on the normalization of $\ol{p}_\ell$ as well as
the normalization of $\chi$, but it is inessential, because by
rescaling the coordinate $z$ we may identify the spaces
$\on{Op}^{\on{irr},c}_{^L \g}(D)$ for all non-zero values of $c$.

This has the following immediate corollary:

\begin{cor}    \label{with psi}
For any $c \in \C^\times$ the algebra ${\mc A}^\psi_\chi$ is
isomorphic to the algebra of functions on the space of $^L \g$-opers
on $\pone$ with regular singularity at $0 \in \pone$ and such that
their restriction to the punctured disc $D_\infty^\times$ at $\infty
\in \pone$ belongs to $\on{Op}^{\on{irr},c}_{^L \g}(D_\infty)$.
\end{cor}

Note that according to the definition of ${\mc A}^\psi_\chi$ (see
\secref{quant fd}), it is isomorphic to the center $Z(\g)$ of
$U(\g)$. In addition, the algebra ${\mc A}_0$ corresponding to
$\chi=0$ is also isomorphic to $Z(\g)$. Thus, we obtain that
$\on{Spec} Z(\g)$ is isomorphic to $\on{Op}^{\on{irr},c}_{^L \g}(D)$
for all $c \in \C$.

We are now ready to give a conjectural description of the spectra of
the KdV Hamiltonians.

\subsection{Spectra of the quantum KdV Hamiltonians}

We start with the case of $\ghat=\su$, which corresponds to the KdV
system.

Recall from \secref{spectra aff} that in the case of regular
semi-simple shift $\chi \in \h$ we expect, by analogy with the
finite-dimensional Gaudin models, that the spectra of the quantum
Hamiltonians are encoded by the second order differential operators on
$\pone$ with global coordinate $z$, which have regular singularity at
$z=0$ with the leading term $-\ell(\ell+1)/z^2$ and irregular
singularity at $\infty$ with the leading term $-\chi^2/s^4$, where
$s=z^{-1}$ is a local coordinate at $\infty \in \pone$. These opers
are (generically) represented by formula \eqref{final form rs}.

We stress again that this $\pone$ is a curve in the {\em spectral
parameter} of the affine Gaudin model.

Now, \corref{with psi} suggests that for the KdV system, in which
$\chi$ is nilpotent, the relevant space of $\su$-opers consists of the
second order differential operators on $\pone$ that have regular
singularity at $z=0$ with the leading term $-\ell(\ell+1)/z^2$ and
irregular singularity at $\infty$, but now of order $3$, that is, of
the form $c/s^3, c \in \C^\times$. By rescaling the coordinate $z$ on
$\pone$, we can make $c=1$, so without loss of generality we will
consider operators with singularity $1/s^3$ at $\infty$. The
corresponding $\su$-opers should therefore have the form
$$
\pa_z + \begin{pmatrix} 0 & v(z) + \la \\
  1 & 0 \end{pmatrix} + \frac{k}{z} d,
$$
where
\begin{equation}    \label{kdv form1}
v(z) = \frac{\ell(\ell+1)}{z^2} + \frac{1}{z} \left( 1 - \sum_{j=1}^m
\frac{k}{w_j} \right) + \sum_{j=1}^{m} \frac{2}{(z-w_j)^2} +
\sum_{j=1}^{m} \frac{k}{w_j} \frac{1}{z-w_j},
\end{equation}
and the coefficients $v_{j,k}$ in the expansion of $v(z)$ in
$z-w_j$ satisfy the equations \eqref{third order},
\begin{equation}    \label{third order again}
\frac{1}{4} \left( \frac{k}{w_j} \right)^3 - \frac{k}{w_j} v_{j,0} +
v_{j,1} = 0, \qquad j=1,\ldots,m.
\end{equation}
As explained in \secref{sec ord}, this is the condition of ``no
monodromy'' at the points $w_j$.

As in \secref{spectra aff}, applying gauge transformation by
$z^{kd}$, we obtain the operator
\begin{equation}
\pa_z + \begin{pmatrix} 0 & v(z) + \la z^k \\
  1 & 0 \end{pmatrix},
\end{equation}
which we rewrite as the following second order differential operator
with spectral parameter:
\begin{equation}    \label{final form kdv}
\pa_z^2 - \frac{1}{z} \left( 1 - \sum_{j=1}^m \frac{k}{w_j} \right) -
\frac{\ell(\ell+1)}{z^2} - \sum_{j=1}^{m} \frac{2}{(z-w_j)^2} -
\sum_{j=1}^{m} \frac{k}{w_j} \frac{1}{z-w_j} - \la z^k.
\end{equation}
The equations \eqref{third order again} are equivalent to the
condition that this operator has no monodromy around $w_j,
j=1,\ldots,m$, that is, no monodromy on $\pone$, except around the
points $0$ and $\infty$, for all values of $\la$.

According to our proposal, the differential operators of this kind
should encode the common eigenvalues of the quantum KdV Hamiltonians
on irreducible highest weight modules over the Virasoro algebra. Let
us recall that the Virasoro algebra obtained by the Drinfeld--Sokolov
reduction (with respect to the loop algebra of $\n_-$) from the affine
Kac--Moody algebra of level $k$ has the central charge
\begin{equation}    \label{ck}
c_k = 1 - \frac{6(k+1)^2}{k+2}
\end{equation}
(see \cite{FF:lmp}). Further, the Drinfeld--Sokolov reduction
of the irreducible $\su$-module with highest weight $\la = 2\ell$ and
level $k$ is either $0$ (which happens if and only if $k-\la \in
\Z_+$), or is the irreducible module over the Virasoro algebra with
highest weight (by which we mean the eigenvalue of $L_0$ on the
highest weight vector)
\begin{equation}    \label{delta}
\Delta_{\ell,k} = \frac{(2\ell+1)^2-(k+1)^2}{4(k+2)}
\end{equation}
(see \cite{FF:lmp}, p. 316).

Let us denote by $L_{\Delta,c}$ the irreducible module over the
Virasoro algebra with highest weight $\Delta$ and central charge $c$
(this is the irreducible quotient of the Verma module with this
highest weight).

Our proposal is then the following:

\begin{conj}
The generic common eigenvalues of the quantum KdV Hamiltonians on the
subspace of $L_{c_k,\Delta_{\ell,k}}$ of conformal dimension
$\Delta_{\ell,k} + m$ are in bijection with the set of operators of
the form \eqref{final form kdv} satisfying the equations \eqref{third
order again}.
\end{conj}

For special values of $\ell$ and $k$ the positions of some of the
poles $w_j$ coalesce or approach $0$, and these opers may degenerate,
as in the case of regular semi-simple $\chi$ discussed in
\secref{spectra aff}. Namely, there may appear points on $\pone$ at
which the relative positions of the two reductions in the
corresponding parabolic Miura oper are given by elements of $W \bs
W_{\on{aff}}$ other than $\ol{s}_0$. If the relative position at a
point $w_j$ is given by the element $\ol{s_0(s_1 s_0)^n} \in W \bs
W_{\on{aff}}$, then the expansion of the function $v(z)$ appearing in
\eqref{final form kdv} around $w_j$ will have the form
$$
v(z) \sim \frac{(2n+1)(2n+2)}{(z-w_j)^2} + \ldots .
$$
The coefficients in the expansion of this function at $w_j$ will also
have to satisfy an analogue of equation \eqref{third order again},
which means that our second order operator has no monodromy at $w_j$
for all $\la$.

In addition, when the Verma module $M_{c_k,\Delta_{\ell,k}}$ becomes
reducible, some of the points $w_j$ may tend to $0$. The common
eigenvalues encoded by the corresponding limits of \eqref{final form
kdv} should then correspond to the eigenvectors which appear in a
submodule of this Verma module generated by singular vectors.

It would be interesting to develop an analogue of the Bethe Ansatz
procedure for constructing eigenvectors of the quantum KdV
Hamiltonians. It is not obvious how to do this. As we have seen above,
the Bethe Ansatz appears when we represent our opers as the Miura
transformation of Cartan connections. This is possible in the Gaudin
models associated to regular semi-simple $\chi$, as discussed in
\secref{spectra aff}. But the opers \eqref{final form kdv} cannot be
represented as Miura transformation of Cartan connections defined
globally on $\pone$. (This becomes clear when we look at the
restriction of \eqref{final form kdv} to the punctured disc at the
point $z=\infty$: the operator has the form $\pa_s^2 + \wt{v}(s)$,
where $\wt{v}(s) = c/s^3 + \ldots$, $c \in \C^\times$, and $s=z^{-1}$
is the local coordinate at $\infty$; this $\wt{v}(s)$ cannot be
represented in the Miura form $u(s)^2 + \pa_s u(z)$, where $u(s)$ is a
Laurent power series in $s$, as easily checked by substitution.)

\subsection{Two different spaces of opers}    \label{two ways}

In the study of the KdV system, second order differential operators,
which we call $\su$-opers, appear in two different ways. First, they
appear as points of the phase space ${\mc M}(\ghat)$ of the classical
KdV system, as described in \secref{class kdv}. These operators have
the form
$$
\pa_t^2 - v(t) - z,
$$
where $t$ is a coordinate on a circle, or a punctured disc, and $v(t)$
is a function on it. (Despite what the notation may suggest, $t$ is
the ``space'', not the ``time'', coordinate of the KdV equation.) On
the other hand, $z$ is the spectral parameter. The classical KdV
Hamiltonians are constructed by expanding the monodromy matrix of this
operator, considered as function of $z$, near $z=0$ and $z=\infty$.

The $\su$-opers appear again as the objects encoding the eigenvalues
of the quantum KdV Hamiltonians. They have the form
$$
\pa_z^2 - v(z) - \lambda z^k,
$$
where $z$ is the spectral parameter of the classical KdV system. Thus,
$z$ is a coordinate on $\pone$, and $v(z)$ is a meromorphic function
on this $\pone$ (see formula \eqref{final form kdv}). On the other
hand $\la$, is a new spectral parameter.

It is important to realize that these two spaces of opers are not
related to each other. In fact, they correspond to two {\em dual}
theories, one of which is classical and the other is quantum. In the
classical theory the affine Kac--Moody algebra $\su$ is the central
extension of $\sw_2\zpart$, where $z$ is the spectral parameter, and
the $\su$-opers ``live'' on a circle with coordinate $t$, which is the
space coordinate of the KdV hierarchy. In the quantum setting, $\su$
is the central extension of $\sw_2\pparl$, and the $\su$-opers
``live'' on a projective line $\pone$, whose coordinate $z$ is the
spectral parameter of the classical system.

Even more important is the fact that if we replace the Lie algebra
$\su$ of the classical theory by a more general affine Kac--Moody
algebra $\ghat$, then the corresponding Lie algebra of the quantum
theory will be not $\ghat$, but its {\em Langlands dual} affine Lie
algebra $^L \ghat$, as we will discuss below. This makes it more clear
that the classical and quantum theories are really different.

\subsection{Comparison with the conjecture of
  Bazhanov, Lukyanov, and Zamolodchikov}

We wish to compare our proposal to the conjectural relation between
the eigenvalues of the quantum KdV Hamiltonians and second order
differential operators on $\pone$ of special kind, proposed by
Bazhanov, Lukyanov, and Zamolodchikov in \cite{BLZ} (generalizing the
proposal made earlier in \cite{DT,BLZ4} for eigenvalues on the highest
weight vectors). Understanding the underlying reasons for this
relation was in fact one of the main motivations for this paper.

Let us assume that $k \neq -2$ and set
$$
\al = -\frac{k+1}{k+2}.
$$
We will apply the change of variables
$$
z = \frac{x^{2\al+2}}{(2\al+2)^2}
$$
to the differential operator \eqref{final form kdv}.

The general transformation formula for a second order operator
$$
\pa_z^2 - v(z): \Omega^{-1/2} \to \Omega^{3/2}
$$
under the change of variables $z=\varphi(x)$ is
\begin{equation}    \label{proj conn}
v(z) \mapsto v(\varphi(x)) \left( \varphi'
\right)^2 - \frac{1}{2} \{ \varphi,x \},
\end{equation}
where
\begin{equation}    \label{schwarzian}
\{\varphi,x\} = \frac{\varphi'''}{\varphi'} - \frac{3}{2}\left(
\frac{\varphi''}{\varphi'} \right)^2
\end{equation}
is the Schwarzian derivative of $\varphi$.\footnote{Note that it is
imperative to consider our second order operators as acting from
$\Omega^{-1/2}$ to $\Omega^{3/2}$, for this is the only way to ensure
that their property of having the principal symbol $1$ and
subprincipal symbol $0$ is coordinate-independent.}

After a straightforward calculation, we find that the new operator
looks as follows:
\begin{equation}    \label{new operator}
\pa_x^2 - \frac{\wt\ell(\wt\ell+1)}{x^2} -
x^{2\al} + 2 \frac{d^2}{dx^2} \sum_{j=1}^m \log(x^{2\al+2}
- z_j) + E,
\end{equation}
where
\begin{align*}
\wt\ell(\wt\ell+1) &= 4(\al+1)^2 \ell(\ell+1) + \al^2 + 2\al +
\frac{3}{4} = 4(\al+1)\Delta_{k,\ell} + \al^2 - \frac{1}{4},
\\ z_j &= (2\al+2)^2 w_j, \\
E &= -(2\al+2)^{\frac{2\al}{\al+1}} \la
\end{align*}

This operator coincides with the Schr\"odinger operators
parameterizing the spectra of the quantum KdV Hamiltonians on
$L_{c_k,\Delta_{\ell,k}}$ introduced in \cite{BLZ} (formula
(1)).\footnote{Note that $\wt\ell$ corresponds to the $\ell$ in
\cite{BLZ}, but it is not equal to our $\ell$.} Moreover, our
condition \eqref{third order again} means that the operator \eqref{new
operator} has no monodromy around the points $z_j$ for all $E$ and is
therefore equivalent to the algebraic equations of \cite{BLZ} (formula
(3)).

Thus, we see that our proposal reproduces the proposal of \cite{BLZ}
in the case of the quantum KdV system! We view this as a confirmation
of the correspondence between the spectra of the Hamiltonians of
quantum soliton systems and opers for the Langlands dual affine
algebras which we propose in this paper.

\subsection{Generalization to the case of $\wh{\sw}_n$}

It is easy to generalize formula \eqref{final form kdv} to the case of
$\wh\sw_n$. In this case the Langlands dual Lie algebra is also
$\wh\sw_n$. The $\wh\sw_n$-opers may be represented as $n$th order
differential operators with spectral parameter. We wish to describe
the differential operators of this kind on $\pone$ that encode the
spectra of the quantum Hamiltonians of the $\wh\sw_n$-KdV system
acting on the irreducible module over the ${\mc W}$-algebra ${\mc
W}_k(\sw_n)$ obtained by the quantum Drinfeld--Sokolov reduction of
the irreducible $\ghat$-module $L_{\nu,k}$ of highest weight $\nu \in
\h^*$ and level $k$. The quantum Drinfeld--Sokolov reduction is
described in \cite{FF:gd,FKW}. In particular, the central charge of
${\mc W}_k(\sw_n)$ is given by the standard formula
$$
c = 1+2n(n+1) - n(n+1)\left( (k+n) + \frac{1}{k+n} \right).
$$

These operators have the following general form:
\begin{equation}    \label{nth order oper}
(-\pa_z)^n - v_1(z) (-\pa_z)^{n-2} - \ldots + v_{n-2}(z) \pa_z -
v_{n-1}(z) - z^k \la,
\end{equation}
and act from $\Omega^{-(n-1)/2}$ to $\Omega^{(n+1)/2}$. This
determines their transformation properties under the changes of
coordinate $z$; in particular, this ensures that the property that
their principal symbol is $(-1)^n$ and the subprincipal symbol is $0$
is preserved by the changes of coordinate. The coefficients $v_r(z)$
will be rational functions on $\pone$ with poles at $z=0, \infty$, and
finitely many other points $w_j, j=1,\ldots,m$ (where $m$ is the
degree of the corresponding module over ${\mc W}_k(\sw_n)$). These
coefficients should satisfy the following properties:

\medskip

\noindent (1) At $z=\infty$ the operator \eqref{nth order oper}
has the mildest possible irregular singularity; namely, we have
\begin{align*}
\wt{v}_r(s) &\sim \frac{\wt{c}_r}{s^{r+1}} + \ldots, \qquad
r=1,\ldots,n-2; \\
\wt{v}_{n-1}(s) &\sim \frac{1}{s^{n+1}} + \ldots,
\end{align*}
where $\wt{v}_r(s)$ are the coefficients of the operator obtained from
\eqref{nth order oper} by the change of variables $z \mapsto s =
z^{-1}$. Here, as in the case of $\sw_2$, we normalize the leading
coefficient of $\wt{v}_{n-1}(s)$ to be equal to $1$. This may always
be achieved by rescaling $z$.

\smallskip

\noindent (2) Near $z=0$ the operator \eqref{nth order oper} has
regular singularity, that is
$$
v_r(z) \sim \frac{c_r(\nu)}{z^{r+1}} + \ldots,
$$
where the coefficients $c_r(\nu)$ are determined by the highest weight
$\nu$ of the $\ghat$-module $L_{\nu,k}$. Namely, representing $\nu$ as
$(\nu_1,\ldots,\nu_n)$, where $\nu_i \in \C$ and $\sum_{i=1}^n \nu_i =
0$, we find the $c_r(\nu)$'s from the following formula:
\begin{equation}    \label{zero}
(-\pa_z)^n - \sum_{i=1}^{n-1} \frac{c_i(\nu)}{z^{i+1}} (-\pa_z)^{n-i-1} =
\left(-\pa_z + \frac{\nu_1}{z} \right) \ldots \left(-\pa_z +
\frac{\nu_n}{z} \right).
\end{equation}

\smallskip

\noindent (3) At the points $w_j$ the operator \eqref{nth order oper}
has regular singularity,
\begin{equation*}
v_r(z) \sim \frac{c_r(\theta)}{(z-w_j)^{r+1}} + \ldots,
\end{equation*}
where $\theta = (1,0\ldots,0,-1)$ is the maximal root of $\sw_n$,
which is the highest weight of the adjoint representation. There is
also the additional requirement that the operator \eqref{nth order
oper} has {\em trivial monodromy} around the point $w_j$ for each
$j=1,\ldots,m$, and all $\la$.

\medskip

Now, our proposal is that the $n$th order differential operators of
this kind should correspond to the generic common eigenvalues of the
quantum $\wh\sw_n$-KdV Hamiltonians on the degree $\Delta+m$ subspace
of the irreducible module obtained by the quantum Drinfeld--Sokolov
reduction of the irreducible $\ghat$-module $L_{\nu,k}$ of highest
weight $\nu \in \h^*$ and level $k$. In general, some of the poles
$w_j$ may coalesce, as in the case of $\su$.

On the other hand, in \cite{BHK} third order differential operators
have been associated to the common eigenvalues of the quantum
$\wh\sw_3$-KdV (Boussinesq) Hamiltonians on the highest weight vectors
in irreducible representations of the ${\mc W}$-algebra associated to
$\wh\sw_3$. Analogous formulas have been proposed in
\cite{DDT,Dorey1,Dorey2} in the case of $\wh\sw_n$. As in the case of
$\su$ discussed in the previous section, we expect that by applying a
change of variables to the above operators when the set of points
$w_j, j=1,\ldots,m$, is empty (this corresponds to the highest weight
vector), we will obtain the operators of \cite{BHK,DDT,Dorey1,Dorey2}.
For $m>0$ we obtain a generalization of the proposal of
\cite{BHK,DDT,Dorey1,Dorey2} to non-highest weight vectors.

Explicitly, this change of variables is given by the formula
$$
z = \frac{x^{n\al+n}}{(n\al+n)^n}
$$
where
$$
\al = -\frac{k+n-1}{k+n}.
$$
It is easy to see that under this change of variables the irregular
term of our differential operator gives rise to the term $x^{n\al}$,
and the term $\la x^k$ gives rise to the new spectral parameter
term
$$
E = (n\al+n)^{\frac{n\al}{\al+1}} \la,
$$
which is independent of $x$. The poles of the
new operator will be at the points $x=0,\infty$ and $z_j,
j=1,\ldots,m$, where
$$
z_j = (n\al+n)^n w_j.
$$

If the affine Lie algebra $^L \ghat$ is of classical type, we may also
realize affine opers as scalar (pseudo)differential operators,
following \cite{DS}. It is possible to describe those of them that
encode the spectra of the quantum Hamiltonians of the $\ghat$-KdV
system in a similar way to the case of $\wh\sw_n$. We expect that the
corresponding operators (in the special case when the set of extra
points $w_j$ is empty) are related to those proposed in
\cite{DDT,Dorey1,Dorey2} after a change of variables. However, if we
wish to describe these affine opers for an arbitrary $^L \ghat$, we
should use instead the canonical form for opers discussed in the next
section.

\subsection{Generalization to an arbitrary affine algebra $\ghat$}

It is also easy to generalize formula \eqref{final form kdv} to the
case of an arbitrary affine Kac--Moody algebra. For simplicity we will
assume here that $^L \ghat$ is an untwisted affine Kac--Moody algebra
$\lhat$, where ${\mf l}$ is a simple Lie algebra, but the general case
may be analyzed along the same lines.

Let us recall (see, e.g., \cite{F:flag}) that ${\mf l}$-opers have
canonical form
$$
\nabla = \pa_z + \ol{p}_{-1} + \sum_{r=1}^\ell v_r(z) \cdot \ol{p}_r,
$$
where $\ol{p}_r, r=1,\ldots,\ell$, are the elements defined in
\secref{class kdv}. Here $v_1(z)$ transforms as a projective
connection (that is, by formula \eqref{proj conn}) and $v_r(z), r>1$,
transforms as a $(d_r+1)$-differential under the changes of coordinate
$z$. Our $\lhat$-opers on $\pone$ will have the form
\begin{equation}    \label{oper spectra}
\pa_z + \ol{p}_{-1} + \sum_{r=1}^\ell v_r(z) \cdot \ol{p}_r + \la z^k
\ol{p}_\ell,
\end{equation}
where each $v_r(z)$ is a rational function on $\pone$ with poles at
$z=0, \infty$, and finitely many other points $w_j,
j=1,\ldots,m$. These coefficients should satisfy the following
properties:

\medskip

\noindent (1) At $z=\infty$ the operator \eqref{oper spectra}
has the form
\begin{align*}
\wt{v}_r(s) &\sim \frac{\wt{c}_r}{s^{d_r+1}} + \ldots, \qquad
r=1,\ldots,\ell-1; \\
\wt{v}_{\ell}(s) &\sim \frac{1}{s^{d_\ell+2}} + \ldots.
\end{align*}
where $\wt{v}_r(s)$ are the coefficients of the operator obtained from
\eqref{oper spectra} by the change of variables $z \mapsto s =
z^{-1}$.

\smallskip

\noindent (2) Near $z=0$ the operator \eqref{oper spectra} has
regular singularity, that is
$$
v_r(z) \sim \frac{c_r(\nu)}{z^{d_r+1}} + \ldots,
$$
where the coefficients $c_r(\nu)$ are determined by the highest weight
$\nu$ of the $\ghat$-module $L_{\nu,k}$. They are determined by the
following rule: the element $$\ol{p}_{-1} + \sum_{r=1}^\ell \left(
c_r(\nu) + \frac{1}{4} \delta_{r,1} \right) \ol{p}_r$$ is the unique
element in the Kostant slice of regular elements,
$$
\ol{p}_{-1} + {\mb v}, \qquad {\mb v} \in \bb,
$$
which is conjugate to $\ol{p}_{-1} - \nu$ (see \cite{F:book},
Sect. 9.1).

\smallskip

\noindent (3) At the points $w_j$ the operator \eqref{oper spectra}
has regular singularity,
\begin{equation*}
v_r(z) \sim \frac{c_r(\theta)}{(z-w_j)^{d_r+1}} + \ldots,
\end{equation*}
where $\theta$ is the maximal coroot of ${\mf l}$. In addition, we
require that the oper have {\em trivial monodromy} around the
point $w_j$ for each $j=1,\ldots,m$, and all $\la$.

\medskip

Our proposal is that the $^L \ghat$-opers \eqref{oper spectra}
satisfying these conditions encode the generic common eigenvalues of
the quantum Hamiltonians of the $\ghat$-KdV system on the irreducible
modules of the corresponding ${\mc W}$-algebra obtained by the
quantum Drinfeld--Sokolov reduction \cite{FF:gd,FKW} from the
irreducible module over $\ghat$ with highest weight $\nu$ and level
$k$.

\section{Affine Gaudin models with regular singularities}
\label{reg sing}

In the previous section we have discussed the affine Gaudin models
corresponding to the shift of argument method. We have seen that the
classical Hamiltonians of these models are integrals of motion of
familiar soliton hierarchies, such as the KdV or AKNS
hierarchies. Hence the quantization of these Gaudin models naturally
leads to the corresponding quantum Hamiltonians and their
diagonalization. However, in the universe of all Gaudin models (see
\cite{FFT}), the shift of argument Gaudin models are very special:
they correspond to two points on $\pone$, one with regular singularity
and one with irregular singularity of order $2$. It is natural to
develop a similar theory for more general, multi-point, Gaudin models.

The simplest of these models are the ones with ``regular
singularities'' at finitely many points of $\pone$, in the terminology
of \cite{FFT}. In this section we study the corresponding affine
Gaudin models. We start out by briefly summarizing the main features
of the Gaudin models with regular singularities in the
finite-dimensional case. We then discuss analogous patterns in the
affine Gaudin models and their potential connection to the quantum KdV
system. In particular, we will suggest an explicit construction of
eigenvectors of the quantum KdV Hamiltonians on representations of the
unitary minimal models via Bethe Ansatz.

\subsection{Gaudin model, finite-dimensional case}    \label{RS fd}

Let $\g$ be a simple finite-dimensional Lie algebra over $\C$. Recall
that it carries the invariant inner product $\ka_0$ normalized in such
a way that the squared length of the maximal root is equal to $2$. Let
$\{ J^a \}$ be a basis of $\g$ and $\{ J_a \}$ the dual basis with
respect to this inner product. The {\em Gaudin Hamiltonians} are the
elements of $U(\g)^{\otimes N}$ given by the following formula:
\begin{equation}    \label{Gaudin ham}
\Xi_i = \sum_{j\neq i} \sum_{a} \frac{J_a^{(i)}
J^{a(j)}}{z_i-z_j}, \qquad i=1,\ldots,N.
\end{equation}
Here $A^{(i)} = 1 \otimes \ldots \otimes A \otimes \ldots \otimes 1$,
with $A$ in the $i$th position. 

As shown in \cite{FFR,F:faro}, the quadratic Gaudin Hamiltonians are
included in a large commutative subalgebra of $U(\g)^{\otimes N}$,
called the {\em Gaudin algebra}. According to \cite{F:faro},
this algebra is isomorphic to the algebra $\on{Fun} \on{Op}_{^L
\g}^{\on{RS}}(\pone)_{(z_i),\infty}$ of functions on the space
$\on{Op}_{^L \g}^{\on{RS}}(\pone)_{(z_i),\infty}$ of $^L \g$-opers on
$\pone$ with regular singularities at the points $z_1,\ldots,z_N$ and
$\infty$. Here $^L \g$ is the Langlands dual Lie algebra of $\g$,
whose Cartan matrix is the transpose of the Cartan matrix of $\g$.

This isomorphism comes about \cite{FFR,F:faro} in the same way as in
the case of the shift of argument Gaudin models, discussed in
\secref{shift fd}. Namely, using the coinvariants (or conformal
blocks) of modules over the affine algebra $\ghat$ of critical level,
we identify the Gaudin algebra with a quotient of the algebra
$\on{End}_{\ghat} \V_{\on{crit}}$, where $\V_{\on{crit}}$ is the
vacuum module of critical level (see formula \eqref{vacuum}). The
latter algebra is isomorphic to the algebra of functions on the space
$\on{Op}_{^L \g}(D)$ of $^L \g$-opers on the (formal) disc
\cite{FF:gd,F:wak}. This implies that the Gaudin algebra is isomorphic
to the algebra $\on{Fun} \on{Op}_{^L
\g}^{\on{RS}}(\pone)_{(z_i),\infty}$ of functions on $\on{Op}_{^L
\g}^{\on{RS}}(\pone)_{(z_i),\infty}$.

As before, the appearance of the Langlands dual Lie algebra $^L \g$
here is important. It signifies the fact that the above isomorphism
may be viewed as a special case of the geometric Langlands duality
(for more on this, see \cite{F:icmp,F:rev,FFT}).

Thus, we obtain that the common eigenvalues of this Gaudin algebra,
including those of the Gaudin Hamiltonians $\Xi_i$, are encoded by
$^L \g$-opers in $\on{Op}_{^L \g}^{\on{RS}}(\pone)_{(z_i),\infty}$.
For example, in the case of $\g=\sw_2$ the corresponding opers are
the second order operators of the form
\begin{equation}    \label{fuchsian}
\pa_z^2 - \sum_{i=1}^N \frac{\Delta_i}{(z-z_i)^2} - \sum_{i=1}^N
\frac{c_i}{z-z_i}.
\end{equation}
In particular, the eigenvalues of the Gaudin algebra on the tensor
product of the Verma modules $\bigotimes_{i=1}^N M_{\nu_i}$ are
encoded by such differential operators, in which
$$
\Delta_i = \frac{1}{4} \nu_i(\nu_i+2);
$$
this is the value of the Casimir operator in $U(\sw_2)$ on the $i$th
factor $M_{\nu_i}$. The operator \eqref{fuchsian} is then
determined by the remaining parameters $c_i, i=1,\ldots,N$, and in
fact they are nothing but the eigenvalues of $\Xi_i, i=1,\ldots,N$
(see \cite{F:faro}).

Thus, in the case of $\g = \sw_2$ the Gaudin algebra is generated by
the $\Xi_i$'s and the Casimir operators along each factor in
$U(\sw_2)^{\otimes N}$. But for all other simple Lie algebras there
are higher order Gaudin Hamiltonians, whose eigenvalues are recorded
by polynomials in the coefficients appearing in the expansions of the
$^L \g$-opers near the points $z_i$. For instance, for $\g=\sw_n$ the
role of \eqref{fuchsian} is played by the $n$th order differential
operators with regular singularities:
$$
(-\pa_z)^n - v_1(z) (-\pa_z)^{n-2} - \ldots + v_{n-2}(z) \pa_z -
v_{n-1}(z),
$$
where $v_j(z)$ has the form
$$
v_j(z) = \sum_{i=1}^N \sum_{k=1}^{j+1} \frac{c_{ijk}}{(z-z_i)^k}.
$$
Thus, such an operator is determined by the coefficients $c_{ijk}$
appearing in its expansion around the singular points $z_i,
i=1,\ldots,N$. The quadratic Gaudin Hamiltonians correspond to
$c_{i,1,1}$, and the coefficients with $j>1$ correspond to higher
Gaudin Hamiltonians of order $j+1$.

\subsection{Bethe Ansatz}

Let us choose a collection of weights $\nu_1,\ldots,\nu_N \in
\h^*$. Then we have the Verma modules $M_{\nu_i}$ with highest
weights $\nu_i, i=1,\ldots,\nu_N$. Eigenvectors of the higher Gaudin
Hamiltonians on the tensor product $\bigotimes_{i=1}^N M_{\nu_i}$ may
be constructed by Bethe Ansatz (see \cite{FFR,F:faro}).

Namely, for any collection of distinct complex numbers
$w_1,\ldots,w_m$ such that $w_j \neq z_i$ for all $i$ and $j$, and a
collection $\al_{i_1},\ldots,\al_{i_m}$ of simple roots of $\g$ (so
that $i_j \in I$, the set of vertices of the Dynkin diagram of $\g$)
consider the following {\em Bethe vector}
\begin{align} \label{genbv}
&\phi(w_1^{i_1},\ldots,w_m^{i_m})\\
&= \sum_{p=(P^1,\ldots,P^N)}
\prod_{k=1}^N \frac{f_{i^k_1}^{(k)} f_{i^k_2}^{(k)} \ldots
f_{i^k_{a_k}}^{(k)}}{(w_{j^k_1}-w_{j^k_2})(w_{j^k_2}-w_{j^k_3}) \ldots
(w_{j^k_{a_k}}-z_k)} v_{\nu_1} \otimes \ldots v_{\nu_N}. \nonumber
\end{align}
Here $f_i^{(k)}$ denotes the generator of the Lie
algebra $\n_- \subset \g$ corresponding to the $i$th simple root,
acting on the $k$th factor of $\bigotimes_{j=1}^N M_{\nu_j}$, and the
summation is taken over all {\em ordered} partitions $P^1 \cup P^2
\cup \ldots \cup P^N$ of the set $\{1,\ldots,m\}$, where $P^k = \{
j^k_1,j^k_2,\ldots,j^k_{a_k} \}$ and $\{ i^k_1,i^k_2,\ldots,i^k_{a_k}
\}$ are the corresponding labels in $I$.

It is proved in \cite{FFR,F:faro} that the vector
$\phi(w_1^{i_1},\ldots,w_m^{i_m})$ is a common eigenvector of all
higher Gaudin Hamiltonians if the following system of equations is
satisfied
\begin{equation} \label{bethe}
\sum_{i=1}^N \frac{\langle \chal_{i_j},\nu_i \rangle}{w_j-z_i} -
\sum_{s\neq j} \frac{\langle \chal_{i_j},\al_{i_s} \rangle}{w_j-w_s} =
0, \quad j=1,\ldots,m,
\end{equation}
where $\langle \cdot,\cdot \rangle$ denotes the pairing between
$\h^*$ and $\h$. These are the {\em Bethe Ansatz equations} of
the Gaudin model.

As shown in \cite{FFR,F:faro}, the common eigenvalues of the higher
Gaudin Hamiltonians on the Bethe vector
$\phi(w_1^{i_1},\ldots,w_m^{i_m})$ given by formula \eqref{genbv}
(provided that it is non-zero) are encoded by the $^L \g$-oper
obtained by applying the {\em Miura transformation} to the Cartan
connection on the $^L H$-bundle $\Omega^\rho$ over $\pone$,
$$
\ol\nabla = \pa_z - \sum_{i=1}^N \frac{\nu_i}{z-z_i} + \sum_{j=1}^m
\frac{\al_{i_j}}{z-w_j},
$$
with regular singularities at the points $z_i$ and $w_j$. The Miura
transformation maps such a connection to an $^L \g$-oper on $\pone$
with regular singularities at these points. This oper has no
singularities at the points $w_j$ if and only if the Bethe Ansatz
equations \eqref{bethe} are satisfied (see
\cite{FFR,F:faro}). Equivalently, we can phrase this as saying that
the Bethe Ansatz equations mean that the Miura transformation of the
connection $\ol\nabla$ have {\em no monodromy} around $w_j,
j=1,\ldots,m$.


There is an additional restriction on the common eigenvalues of the
Gaudin Hamiltonians that may occur in the tensor product of
finite-dimensional representations, and that is the {\em condition of
no monodromy around the points $z_i, i=1,\ldots,N$, and
$\infty$}. Denote by $V_\nu$ the irreducible finite-dimensional
$\g$-module with highest weight $\nu$, which is dominant and
integral. The higher Gaudin Hamiltonians act on $\bigotimes_{i=1}^N
V_{\nu_i}$ and commute with the diagonal $\g$-action. Let us describe
the joint eigenvalues of these Hamiltonians on the isotypical
component
$$
\left( \bigotimes_{i=1}^N V_{\nu_i} \otimes V_{\nu_\infty}
\right)^{\g},
$$
where $\nu_\infty$ is another dominant integral weight. According to
\cite{F:flag}, Corollary 4.8 (see also \cite{F:faro}, Theorem
2.7,(3)), the set of $^L \g$-opers representing them (counted without
multiplicities) is a subset of the set of $^L \g$-opers on $\pone$
with regular singularities at $z_i, i=1,\ldots,N$, and $\infty$, with
fixed residues determined by $\nu_i, i=1,\ldots,N$, and $\nu_\infty$,
and with trivial monodromy representation (such opers may equivalently
be described as $\nu_i$-regular at $z_i$ and $\nu_\infty$-regular at
$\infty$, see below). Furthermore, according to Conjecture 1 of
\cite{F:faro}, this subset is in fact the entire set of such opers. In
other words, all opers of this kind correspond to common eigenvalues
of the Gaudin Hamiltonians.

\subsection{Generalization to the affine case}    \label{gen aff}

Now we would like to generalize the above results on the Gaudin models
to the case when a finite-dimensional simple Lie algebra $\g$ is
replaced by an affine Kac--Moody algebra $\gtil$, by analogy to the
shift of argument models discussed above. We begin by summarizing the
known facts about the affine Gaudin model.

Let $\wt\De$ be the set of roots of an affine Kac--Moody algebra
$\gtil$ (see \secref{aff km}). The Lie algebra $\gtil$ has a canonical
invariant inner product $\wt\ka_0$ whose restriction to $\g \otimes 1$
is equal to $\ka_0$. Let us choose dual bases $\{ J^{a} \}$ and $\{
J_{a} \}$ of $\g$ with respect to $\ka_0$.  Then we choose dual bases
$\{J^{\wt{a}}\},\{J_{\wt{a}} \}$ of $\gtil$ with respect to $\wt\ka_0$
in the following way. The former consists of $J^{a} \otimes t^n, n \in
\Z$, ${\mb 1}$ and $d$, and the latter consists of $J_{a} \otimes t^n,
n \in \Z$, $d$ and ${\mb 1}$.

Let $z_1,\ldots,z_N$ be a set of $N$ distinct complex numbers.
The quadratic Gaudin Hamiltonians associated to $\gtil$ are given by
the same formula \eqref{Gaudin ham} as in the finite-dimensional case:
$$
\Xi_i = \sum_{j\neq i} \sum_{\wt{a}} \frac{J_{\wt{a}}^{(i)}
J^{\wt{a}(j)}}{z_i-z_j}, \qquad i=1,\ldots,N
$$
(in fact, such Hamiltonians may be defined for any symmetrizable
Kac--Moody algebra). The elements $\Xi_i$ belong to a completion of
the $N$-fold tensor power of the universal enveloping algebra
$U(\gtil)$. It gives rise to a well-defined linear operator on the
tensor product $\bigotimes_{i=1}^N M_i$, if each $M_i$ is a {\em
smooth} $\gtil$-module. By definition, a $\gtil$-module $M$ is called
smooth if for any $v\in M$ we have $\g \otimes t^m\C[[t]] v=0$ for
sufficiently large $m\in\Z_+$. In particular, highest weight modules,
or modules from the category ${\mc O}$, are smooth.

It is clear that the Hamiltonians $\Xi_i$ commute with the diagonal
action of $\gtil$ and that
\begin{equation}    \label{zero sum}
\sum_{i=1}^N \Xi_i = 0.
\end{equation}

For a weight $\nu \in \htil^*$, denote by $M_\nu$ the Verma
module
$$
M_\nu = \on{Ind}_{\ntil_+ \oplus \htil}^{\gtil} \C_\nu.
$$
Consider the problem of diagonalization of the quadratic Gaudin
Hamiltonians on the tensor product of Verma modules
$\bigotimes_{i=1}^N M_{\nu_i}$. The Bethe Ansatz presented in the
previous section generalizes in a straightforward way to
the case of affine algebras (and, moreover, to all symmetrizable
Kac--Moody algebras).

More precisely, for any collection of distinct complex numbers
$w_1,\ldots,w_m$ such that $w_j \neq z_i$ for all $i$ and $j$, and a
collection $\al_{i_1},\ldots,\al_{i_m}$ of simple roots of $\gtil$ (so
that $i_j \in \wt{I}$, the set of vertices of the Dynkin diagram of
$\gtil$) we have the {\em Bethe vector} given by formula
\eqref{genbv}. Then, according to the results of \cite{FFR}, Sect. 5
(based on the computations of \cite{SV}), this vector is an
eigenvector of the quadratic Gaudin Hamiltonians $\Xi_i,
i=1,\ldots,N$, if the Bethe Ansatz equations \eqref{bethe} are
satisfied.\footnote{This may also be proved along the lines of
\cite{RV}.} This is in complete analogy with the finite-dimensional
case.

Thus, so far in the affine case we have the quadratic Gaudin
Hamiltonians and the Bethe vectors \eqref{bethe}. However, by analogy
with the finite-dimensional case (and the affine shift of argument
Gaudin models), we expect they are included in a large commutative
subalgebra of a completion of $U(\gtil)^{\otimes N}$. This raises the
following questions:

\begin{itemize}

\item[(1)] Are there higher order affine Gaudin Hamiltonians?

\item[(2)] What parametrizes the joint eigenvalues of the affine
  Gaudin Hamiltonians?

\end{itemize}

In the next section we will take up the first question in the simplest
case of two points $z_1, z_2$ (or, more properly, three points, since
we should count $\infty$ as an extra point). In this case there
is only one quadratic Gaudin Hamiltonian, $\Xi_1$, because $\Xi_2 = -
\Xi_1$ due to the identity \eqref{zero sum}. We will show that $\Xi_1$
is nothing but the operator $L_0$ of the {\em coset Virasoro algebra}
acting on the tensor product of two $\ghat$-modules $M_1 \otimes
M_2$. Therefore all quantum KdV Hamiltonians $H_{2i+1}, i \geq 0$,
written with respect to this Virasoro algebra, commute with $\Xi_1$
(which is equal to $L_0 = H_1$). We will conjecture they these are
precisely the higher (local) affine Gaudin Hamiltonians.

In the case of multiple points we do not have an analogous formula for
the local Gaudin Hamiltonians, but we will give a procedure for
constructing the {\em non-local} Gaudin Hamiltonians (by analogy with
the above construction of the non-local quantum AKNS and KdV
Hamiltonians). We will also consider the quasi-classical limit of this
construction, which is an affine analogue of the Hitchin system.

We will then discuss the second question of describing the spectra of
the affine Gaudin Hamiltonians. By analogy with the finite-dimensional
case, and the affine shift of argument Gaudin models, we will
conjecture that the joint eigenvalues of the affine Gaudin
Hamiltonians are parameterized by the {\em affine
$^L\ghat$-opers} with regular singularities at the points
$z_1,\ldots,z_N$ and $\infty$.

\subsection{The case of two points}    \label{two points local}

Let us consider the simplest case of two distinct points $z_1, z_2$.
Thus, we have two smooth $\gtil$-modules $M_1$ and $M_2$ and the
quadratic Gaudin Hamiltonian
$$
\Xi = \sum_{\wt{a}} J_{\wt{a}} \otimes J^{\wt{a}}
$$
acting on $M_1 \otimes M_2$. Note that $\Xi_1 = - \Xi_2 =
\Xi/(z_1-z_2)$. It commutes with the diagonal action of $\gtil$ on
$M_1 \otimes M_2$. Our question is to construct the higher affine
Gaudin Hamiltonians, which commute with the diagonal action
of $\ghat$ and with $\Xi$.

We will construct them using the coset construction of
\cite{GKO}. Suppose that the levels $k_1$ and $k_2$ of the modules
$M_1$ and $M_2$ are non-critical, i.e., $k_i \neq -h^\vee$ and in
addition $k_1 + k_2 \neq -h^\vee$. Then the GKO construction gives us
the Virasoro algebra which commutes with the diagonal action of
$\ghat$ in the completed tensor product $U_{k_1}(\ghat) \wh\otimes
U_{k_2}(\ghat)$, where $U_k(\ghat) = U(\ghat)/({\mb 1} - k)$. This
algebra is generated by $L^{\on{GKO}}_n, n \in \Z$, such that
$$
T^{\on{GKO}}(w) = \sum_{n \in \Z} L^{\on{GKO}}_n w^{-n-2} =
T^{\on{SS}}_1(w) + T^{\on{SS}}_2(w) - T^{\on{SS}}_{\on{diag}}(w).
$$
Here $T^{\on{SS}}(w)$ is the Segal--Sugawara current at level $k \neq
- h^\vee$,
\begin{align*}
T^{\on{SS}}(w) &= \sum_{n \in \Z} L^{\on{SS}}_n w^{-n-2} =
\frac{1}{2(k+h^\vee)} \sum_{a} \Wick J^{a}(w) J_{a}(w)
\Wick \, \\ J^{a}(w) &= \sum_{n \in \Z} J^{a}_n w^{-n-1},
\qquad J^{a}_n = J^{a} \otimes t^n,
\end{align*}
and the lower indices refer to the first, the second or the diagonal
$\ghat$ in $U_{k_1}(\ghat) \wh\otimes U_{k_2}(\ghat)$.

Without loss of generality we may, and will, assume that the action
of $L^{\on{SS}}_0$ on the $\gtil$-modules $M_1$ and $M_2$ coincides
with the action of $-d \in \gtil$. Then we find that
\begin{align*}
-(k_1+k_2+h^\vee) L^{\on{GKO}}_0 &= \sum_{a} \sum_{n \in \Z} J_{a,n}
\otimes J^{a}_n - k_1 \otimes L^{\on{SS}}_0 - L^{\on{SS}}_0 \otimes
k_2 \\ &= \sum_{\wt{a}} J_{\wt{a}} \otimes J^{\wt{a}} = \Xi.
\end{align*}

Thus, up to a non-zero scalar, the quadratic Gaudin Hamiltonian
coincides with the generator $L^{\on{GKO}}_0$ of the coset Virasoro
algebra. Since the entire coset Virasoro algebra commutes with the
diagonal action of $\ghat$, the first requirement on the higher Gaudin
Hamiltonians is satisfied for any element of its completed universal
enveloping algebra. Now we search for a commutative subalgebra of this
enveloping algebra containing $L^{\on{GKO}}_0$.

Consider the case of $\gtil = \wt\sw_2$. Then a natural candidate for
this commutative subalgebra is provided by the algebra of quantum KdV
Hamiltonians. Let us recall that the quantum KdV Hamiltonians are
certain elements $H_{2i+1}, i\geq 0$, of the completed universal
enveloping algebra of the Virasoro algebra whose symbols are the
classical KdV Hamiltonians (see \secref{class ham kdv}). The first of
them, $H_1$, is nothing but $L_0$, and the next one is
$$
H_3 = \frac{1}{2} L_0^2 + \sum_{n > 0} L_{-n} L_n = \frac{1}{2} \int
\Wick T(w)^2 \Wick \, w^3 dw.
$$
The other Hamiltonians have the form
$$
H_{2i+1} = \int \Wick T(w)^{i+1} \Wick \, w^{2i+1} dw + \on{lower } \;
\on{ order } \; \on{ terms}.
$$
Explicit formulas for the first few of them may be found in
\cite{BLZ1}. The existence of these Hamiltonians has been established in
\cite{FF:toda,FF:laws}. We conjecture that

\bigskip

\hspace*{4mm} \parbox{.9\linewidth}{\em
The quantum KdV Hamiltonians of the coset Virasoro algebra are the
higher Gaudin Hamiltonians of the two-point $\su$ affine Gaudin
model with $N=2$.}

\bigskip

More concretely, we conjecture the following:

\begin{conj}
The affine Bethe vectors \eqref{genbv} (for $\su$ and $N=2$) are
simultaneous eigenvectors of all quantum KdV Hamiltonians of the coset
Virasoro algebra, provided that the Bethe Ansatz equations
\eqref{bethe} are satisfied.
\end{conj}

Consider the special case when the module $M_1$ is the vacuum
module $V_{(0,1)}$ of level $k_1=1$ and $M_2$ is an integrable module
$V_{(n,k)}$ of level $k_2 = k \in \Z_+$ and highest weight $n \in
\Z_+, n\leq k$, with respect to the one-dimensional Cartan subalgebra
$\h \subset \htil$. Then, according to \cite{GKO}, we have the
decomposition
$$
V_{(0,1)} \otimes V_{(n,k)} \simeq \bigoplus_{m=0}^{k+1} V_{(m,k+1)}
\otimes L_{n,m},
$$
with respect to the diagonal action of $\su$, where $L_{n,m}$ is the
irreducible module over the Virasoro algebra from the $(k+2,k+3)$
unitary minimal model with the central charge $1-6/(k+2)(k+3)$.

If the above conjecture is true, then we obtain an effective method
for constructing eigenvectors of the quantum KdV Hamiltonians on the
irreducible modules of the Virasoro algebra from the unitary minimal
models.

In the case of a general simply-laced Lie algebra $\g$, we also expect
that the higher affine Gaudin Hamiltonians coincide with the quantum
$\ghat$-KdV Hamiltonians, which generate a commutative subalgebra in
the coset $\W$-algebra. If this is true, then the Gaudin Bethe vectors
would give us eigenvectors of the quantum $\ghat$-KdV Hamiltonians.

In the above discussion we considered the ``local'' Gaudin
Hamiltonians, that is, those that may be expressed as infinite sums of
monomials in the generators of $\gtil$ of {\em bounded} order (for
example, the Hamiltonian $\Xi$ has order $2$). However, in addition to
these local Hamiltonians there are also non-local ones. In order to
explain the appearance of these non-local Hamiltonians, we consider
the quasi-classical analogues of the affine Gaudin models: the affine
Hitchin systems.

\subsection{Affine Hitchin systems}    \label{hit sys}

The Gaudin models associated to simple Lie algebras have
quasi-classical versions, which may be interpreted as the
Beauville--Hitchin systems on $\pone$ with ramification (see
\cite{Hitch,Beauville,Markman,Bot}). In the case of Gaudin models with
regular singularities the algebra of quantum Hamiltonians is a
commutative subalgebra of $U(\g)^{\otimes N}$. Its quasi-classical
counterpart is then a Poisson commutative subalgebra in
$S(\g)^{\otimes N} = \on{Fun}((\g^*)^{\oplus N})$. Let us identify the
Poisson manifold $(\g^*)^{\oplus N}$ (with its Kirillov--Kostant
Poisson structure) with the space of one-forms on $\pone$ with poles
of order one at the points $z_1,\ldots,z_N$, and $\infty$,
\begin{equation}    \label{eta}
L = \sum_{i=1}^N \frac{A_i}{z-z_i} dz, \qquad A_i \in \g^*.
\end{equation}
Let $\ol{P}_1,\ldots,\ol{P}_\ell$ be generators of the algebra
$\on{Inv} \g^*$ of $G$-invariant functions on $\g^*$, of degrees
$d_1+1,\ldots,d_\ell+1$, where the $d_i$'s are the exponents of
$\g$. Then we have the Hitchin map
\begin{align*}
p: (\g^*)^{\oplus N} &\to {\mc H} = \bigoplus_{i=1}^N
H^0\left(\pone,(\Omega(z_1+\ldots+z_N+\infty))^{\otimes
  (d_i+1)}\right), \\
L &\mapsto (P_1(L),\ldots,P_\ell(L)).
\end{align*}

The pull-backs of polynomial functions on ${\mc H}$ to $(\g^*)^{\oplus
N}$ via the map $p$ give rise to Poisson commuting Hamiltonians on
$(\g^*)^{\oplus N}$. They generate a Poisson commutative subalgebra of
$S(\g)^{\otimes N}$. According to the results of \cite{F:faro,FFT},
this is precisely the algebra of symbols of the higher Gaudin
Hamiltonians.

Now we wish to generalize this construction to the affine case, where
we will have commuting quantum Hamiltonians in a completion of
$U(\gtil)^{\otimes N}$. Their symbols should therefore belong to a
completion of the algebra of functions on $(\gtil^*)^{\oplus N}$. It
is natural to identify the latter with the space of one-forms $L$ as
in \eqref{eta}, but with each $A_i$ now being an element of
$\gtil^*$. Actually, as in the case of the shift of argument
subalgebra, we will consider instead of the $i$th copy of $\gtil^*$,
the hyperplane $\ghat^*_{k_i}$ in $\ghat^*$ which consists of the
linear functionals on $\ghat$ taking value $k_i$ on the central
element ${\mb 1}$. We will identify this hyperplane with the space of
$k_i$-connections on the punctured disc
$$
k_i \pa_t + A_i(t), \qquad A_i(t) \in \g\ppart \overset{\ka_0}\simeq
\g^*\ppart,
$$
so that the coadjoint action of $G\ppart$ becomes the gauge action on
the $k_i$-connections (compare with \secref{aff rs}). This hyperplane
carries a canonical Kirillov--Kostant Poisson structure.

Then the $L$-operator of our affine Hitchin system may be represented
in the form
\begin{equation}    \label{L}
L = \sum_{i=1}^N \frac{k_i \pa_t + A_i(t)}{z-z_i} dz, \qquad A_i \in
\g\ppart.
\end{equation}
In what follows we will assume that $k_i \neq 0$.

Equivalently, we may rewrite this as a connection on the punctured
disc, which depends on the spectral parameter $z \in \pone \bs \{
z_1,\ldots,z_N,\infty \}$:
\begin{equation}    \label{conn x}
L' = \pa_t + \frac{1}{f(z)} \sum_{i=1}^N \frac{A_i(t)}{z-z_i} dz,
\end{equation}
where
\begin{equation}    \label{fx}
f(z) = \sum_{i=1}^N \frac{k_i}{z-z_i}.
\end{equation}

Now we will define commuting Hamiltonians of the affine Hitchin system
by evaluating on $L'$ gauge invariant functions on the space of
connections $\pa_t + A(t)$.

\subsection{Commuting Hamiltonians from the monodromy matrix}

As in the case of the shift of argument model (see \secref{aff rs}),
the invariant functions are obtained by taking the invariants of the
monodromy of the connection with respect to the adjoint action. Denote
by $M(z) \in G$ the formal monodromy matrix of the connection operator
$L'$ given by formula \eqref{conn x}. As before, $M(z)$ is defined as
follows: its action on any finite-dimensional representation $V$ of
$G$ is obtained by computing the monodromy of the corresponding matrix
differential operator $\pa_t + A(t)$ acting on $V$.

For any function $\varphi$ on $G$ which is invariant under the adjoint
action of $G$, the Hamiltonian
$$
H_{\varphi}(z) = \varphi(M(z))
$$
is a gauge invariant function of $L'$ (or of $L$ given by formula
\eqref{L}), depending on the spectral parameter $z$. One shows in
the same way as in the finite-dimensional case that
\begin{equation}    \label{comm ham1}
\{ H_{\varphi}(z),H_{\psi}(y) \} = 0, \qquad \forall
\varphi,\psi \in \on{Fun}(G)^{G}, \; z,y \in \C \bs \{
z_1,\ldots,z_N \}.
\end{equation}
Therefore we may obtain families of commuting Hamiltonians of the
affine Hitchin model by expanding the Hamiltonians $H_{\varphi}(z)$
around points of $\pone$.

Let us first consider the expansion around one of the points $z_i,
i=1,\ldots,N$. The expansion of the operator $k_i L'$ at $z_i$ reads
as follows:
$$
k_i \pa_t + A_i(t) + (z-z_i) \left( - A_i(t) \sum_{j \neq
i} \frac{k_j}{k_i(z_i-z_j)} + \sum_{j \neq i}
\frac{A_j(t)}{z_i-z_j} \right) + \ldots,
$$
where the dots stand for higher order terms in $(z-z_i)$.

We now apply \lemref{linear term} in this case, using $(z-z_i)$ as the
parameter $\ep$. Let us write
$$
A_i(t) = \sum_a J_{i,a}(t) J^a,
$$
and denote by
$$
\Psi_i(t) = \sum_a \Psi_{\varphi,i}^a(t) J^a
$$
the single-valued solutions of the equation
\begin{equation}    \label{eq adj}
k_i \pa_t \Psi_i(t) + [A_i(t),\Psi_{i}(t)] = 0,
\end{equation}
corresponding to $\varphi \in (\on{Fun} G)^G$, as before. Then we
find, using \lemref{linear term}, that the $(z-z_i)$-linear term in
the expansion of $\varphi(M(z))$, where $M(z)$ is the monodromy of the
operator $k_i L'$, is equal to
\begin{equation}    \label{Hia}
H_{i,\varphi} = \int \left( - J_{i,a}(t) \sum_{j \neq
i} \frac{k_j}{k_i(z_i-z_j)} + \sum_{j \neq i}
\frac{J_{j,a}(t)}{z_i-z_j} \right) \Psi_{i,\varphi}^a(t) dt
\end{equation}
(here and below the summation over the repeated index $a$ is
understood). The commutativity \eqref{comm ham1} of the Hamiltonians
$H_{\varphi}$ then implies that
$$
\{ H_{i,\varphi},H_{k,\psi} \} = 0, \qquad i=1,\ldots,N, \quad \varphi,
\psi \in (\on{Fun} G)^G.
$$

Higher terms in the expansion of $\varphi(M(z))$ in $(z-z_i)$ are
given by more complicated formulas, involving multiple integrals of
solutions of \eqref{eq adj}. All of these Hamiltonians are non-local,
in the sense that they cannot be realized as integrals of differential
polynomials in the $A_i(t)$.

In order to obtain local commuting Hamiltonians, we need to consider
the expansion of the operator $L'$ at the points
$z=y_1,\ldots,y_{n-1}$, which satisfy $f(y_i) = 0$, where $f(z)$ is
given by formula \eqref{fx} (this is similar to what happens in the
more familiar integrable systems; see, e.g., \cite{DS,Reyman})). We
expect that for each of these points there is an infinite series of
Hamiltonians whose degrees are all positive integers equal to the
exponents of $\gtil$ modulo the Coxeter number. These Hamiltonians
should define interesting hierarchies of commuting flows on our phase
space $\prod_{i=1}^N \ghat_{k_i}^*$.

We have not seen these particular hierarchies discussed in the
literature, but some closely related hierarchies have been studied
before. For example, consider the special case when $N=2$. Then the
limit of the $L$-operator \eqref{conn x} when $k_2 \to 0$ is closely
related to the $L$-operator of the principal chiral field model. We
expect that the corresponding limits of our Hamiltonians coincide with
the Hamiltonians of that model constructed in \cite{Evans} (in
particular, they have the right degrees). We plan to discuss these
Hamiltonians in more detail elsewhere.

\subsection{Quantum non-local Hamiltonians}    \label{nl Ham}

The classical Hamiltonians of the affine Hitchin system may be
quantized, giving rise to quantum (non-local) Hamiltonians of the
affine Gaudin model.

Let us consider first the case of two points, so that $N=2$. In this
case we have the following classical Hamiltonians:
\begin{align*}
H_{1,\varphi} &= \int \left( - \frac{k_2}{k_1} J_{1,a}(t) + J_{2,a}(t)
\right) \Psi_{1,\varphi}^a(t) dt, \\ H_{2,\varphi} &= \int \left( -
\frac{k_1}{k_2} J_{2,a}(t) + J_{1,a}(t) \right) \Psi_{2,\varphi}^a(t)
dt
\end{align*}
(we have multiplied them by $(z_1-z_2)$). To quantize these
Hamiltonians, we replace $J_{i,a}(t)$ by the current $J^{(i)}_a(z)$
and $\Psi_{i,\varphi}^a(t)$ by the component $\Phi^{a(i)}_p(z)$ of
the primary field $\Phi^{(i)}_p(z)$, introduced in \secref{affine
shift quantum},
$$
\Phi^{(i)}_p(w) = \sum_a \Phi^{a(i)}_p(w) J_a.
$$
Here, as before, the upper index $(i)$ indicates that this operator
acts along the $i$th factor of the tensor product $M_1 \otimes
M_2$. We also need to adjust the levels, as in the equation
\eqref{quantum sol} (and the Segal--Sugawara formula).

The resulting quantum Hamiltonians read
\begin{align*}
H_{1,p} &= - \frac{k_2}{k_1+h^\vee} \sum_{n \in \Z} \Wick J_{a,n}
\Phi^{a}_{p,(-n)} \Wick \otimes 1 + \sum_{n \in \Z}
\Phi^{a}_{p,(-n)} \otimes J_{a,n}, \\ H_{2,p} &= -
\frac{k_1}{k_2+h^\vee} \sum_{n \in \Z} 1 \otimes \Wick J_{a,n}
\Phi^{a}_{p,(-n)} \Wick + \sum_{n \in \Z} J_{a,n} \otimes
\Phi^{a}_{p,(-n)}.
\end{align*}
They act on the tensor product $M_{\la_1,k_1} \otimes M_{\la_2,k_2}$
of Verma modules over $\ghat$ (or more general modules $M_1 \otimes
M_2$, provided that the primary fields $\Phi^{(i)}_p(w), i=1,2$,
are well-defined on them).

It is easy to see that they commute with the diagonal action of
$\ghat$. Indeed, consider again the extended affine algebra $\gtil =
\C d \ltimes \ghat$. We choose dual bases $\{ J^{\wt{a}} \}$ and $\{
J_{\wt{a}} \}$ in $\gtil$ as in \secref{gen aff}. Recall that they
consist of elements of the form $J^{a,n}$ (resp., $J_{a,n}$), ${\mb
1}$ and $d$ (resp., $d$ and ${\mb 1}$). Suppose that we have two
homomorphisms of $\gtil$-modules $\rho_i: \gtil \to \on{End}_\C M_i,
i=1,2$. Then the endomorphism
\begin{equation}    \label{end diag}
\sum_{\wt{a}} \rho_1(J_{\wt{a}}) \otimes \rho_2(J^{\wt{a}})
\end{equation}
of $M_1 \otimes M_2$ commutes with the diagonal action of $\gtil$.

In our case, we choose as the first homomorphism $\rho_1$ the action
of $\gtil$ on $M_1$. Under this action ${\mb 1}$ maps to $k_1
\on{Id}$ and $d$ maps to the Sugawara operator $-L_0$. We define the
second homomorphism as follows:
\begin{align*}
J^{a,n} &\mapsto \Phi^{a}_{p,(n)}, \\
d &\mapsto - \frac{1}{k_2+h^\vee} \sum_{n \in \Z} \Wick J_{a,n}
\Phi^{a}_{p,(-n)} \Wick \, , \\
{\mb 1} &\mapsto 0.
\end{align*}
One checks, in the same way as for the Segal--Sugawara formula, that
this indeed defines a homomorphism of $\gtil$-modules. The resulting
operator \eqref{end diag} is the Hamiltonian $H_{2,p}$. Hence it
commutes with the diagonal action of $\gtil$. Since the Hamiltonian
$H_{1,p}$ is obtained from $H_{2,p}$ by switching the modules
$M_1$ and $M_2$, we find that it also commutes with $\gtil$.

Furthermore, we claim that the Hamiltonians $H_{i,p}$ commute with
each other for all $i=1,2$ and $p = 1,\ldots,\ell$. This may be
checked by a direct calculation which will be presented in
\cite{FF:next}.

\medskip

The above formulas are generalized in a straightforward way to the
case of $N$ points. We now have the following commuting Hamiltonians
quantizing the classical Hamiltonians \eqref{Hia}:
\begin{equation}    \label{N pts}
H_{i,p} = \sum_{n \in \Z} \Wick \left( -
\frac{J_{a,n}^{(i)}}{k_i+h^\vee} \sum_{j \neq i} \frac{k_j}{z_i-z_j} +
\sum_{j \neq i} \frac{J_{a,n}^{(j)}}{z_i-z_j} \right)
\Phi^{a(i)}_{p,(-n)} \Wick \, ,
\end{equation}
where $i=1,\ldots,N, p = 1,\ldots,\ell$. They act on the tensor
product $\bigotimes_{i=1}^N M_{\la_i,k_i}$ of Verma modules, or more
general highest weight $\ghat$-modules of levels $k_i, i=1,\ldots,N$
(recall that $k_i \neq -h^\vee$ by our assumption). These are the
simplest {\em non-local quantum affine Gaudin Hamiltonians}.

Other non-local quantum Gaudin Hamiltonians may be constructed by
generalizing the procedure of \cite{BLZ1,BHK,FFS}, as outlined in
\secref{affine shift quantum}. We will discuss this in more detail in
\cite{FF:next}.

\subsection{Hamiltonians in the limit $z_2 \to \infty$}

It is important to observe that the Hamiltonians of the shift of
argument Gaudin model may be obtained from the Hamiltonians of the
Gaudin model with regular singularities and $N=2$ in a certain limit
as $z_2 \to \infty$ (we also set $z_1=0$). In the finite-dimensional
case this was explained by L. Rybnikov in \cite{Ryb}. The Hamiltonians
of the latter model generate a commutative subalgebra in $U(\g)
\otimes U(\g)$. Consider the limit $z_2 \to \infty$, in which we also
degenerate the second factor $U(\g)$ (corresponding to the point
$z_2$) to $S(\g)$. Then we obtain a commutative subalgebra in $U(\g)
\otimes S(\g)$. Applying the homomorphism $S(\g) \to \C$,
corresponding to evaluation at $\chi \in \g^*$, we then obtain a
commutative subalgebra in $U(\g)$. It turns out to be the quantum
shift of argument subalgebra ${\mc A}_\chi$, as shown in \cite{Ryb}.

We may apply the same limiting procedure in the affine case. The above
quantum Hamiltonians of the Gaudin model with regular singularities
and $N=2$ generate a commutative subalgebra in the (completed) tensor
product $U(\gtil) \wh\otimes U(\gtil)$. Taking the limit $z_2 \to
\infty$ and at the same time degenerating $U(\gtil)$, we obtain a
commutative subalgebra in $U(\gtil) \wh\otimes S(\gtil)$. We then need
to evaluate our Hamiltonians at $\chi \in \gtil^*$ along the second
factor. For example, let us apply this procedure to the Hamiltonians
$H_{1,p}$ and $H_{2,p}$ constructed above. Let us set $k_2=0$. Then in
our limit the generators $J_{a,n}$ of $\gtil$ becomes classical, that
is, linear functionals on $\gtil^*$. Evaluating them on $\chi \in \h^*
\subset \gtil^*$, we obtain precisely the non-local Hamiltonians
$$
H_p(\chi) = \sum_a \Phi^a_{p,(0)} \langle \chi,J_{a,0} \rangle
$$
introduced in \secref{affine shift quantum}.

On the other hand, if we apply this procedure to $H_{2,p}$, then the
primary fields $\Phi^a_{p}(w)$ will become classical, i.e., will
degenerate to $\Psi^a_\varphi(w)$, and we will need to evaluate them
on $\chi$. The resulting Hamiltonians will be linear combinations of
$J_{a,0}$ that belong to $\h \subset \gtil$, where $\h$ is the
centralizer of $\chi$ in $\g$. Thus, we obtain local (in fact, linear)
Hamiltonians of the shift of argument affine Gaudin model.

Similar formulas will appear in general. The non-local Hamiltonians in
the Gaudin model with regular singularities and $N=2$ have, roughly,
the following form:
$$
\sum_i A_i \otimes B_i \qquad \on{or} \qquad \sum_i B_i \otimes A_i,
$$
where the $A_i$ are non-local expressions which are quantizations of
the derivatives of the ``quantum monodromy matrix'', such as
$\Phi^a_{p,(-n)}$, and the $B_i$ are local expressions, that is, of
finite degree in $J^a_n$. When we degenerate $U(\gtil)$ to
$S(\gtil)$ along the second factor, either $B_i$ or $A_i$ become
classical, and we may then evaluate them on $\chi \in \gtil^*$. In the
first case we obtain non-local Hamiltonians of the shift of argument
affine Gaudin model (such as the $\g$-AKNS model). In the second case,
we obtain local Hamiltonians.

A similar limiting procedure may be constructed for the Hamiltonians
of the quantum $\ghat$-KdV systems. This suggests a possible way to
construct local quantum Hamiltonians of the $\g$-AKNS and $\ghat$-KdV
systems.

\subsection{Spectra of the quantum Hamiltonians and affine opers}

In \cite{F:flag} to each solution of the Bethe Ansatz \eqref{bethe}
equations for the affine Lie algebra $\ghat$ the following connection
on the $^L \ol{H}$-bundle $\Omega^\rho$ was attached:
\begin{equation}    \label{new miura}
\ol\nabla = \pa_z - \sum_{i=1}^N \frac{\nu_i}{z-z_i} +
\sum_{j=1}^m \frac{\al_{i_j}}{z-w_j}.
\end{equation}
It was shown in \cite{F:flag} that the Bethe Ansatz equations are
equivalent to the condition that the $^L \ghat$-oper $\nabla =
\ol\nabla + p_{-1}$ obtained by applying the Miura transformation to
the connection \eqref{new miura} is regular at the points $w_j,
j=1,\ldots,m$, or, equivalently, has no monodromy around these
points. This oper has regular singularities at the points $z_i,
i=1,\ldots,N$ and $\infty$ with the residues determined by the highest
weights $\nu_i, i=1,\ldots,N$.

This suggests that the true parameters of the Bethe vectors are not
the Cartan connections \eqref{new miura}, but the corresponding $^L
\ghat$-opers obtained by applying the Miura transformation. Thus, we
propose the following answer to question (2) of the previous section:
what parameterizes common eigenvalues of the affine Gaudin
Hamiltonians?

\bigskip

\hspace*{4mm} \parbox{.9\linewidth}{{\em The common eigenvalues of the
affine Gaudin Hamiltonians on the tensor product $\bigotimes_{i=1}^N
M_{\nu_i}$ of Verma modules over $\ghat$ are encoded by $^L
\ghat$-opers on $\pone$ with regular singularities at the points
$z_1,\ldots,z_N$ and $\infty$.}}

\bigskip

Again, we emphasize that this $\pone$ is a curve in the {\em spectral
parameter} $z$ of the affine Gaudin model.

\medskip

To see what these $^L \ghat$-opers look like, consider the case when
$\ghat = {}^L \ghat = \su$. We will use the conventions of
\secref{spectra aff}. We again divide the set of points
$w_1,\ldots,w_m$ into two subsets $w^1_1,\ldots,w^1_{m_1}$, and
$w^0_1,\ldots,w^0_{m_0}$, corresponding to the simple roots $\al_1$
and $\al_0$, respectively. The connection \eqref{new miura} reads
$$
\pa_z + \begin{pmatrix} u(z) & \la \\ 1 & -u(z) \end{pmatrix} +
\sum_{i=1}^N \frac{k_i}{z-z_i} d,
$$
where
\begin{equation}    \label{ut}
u(z) = - \sum_{i=1}^N \frac{\ell_i}{z-z_i} +
\sum_{j=1}^{m_1} \frac{1}{z-w^1_j} - \sum_{j=1}^{m_0}
\frac{1}{z-w^0_j}.
\end{equation}
Applying the gauge transformation by
$$
\begin{pmatrix} 1 & -u(z) \\ 0 & 1 \end{pmatrix},
$$
we obtain the operator
\begin{equation}    \label{oper bethe}
\pa_z + \begin{pmatrix} 0 & v(z) + \la \\ 1 & 0 \end{pmatrix} +
\sum_{i=1}^N \frac{k_i}{z-z_i} d,
\end{equation}
where
\begin{align} \notag
v(z) &= u(z)^2 + \pa_z u(z) \\ \label{form of conn} &= \sum_{i=1}^N
\frac{\ell_i(\ell_i+1)}{(z-z_i)^2} + \sum_{i=1}^N \frac{c_i}{z-z_i} +
\sum_{j=1}^{m_0} \frac{2}{(z-w^0_j)^2} + \sum_{j=1}^{m_0}
\frac{v_{j,-1}}{z-w^0_j}, \\ \notag c_i &= 2\ell_i \left( \sum_{j \neq
i} \frac{\ell_j}{z_i-z_j} - \sum_{j=1}^{m_1} \frac{1}{z_i-w^1_j} +
\sum_{j=1}^{m_0} \frac{1}{z_i-w^0_j} \right), \\
\label{gaj} v_{j,-1} &= \sum_{i=1}^N \frac{k_i}{w_j^0-z_i}.
\end{align}
These formulas follow from the Bethe Ansatz equations
\eqref{bethe}, which read
\begin{align}    \label{beth1}
& \sum_{i=1}^N \frac{\ell_i}{w_j^1-z_i} -
\sum_{s \neq j} \frac{1}{w_j^1-w^1_s} + \sum_{s=1}^{m_0}
\frac{1}{w^1_j-w^0_s} = 0, \qquad j=1,\ldots,m_1, \\ \label{beth2}
&\sum_{i=1}^N \frac{\frac{k_i}{2} - \ell_i}{w_j^0-z_i} +
\sum_{s=1}^{m_1} \frac{1}{w_j^0-w^1_s} - \sum_{s \neq j}
\frac{1}{w^0_j-w^0_s} = 0, \qquad j=1,\ldots,m_0.
\end{align}
These equations imply, in particular, that $v(z)$ has no singularities
at the points $w^1_j, j=1,\ldots,m_1$.

In addition, because $v(z)$ is obtained as the Miura transformation of
$u(z)$ given by formula \eqref{ut}, the coefficients of its expansion
at each point $w_j^0, j=1,\ldots,m_0$,
$$
v(z) \sim \frac{2}{(z-w^0_j)^2} + \sum_{n \geq -1} v_{j,n} (z-w^0_j)^n
$$
satisfy the equation \eqref{third order}. This is the condition that
the operator \eqref{oper bethe} has no monodromy around the point
$w^0_j$. For instance, if all $k_i=0$, then $v_{j,-1}=0$ and the
equation \eqref{third order} reduces to $v_{j,1}=0$ for all
$j=1,\ldots,m_0$. Hence we obtain that in this case the condition is
that the coefficients in front of $(z-w^0_j)^{-1}$ and $(z-w^0_j)$ in
the expansion of $v(z)$ in Laurent series in $(z-w^0_j)$ should vanish
for all $j=1,\ldots,m_0$.

Thus, a connection $\pa_z + u(z)$ given by formula \eqref{ut} with the
numbers $w^0_j, w^1_j$ satisfying the Bethe Ansatz equations
\eqref{beth1}--\eqref{beth2} gives rise to an operator \eqref{oper
bethe} such that its coefficients in the expansion at $w^0_j,
j=1,\ldots,m_0$, satisfy the equations \eqref{gaj} and \eqref{third
order}. (However, we do not expect the converse to be true in
general.) Our proposal is that the generic common eigenvalues of the
affine Gaudin Hamiltonians on the tensor product $\bigotimes_{i=1}^N
M_{\nu_i}$ of Verma modules over $\su$ are encoded by $\su$-opers of
this kind. In general, some degenerations of these opers will also
correspond to some eigenvectors, as in the shift of argument Gaudin
model (see \secref{spectra aff}).

It is useful to express the operator \eqref{oper bethe} as a second
order differential operator (compare with \secref{spectra
aff}). Applying the (formal) gauge transformation by the element
$\prod_{i=1}^N (z-z_i)^{k_i d}$, we obtain
$$
\pa_z + \begin{pmatrix} 0 & v(z) + \la \prod_{i=1}^N (z-z_i)^{k_i} \\
  1 & 0 \end{pmatrix},
$$
which may be rewritten as the following second order operator with
spectral parameter
$$
\pa_z^2 - v(z) - \la \prod_{i=1}^N (z-z_i)^{k_i}.
$$

Alternatively, we may rewrite \eqref{oper bethe} as the following
partial differential operator:
$$
\left(\pa_z + \sum_{i=1}^N \frac{k_i}{z-z_i} \la \pa_\la \right)^2 -
v(z) - \la.
$$

\subsection{Integrable modules}

Next, we consider the common eigenvalues of the affine Gaudin
Hamiltonians on the tensor product of integrable
$\ghat$-modules. Recall that for each integral dominant weight $\nu$
we have an integrable module $V_\nu$, which is the irreducible
quotient of the Verma module $M_\nu$. By analogy with the
finite-dimensional case, we propose the following description of the
common eigenvalues of the affine Gaudin Hamiltonians on the tensor
product $\bigotimes_{i=1}^N V_{\nu_i}$:

\begin{conj}    \label{aff reg bethe}
There is a bijection between the set of common eigenvalues of the
affine Gaudin Hamiltonians on the tensor product $\bigotimes_{i=1}^N
V_{\nu_i}$ of integrable modules over $\ghat$ (counted without
multiplicities) and the $^L \ghat$-opers on $\pone$ which are
$\nu_i$-regular at $z_i, i=1,\ldots,N$, and have regular singularity
at $\infty$.
\end{conj}

Note that such opers necessarily have {\em trivial monodromy} around
each of the points $z_i, i=1,\ldots,N$, and hence they have no
monodromy on $\pone$.

If the Bethe Ansatz equations \eqref{bethe} are satisfied and the
projection of the Bethe vector \eqref{genbv} onto $\bigotimes_{i=1}^N
V_{\nu_i}$ is non-zero, then we obtain an eigenvector of the affine
Gaudin Hamiltonians in $\bigotimes_{i=1}^N V_{\nu_i}$. A necessary
condition for this is that its weight
\begin{equation}    \label{ga}
\ga = \sum_{i=1}^N \nu_i - \sum_{j=1}^m \al_{i_j}
\end{equation}
is dominant. The corresponding eigenvalue should then be encoded by
the $^L \ghat$-oper obtained by applying the Miura transformation to
the Cartan connection \eqref{new miura}. As shown in \cite{F:flag},
such an oper satisfies the conditions of \conjref{aff reg bethe}. This
provides the first consistency check for this conjecture.

Note that we allow the possibility that there exist eigenvectors other
than the Bethe vectors in $\bigotimes_{i=1}^N V_{\nu_i}$. In this case
we still believe that the common eigenvalues of the affine Gaudin
Hamiltonians are encoded by the $^L \ghat$-opers as in \conjref{aff
reg bethe}. However, the eigenvalues corresponding to the non-Bethe
eigenvectors should be encoded by opers that cannot be expressed as
the Miura transformation of a connection of the form \eqref{new
miura}, of the kind described at the end of \secref{spectra aff}.

Finally, if the weight \eqref{ga} is not dominant, then the
corresponding Bethe vector in $\bigotimes_{i=1}^N V_{\nu_i}$ is
zero. However, it is still interesting to consider the corresponding
space of solutions of the Bethe Ansatz equations \eqref{bethe}. As
shown in \cite{F:flag}, Theorem 5.7, for each fixed weight $\ga$ this
set may be identified with a disjoint union of open subsets of a
finite-dimensional $^L \ol{B}_-$-orbit (Schubert cell) in the affine
flag variety $^L \ol{G}/{}^L \ol{B}_-$. We will consider a special
case of this correspondence in the next subsection.

\subsection{Special case: rational solutions of the KdV hierarchy}

Consider the special case of the above construction when $\ghat=\su$,
$N=1$, the highest weight $\nu_1=0$, and the level $k_1=0$. According
to formula \eqref{form of conn}, the corresponding generic $\su$-opers
may be represented by the following second order differential
operators:
\begin{equation}    \label{rat kdv}
\pa^2_z - \sum_{j=1}^{m_0} \frac{2}{(z-w_j^0)^2} - \la,
\end{equation}
which satisfy the condition
\begin{equation}    \label{rat kdv cond}
\sum_{s \neq j} \frac{1}{(w^0_s-w^0_j)^3} = 0, \qquad j=1,\ldots,m_0.
\end{equation}
These equations (analogous to the Bethe Ansatz equations
\eqref{bethe}) mean that the $(z-w^0_j)$-coefficient $v_{j,1}$ of the
expansion of $v(z)$ at $z=w^0_j$ vanishes, which is equivalent to
\eqref{third order} in this case. This condition means that the
operator \eqref{rat kdv} has {\em no monodromy} on $\pone$ for all
values of the spectral parameter $\la$.

The operators \eqref{rat kdv} were introduced by
Airault--McKean--Moser in \cite{AMM} as the generic rational solutions
of the KdV hierarchy. They were also encountered by Duistermaat and
Gr\"unbaum \cite{DG} in their study of the bispectral problem for
second order differential operators. Here we have obtained these
operators in a different way; namely, in the framework of the affine
Gaudin model in the special case when $N=1, \nu_1=0$ and
$k_1=0$. Following the argument of \cite{F:flag}, we obtain that the
set of these operators may be identified with an open subset of the
$SL_2[[\la]]$-orbit
$$
SL_2[[\la]] \begin{pmatrix} \la^{m_0} & 0 \\ 0 & \la^{-m_0}
\end{pmatrix} SL_2[[\la]]
$$
in the affine Grassmannian $SL_2\pparl/SL_2[[\la]]$.

The points outside this locus correspond to the $\su$-opers in which
some of the points $w_j^0$ coalesce. They are represented by the
more general second order operators of the form
\begin{equation}    \label{rat kdv1}
\pa^2_z - \sum_{j=1}^{m_0} \frac{\ell_j(\ell_j+1)}{(z-w_j^0)^2} - \la,
\end{equation}
where $\ell_j \in \Z_{>0}$, which satisfy the no-monodromy condition
\begin{equation}    \label{rat kdv cond1}
\sum_{s \neq j} \frac{\ell_s(\ell_s+1)}{(w^0_s-w^0_j)^{2p_j+1}} = 0,
\qquad j=1,\ldots,m_0; p_j=1,\ldots,\ell_j.
\end{equation}
These are precisely the operators constructed in \cite{DG} (see
formulas (3.31) and (3.32)), where it was shown that they correspond
to the most general rational solutions of the KdV hierarchy decaying
at $\infty$. The fact that the set of all operators of this form is in
bijection with the affine Grassmannian $SL_2\pparl/SL_2[[\la]]$ is in
agreement with the well-known description of the polynomial
tau-functions of the KdV hierarchy as points of the affine
Grassmannian.\footnote{An identification of the space of rational
solutions of the KdV hierarchy with the space of affine opers of this
kind may also be obtained from the results of \cite{BF}.}

The rational solutions of the $\ghat$-KdV hierarchies for other affine
Kac--Moody algebras may be interpreted in a similar way.

\section{Concluding remarks}    \label{conclusion}

Let us summarize what we have learned so far. Our point of departure
was the observation that soliton systems, such as the KdV and AKNS
hierarchies, are special cases of the Gaudin models associated to
affine Kac--Moody algebras. Therefore the quantization of these models
should be understood in this framework. In order to better understand
these models, we look at their simplified versions; namely, the Gaudin
models associated to finite-dimensional simple Lie algebras. We then
use this analogy to describe various elements of the affine Gaudin
models, and in particular the quantum soliton systems. In this paper
we have made the first steps in this direction.

Let us first discuss the classical soliton systems. These are
integrable systems defined on the space of connections
$$
\pa_t + A(t), \qquad A(t) \in \g\ppart.
$$
We apply a shift of argument and consider the space of connections of
the form
\begin{equation}    \label{conn with shift}
\pa_t + A(t) + z\chi, \qquad A(t) \in \g\ppart.
\end{equation}
The shift parameter $z$ plays the role of a new variable, and so these
integrable systems are properly understood as integrable systems on
the {\em double loop algebra} in variables $t$ and $z$, as was already
noted in \cite{RS}. The classical commuting Hamiltonians are obtained
from the functions $\varphi(M(z))$, where $M(z)$ is the monodromy
matrix of \eqref{conn with shift} and $\varphi$ is a function on $G$
invariant under conjugation. Expanding these functions in $z$, we
obtain local and non-local Hamiltonians. It is important to realize
that, when viewed as a function on the entire double loop algebra,
$\varphi(M(z))$ is actually central, with respect to the natural
Kirillov--Kostant structure.

How can we quantize this system? It is clear that the quantum
Hamiltonians should be obtained from quantum versions of the functions
$\varphi(M(z))$, which should be in the center of a proper completion
of the enveloping algebra of the double loop algebra. Therefore to
describe the spectrum of these Hamiltonians we need to understand the
structure of this center.

At this point it is instructive to look more closely at the ``baby
version'' of this model; namely, the Gaudin model associated to a
finite-dimensional simple Lie algebra $\g$. In this case the naive
phase space is $\g^*$, but because we make the shift of argument, the
true phase space is the dual of the loop algebra $\g\zpart$. The
classical Hamiltonians appear from the $z$-expansion of the central
functions on $\g\zpart^*$, obtained from invariant polynomials on
$\g^*$. Naively, we expect that the quantum versions of these
functions are central elements in a completed enveloping algebra of
$\g\zpart$. But there is an important twist to this story: we need to
replace $\g\zpart$ by its central extension; its completed enveloping
algebra then has a large center at the critical level.

This center turns out to be isomorphic to the algebra of functions on
the space $\on{Op}_{^L \g}(D^\times)$ of $^L \g$-opers on the
punctured disc, where $^L \g$ is the Langlands dual Lie algebra of
$\g$ \cite{FF:gd,F:book}. This is a non-trivial result, which is
closely related to the geometric Langlands correspondence
\cite{BD,F:rev}. This implies \cite{F:faro,FFT} that the spectrum of
the algebra of quantum Hamiltonians of the Gaudin models associated to
$\g$ is isomorphic to the space of $^L \g$-opers on $\pone$ with
prescribed singularities at finitely many points (depending on the
particular model).

Now, we conjecture that the algebras of quantum Hamiltonians of the
affine Gaudin models exhibit similar features. Namely, their spectra
should be expressed in terms of $^L \ghat$-opers on $\pone$ with
prescribed singularities at finitely many points, where $^L \ghat$ is
the Langlands dual Lie algebra of $\ghat$. These $^L \ghat$-opers may
be described more concretely as differential operators. For example,
for $\ghat = \su$ these are Schr\"odinger operators with spectral
parameter.

We remark that in this picture opers appear in two different contexts,
and this may cause some confusion. First of all, the classical
$\ghat$-KdV system has the space of $\ghat$-opers (or $\g$-opers) on
the punctured disc as the phase space. Then we quantize it and look at
the spectra of the corresponding algebra of quantum KdV
Hamiltonians. We conjecture that those are encoded by $^L \ghat$-opers
on $\pone$. So opers appear again, but they have nothing to do with
the $\ghat$-opers entering the definition of the classical KdV system!
For one thing, they correspond to an {\em a priori} different Lie
algebra; namely, the Langlands dual Lie algebra $^L \ghat$. Besides,
they are defined on $\pone$ (on which the spectral parameter $z$ is a
global coordinate), not on a punctured disc.

In this paper we have presented explicit conjectures for the $^L
\ghat$-opers encoding the spectra of the quantum Hamiltonians of the
affine Gaudin models associated to $\ghat$. We have checked that our
proposal coincides with the proposal of \cite{BLZ} in the case of the
quantum KdV Hamiltonians.

What can we learn from these conjectures? First of all, placing the
quantum soliton systems in the context of affine Gaudin models is
useful, because we may apply the tools of the Gaudin models to
investigate the soliton systems. An example of such a tool is Bethe
Ansatz, which we use here to produce explicit formulas for the
eigenvectors of the quantum Hamiltonians. Second, this points us to a
conjectural description of the spectra of the quantum Hamiltonians in
terms of opers, which is consistent with and generalizes the proposals
obtained by other methods. Finally, we place the quantum soliton
systems in the context of Langlands duality. We hope that this will
allow us to understand better the true meaning of quantum soliton
system and their relations to other areas of mathematics and physics.

\section{Appendix. Proof of \lemref{comm DMT}}

For any simple Lie algebra $\g$ we have the Segal--Sugawara current
$$
S(w) = \frac{1}{2} \sum_{\al \in \Delta_+} \left( \Wick
e_\alpha(w) f_\alpha(w) \Wick + \Wick f_\alpha(w) e_\alpha(w) \Wick
\right) + \frac{1}{2} \sum_{i=1}^\ell \Wick h_i(w) h^i(w) \Wick \, ,
$$
where $\{ h^i \}$ is the dual basis to the coroot basis $\{ h_i \}$ of
$\h$. Let us set
$$
\wt{S}(w) = \sum_{\al \in \Delta_+} \Wick e_\alpha(w) f_\alpha(w)
\Wick \,
$$
and
$$
\ol{S}(w) =  \frac{1}{2} \sum_{i=1}^\ell \Wick h_i(w) h^i(w) \Wick
\, .
$$

Now denote by
$$
S_\al(w) = \sum_{n \in \Z} S_{\al,n} w^{-n-2}
$$
the element $S(w)$ associated to the $\sw_2$ subalgebra $\g_\al
\subset \g$ generated by $e_\al$ and $f_\al$, and by
$$
S_{\al\beta}(w) = \sum_{n \in \Z} S_{\al\beta,n} w^{-n-2}
$$
the element $S(w)$ associated to the rank two Lie subalgebra
$\g_{\al\beta} \subset \g$ whose root system $\Delta_{\al\beta}
\subset \Delta$ is spanned by the roots $\al,\beta \in \De_+$. We will
use a similar notation for the corresponding currents $\wt{S}(w)$ and
$\ol{S}(w)$.

Adapting the argument proving the commutativity of the DMT
Hamiltonians in the finite-dimensional case (see \cite{TL}), we reduce
the statement of part (1) of \lemref{comm DMT} to checking the
following identity:
\begin{equation}    \label{ident1}
[\wt{S}_{\gamma,0},\wt{S}_{\al\beta,0}] = 0, \qquad \forall \gamma \in
\Delta^+_{\al\beta}.
\end{equation}

We have
\begin{align*}
\wt{S}_{\gamma,0} &= S_{\gamma,0} - \ol{S}_{\gamma,0} - \ga_0, \\
\wt{S}_{\al\beta,0} &= S_{\al\beta,0} - \ol{S}_{\al\beta,0} -
2\rho_{\al\beta,0}.
\end{align*}
Here $\gamma$ (resp., $\rho$) denotes the element of $\h$
corresponding to the root $\ga$ (resp., the half-sum of positive roots
of $\g_{\al\beta}$) in $\h^*$ under the identification $\h \simeq
\h^*$ corresponding to $\ka_0$. We write $\ga_0$ (resp., $\rho_0$) for
$\ga \otimes 1$ (resp., $\rho \otimes 1$) in $\ghat$.

Since all summands in the right hand sides of these formulas commute
with $\h \otimes 1 \subset \ghat$, and so in particular with
$\rho_{\al\beta,0}$, we obtain that \eqref{ident1} is equivalent to
$$
[\wt{S}_{\gamma,0},S_{\al\beta,0} - \ol{S}_{\al\beta,0}] = 0.
$$
Actually, we will show that
$$
[\wt{S}_{\gamma,0},S_{\al\beta,0}] =
[\wt{S}_{\gamma,0},\ol{S}_{\al\beta,0}] = 0.
$$
By using the formalism of operator product expansion, it is easy to
derive these identities from the following identities:
\begin{align}    \label{ident}
& S_{\al\beta,-1} \wt{S}_{\gamma,-2} v_k \in \on{Im} T, \\
\label{ident2}
& \ol{S}_{\al\beta,-1} \wt{S}_{\gamma,-2}  v_k \in \on{Im} T
\end{align}
in the vertex algebra $V_k(\g)$, where $v_k$ is the vacuum vector and
$T$ is the translation operator, in the notation of \cite{FB}. (The
sought-after commutation relations will follow after we make the
change of variables $w=e^u-1$ in the corresponding operator product
expansion.)

To check \eqref{ident1}, recall that the Segal--Sugawara operator
$S_{\al\beta,-1}$ acts on the completed enveloping algebra of
$\ghat_{\al\beta}$ as a multiple of the derivation $T = -t \pa_t$.
This implies formula \eqref{ident}. To prove formula \eqref{ident2},
note that $\wt{S}_{\gamma,-2} v_k = e_{\ga,-1} f_{\ga,-1}
v_k$ and $\ol{S}_{\al\beta,-1}$ is a linear combination of
\begin{equation}    \label{ga and mu}
\sum_{n \in \Z} \check\ga_{-1-n} \check\ga_{n} \quad \on{and} \quad
\sum_{n \in \Z} \check\mu_{-1-n} \check\mu_{n},
\end{equation}
where $\check\mu \in \h$ is orthogonal to $\check\ga$. Therefore the
second term, applied to $\wt{S}_{\gamma,-2} v_k$, gives $0$, and the
first term gives
$$
\sum_{n \in \Z} \check\ga_{-1-n} \check\ga_{n} \cdot e_{\ga,-1}
f_{\ga,-1} v_k = 2 \check\ga_{-1} \check\ga_{-2} v_k = T \cdot
(\check\ga_{-1})^2 v_k.
$$
This implies formula \eqref{ident2} and hence proves part (1) of
\lemref{comm DMT}.

The commutativity of the operators $\wt{T}_\ga$ given by formula
\eqref{wt T} follows from the corresponding analogue of the identity
\eqref{ident1}, which is proved in the same way as above. In order to
complete the proof of part (2) of \lemref{comm DMT}, it remains to
show that the operators $\wt{T}_\ga$ commute with $\wh\h$. This, in
turn, follows from the formula
$$
h_{i,n} \cdot \left( e_{\alpha,-1}
f_{\alpha,-1} + f_{\alpha,-1} e_{\alpha,-1} - \frac{1}{k} h_{\al,-1}^2
\right) v_k = 0, \qquad n \geq 0,
$$
which is verified in a straightforward fashion.

\end{document}